\documentclass{amsart}

\usepackage{amssymb,amsmath,amsfonts,amsthm}
\usepackage{graphicx}

\newtheorem{thm}{Theorem}[section]
\newtheorem{prop}[thm]{Proposition}
\newtheorem{cor}[thm]{Corollary}
\newtheorem{lemma}[thm]{Lemma}
\newtheorem{defn}[thm]{Definition}
\newtheorem{conj}[thm]{Conjecture}
\newtheorem{cond}[thm]{Condition}
\theoremstyle{definition}
\newtheorem{ex}[thm]{Example}
\newtheorem{rem}[thm]{Remark}

\def\A{\mathcal{A}}
\def\a{\alpha}
\def\B{\mathcal{B}}
\def\b{\beta}
\def\C{\mathbb{C}}
\def\centra{\mathcal{C}_m(n)}
\def\ca{\check{\a}}

\def\ceil#1{\lceil{#1}\rceil}
\def\cent{\mathcal{Z}}

\def\e{\epsilon}
\def\F{\mathbb F}
\def\flip{\mathcal{F}}
\def\floor#1{\lfloor{#1}\rfloor}
\def\g2#1#2{\setlength{\unitlength}{1cm}
\begin{picture}(3,1)
\put(1,0.1){\circle{0.2}}
\put(1,0){\line(1,0){1.2}}
\put(1.1,0.1){\line(1,0){1}}
\put(1,0.2){\line(1,0){1.2}}
\put(2.2,0.1){\circle{0.2}}
\put(1.6,0){$\rangle$}
\put(0.9,0.3){#1}
\put(2.1,0.3){#2}
\end{picture} }
\def\G{\mathcal{G}}
\def\ga{\gamma}
\def\H{\mathcal{H}}
\def\half{\frac{1}{2}}
\def\halfZ{\frac{1}{2}\Z}
\def\Hrcc{\hat{\H}^{\tiny{\mbox{temp}}}_{\tiny{\mbox{rcc}}}}

\def\Hgr{\mathbb{H}}
\def\Hgrrcc{\hat{\Hgr}^{\tiny{\mbox{temp}}}_{\tiny{\mbox{rcc}}}}
\def\I{\mathcal{I}}
\def\J{\mathcal{J}}
\def\L{\mathcal{L}}

\def\l{\lambda}

\def\N{\mathbb N}

\def\P{\mathcal{P}}

\def\pf{\noindent{\it Proof}}
\def\qed{\hfill $\square$\par\medskip}
\def\Q{\mathbb{Q}}

\def\real{\mathfrak{a}}
\def\R{\mathbb{R}}
\def\Rdatum{\mathcal{R}}

\def\S{\mathcal{S}}

\def\t{\mathfrak{t}}
\def\trm{\mbox{tr}_m-\mbox{Ind}}
\def\U{\mathcal{U}}
\def\vet#1{\mbox{\boldmath$#1$}}

\def\Z{\mathbb{Z}}

\begin{document}
\title{Generalized Green functions and graded Hecke algebras}
\author{K. Slooten}
\date{\today}
\address{Institut Gaspard Monge, Universit\'e de Marne-la-Vall\'ee, 5 Boulevard Descartes, Champs-sur-Marne, 77454 Marne-la-Vall\'ee Cedex 2,  France}
\email{slooten@univ-mlv.fr}
\begin{abstract} We state a conjecture which gives a combinatorial parametrization of the irreducible tempered representations with real central character of a
graded Hecke algebra with unequal labels, associated to a root sytem of type $B$ or $C$. This
conjecture is based on a combinatorial generalization of the
Springer correspondence in the classical (equal label) case. In
particular, the described modules turn out to have a natural
grading for the action of $W_0$, and are completely determined by
their central character together with the $W_0$-representation in the top
degree. This latter is an irreducible $W_0$-character which we call Springer correspondent.
\end{abstract} \maketitle \begin{small}\tableofcontents\end{small}
\section{Introduction}

In this article we study combinatorics which conjecturally describe the irreducible tempered representations with real central character, for (specialized) graded Hecke algebras of type $B$ and $C$. The graded Hecke algebra was introduced by Lusztig in \cite{lusztig}, where he showed it to play an important role in the representation theory of affine Hecke algebras.

Let ${\mathcal R}=(R_0,X,\check{R}_0,Y,\Pi_0)$ be a root datum. Let the associated affine Weyl group be $W$, whose set of simple reflections we denote by $S$. Let $l$ be the length function on $W$. Given a number ${\bf q}>0$ and for every $s \in S$ a number $f_s\in \R$ such that $f_s=f_{s'}$ if $s \sim_W s'$, put $q(s)={\bf q}^{f_s}$. To these data, one can associate a unital associateve $\C$ algebra $\H=\H(W,q)$ called the affine Hecke algebra. It has a basis $T_w;w\in W$, subject to the relations $T_wT_{w'}=T_{ww'}$ if $l(w)+l(w')=l(ww')$ and $(T_s+1)(T_s-q(s))=0$ for all $s \in S$. It is well known that such an object exists and is unique.

These algebras arise naturally in the representation theory of a reductive $p$-adic group $\G$. In this setting, the best known case is when $f_s=1$ for all $s$ and ${\bf q}=q$, the cardinality of the residue class field of the $p$-adic field. Then the irreducible representations of $\H$ are in bijection with those of $\G$ having an Iwahori-fixed vector. 
Morover, there exists a Plancherel measure preserving bijection between the Iwahori-spherical part of $\hat{\G}$ and $\hat{\H}$ (cf. \cite{HO2}).


More generally, affine Hecke algebras with non-equal labels $q(s)$ arise as well, if one considers unipotent rather than Iwahori-spherical representations of $\G$ (see \cite{lusuni1}). In this case as well, one is lead to the question for a parametrization of the irreducible tempered representations, as well as for the Plancherel measure.

In the equal label case where $\G$ is of adjoint type, a parametrization of the irreducible representations of $\H$ was obtained by Kazhdan and Lusztig, including criteria to select the tempered ones. For affine Hecke algebras with arbitrary labels $q(s)$, Opdam has determined the Plancherel formula (i.e., the measure on $\hat{\H}$ giving the decomposition of the natural trace on $\H$) in \cite{O1}, including the classification of the central characters in the support of the Plancherel measure. It remains an open question, however, to determine the set $\Delta(r)$ of irreducible tempered representations of $\H$ which have central character $r$.

A first step is to consider real central character. In this case, by results of Lusztig, one can solve the problem of finding the irreducible tempered representations of $\H$ by solving the corresponding problem for a related (correctly specialized) graded Hecke algebra associated to the same root system. The graded Hecke algebra contains parameters $k_\a$, for simple $\a \in R_0$, depending on the parameters $q(s)$(see \eqref{kq} below). Let $\Hgr$ denote a graded Hecke algebra associated to the root system $R_0$ and parameters $k_\a$. That is, if $Q$ is the root lattice we put $V=Q\otimes_\Z\R$, and define $\Hgr=\C[W_0] \tilde{\otimes} S(V_\C)$ where the cross relations are given by $xs_\a-s_\a s_\a(x)=k_\a\langle x,\ca\rangle$ for $x \in V_\C$, and $\a$ simple. If all $k_\a=k$, an important role in the classification of the irreducible tempered representations of $\Hgr$ is played by Springer modules. Indeed, then the set $\Hgrrcc$ of equivalence classes of irreducible tempered $\Hgr$- representations with real central character is given by a set of representing modules $\{M_\chi\mid \chi \in \hat{W}_0\}$, such that 
\begin{itemize}
\item[(1)]{The module $M_\chi$ is a naturally graded $W_0$-module whose top degree is $\chi$;}
\item[(2)]{The bijection $\chi \mapsto M_\chi$ is completely determined by requiring that $\chi$ occurs in $M_\chi\mid_{W_0}$;}
\item[(3)]{$M_\chi$ and $M_{\chi'}$ have the same central character if and only if $\chi$ and $\chi'$ are Springer correspondents of the same unipotent class of $\hat{G}$, the complex simply connected reductive group with root system $R_0$;} 
\item[(4)]{Let $M_\chi^i$ be the degree-$i$ part of $M_\chi\mid W_0$. The mutiplicity $(M^i_\chi,\psi)$ is given by (a coefficient of) a Green function;}
\end{itemize}

In this article, we give a conjecture for the special case where the root system $R_0$ is of type $B_n$ (or, equivalently, $C_n$). Suppose that $k_1$ is the label of the long roots and $k_2$ the label of the short roots. We will assume that $k_1 \neq 0$ and $k_2=mk_1$ for $m \in \half \Z$. The conjecture basically states that (1)-(4) still hold, although of course in (3) the set of unipotent class must be replaced by a combinatorial analogue.

The structure of this article is as follows. First, in section \ref{hecke} we explain in more detail the relation between the affine and the graded Hecke algebra, and the determination of the modules $M_\chi$. Although this article focuses on the graded Hecke algebra, we need the determination of the set of central characters of irreducible tempered representations of the affine Hecke algebra, to be able to state the anologous description for the graded Hecke algebra. In section \ref{greenfunctions} we give a combinatorial analogue for arbitrary $k_2=mk_1$ of Lusztig's symbols, which lead to an ordering on $\hat{W}_0$, and a definition of Green functions (as they have been introduced by Shoji, in the context of complex reflection groups). In section \ref{conjecture} we state a conjecture on how these Green functions compute the $W_0$-structure of the modules in $\Hgrrcc$. 

In section \ref{combi} we proceed to develop the combinatorial machinery to compute the central character of each such representation. This is done by constructions which generalize the description of the Springer correspondence in the equal label case. We associate a set of Springer correspondents to each central character $W_0c$ of $\Hgrrcc$. This set is shown to be in bijection with the set of central characters of the generic algebra which specialize into $W_0c$ upon specialization of the algebra. This bijection yields a combinatorial explanation for Lusztig's symbols (and their generalization), in terms of split Young tableaux. It supports the belief that generically, the modules in $\Hgrrcc$ are separated by their central character (for type $B$ and $C$). Then, we define a set of partitions which replace the unipotent conjugacy classes of $G$, and a (Springer-like) map associating a set of Weyl group characters to each such ``unipotent class''. We show that, as in the equal label case, the unipotent classes are in bijection with the set of central characters of $\Hgrrcc$. Finally we refine the conjecture to include also the central character, in Conjecture \ref{finalconj}. 

We end by giving some explicit examples of the combinatorial constructions and a particular case of the conjectures.

{\bf Acknowledgements.} The author would like to thank Eric Opdam for his invaluable mathematical and moral support. Part of this research was carried out with the support of grant number RTN2-2001-00059 of the Algebraic Combinatorics in Europe Network.

\section{The graded and the affine Hecke algebra}\label{hecke}

In this section we recall the definitions of both types of Hecke algebras, in order to describe the set of central characters of irreducible tempered representations of the graded Hecke algebra. This is a translation of the corresponding description, due to Opdam in \cite{O1}, for the affine Hecke algebra. This translation is done by using Lusztig's results in \cite{lusztig}, which describe the relation between the (representation theory of) the graded and the affine Hecke algebra.


\subsection{The graded Hecke algebra}We begin by fixing some notation and introding the graded Hecke algebra.\subsubsection{Definition}
Let $\Rdatum=(R_0,X,\check{R}_0,Y,\Pi_0)$ be a reduced, irreducible root
datum (notation as in \cite{lusztig}, or \cite{O1}). In
particular, $X$ and $Y$ are finitely generated free abelian groups
with a perfect pairing $\langle \cdot,\cdot \rangle$. The set of
roots $R_0 \subset X$ satisfies the usual axioms, and generates
the root lattice $Q=\Z R_0\subset X$; suppose $\mbox{rank}(Q)=n$.
The set $\check{R_0}\subset Y$ (the coroots) is in bijection with $R_0$,
where this bijection $\a \leftrightarrow \ca$ is such that
$\langle \a,\ca\rangle=2$. In particular, for every $\a \in R_0$
we have a linear map $s_{\ca}: X \to X: x \mapsto x-\langle x,\ca
\rangle \a$. Finally, the set $\Pi_0$ is a choice of $n$ roots
$\a_1,\dots,\a_n$ such that every root in $R_0$ can be written as
a non-negative or non-positive linear combination of the $\a_i \in
\Pi_0$. Then the Weyl group $W_0=\langle s_{\ca} \mid \a \in
R_0\rangle$ satisfies $W_0=\langle s_{\ca}\mid \a \in \Pi_0\rangle$.
Put $I=\{1,2,\dots,n\}$.

Let $\t^*=X \otimes \C$ and $\t=Y \otimes \C$. Inside $\t^*$ we have the distinguished real vector space $\real^*=X\otimes_Z\R$, and similarly for $\real=Y\otimes_\Z \R \subset \t$. We define the degenerate root datum associated to $\Rdatum$ to be $\Rdatum^{deg}=(R_0,\real^*,\check{R}_0,\real,\Pi_0)$.

Choose formal
parameters $\vet{k}_\a$ for $ \a\in \Pi_0$ such that
$\vet{k}_\a=\vet{k}_\beta$ if $\a$ and $\beta$ are conjugate under
$W_0$. We denote by $\vet{k}$ the function $\a \mapsto \vet{k}_\a$. Then the graded (also called degenerate) Hecke algebra
$\tilde{\Hgr}=\tilde{\Hgr}(\Rdatum^{deg},\vet{k})$ is by
definition the tensor product of algebras
\begin{equation}
\tilde{\Hgr}=\C[W_0] \otimes S[\t^*] \otimes \C[\vet{k}_\a],
\end{equation}
subject to the cross relations that the $\vet{k}_\a$ are central and that
\begin{equation}\label{lusrelgr}
x\cdot s_\a- s_\a\cdot s_\a(x)=\vet{k}_\a\langle x,\ca \rangle
\end{equation}
for all $x \in \t^*, \a \in \Pi_0$.

Notice that $\tilde{\Hgr}$ is a graded algebra if we put $\t^*$ and the ${\vet k}_\a$ in degree one and $\C[W_0]$ in degree zero.

\subsubsection{Parabolic algebras}\label{para}Let $\Pi_P \subset \Pi_0$ generate the standard parabolic root system $R_P \subset R_0$ and let $\real_P \subset \real$ be the real span of $\check{R}_P$. Then we define two degenerate root data $\Rdatum_P^{deg}=(R_P,\real_P^*,\check{R}_P,\real,\Pi_P)$ and $\Rdatum^{P,deg}=(R_P,\real^*,\check{R}_P,\real,\Pi_P)$. Let $k_P$ denote the restriction of $k$ to $R_P$. Given $\Hgr=\Hgr(\Rdatum,k)$, we define $\Hgr_P=\Hgr(\Rdatum_P^{deg},k_P)$ and $\Hgr^P=\Hgr(\Rdatum^{P,deg},k_P)$.

\begin{rem}\label{scaling}
Consider the generic graded Hecke algebra $\Hgr(\Rdatum^{deg},\vet{k})$. Then $R_0$ is a root system in the real vector space $\real^*$. Suppose that $R_0'$ is another root system inside $\real^*$, whose Weyl group $W_0'$ satisfies $W_0=W_0'$. In particular, the reflecting hyperplanes must be the same for both Weyl groups and it follows that there is a unique $\phi:R_0 \to R_0'$ such that $\phi(\a) \in \R_{>0}\a$. Then there is a unique $w \in W_0$ such that $\Pi'_0=w(\phi(\Pi_0))$. We let $\a'=\phi(\a)=\l_\a\a$. Clearly $\l_\a=\l_\beta$ if $W_0\a=W_0\b$. 
We put $\vet{k'}_{\a'}=\l_\a^{-1}{\vet{k}}_\a$. 

We choose a root datum $\Rdatum'=(R'_0,X',\check{R}'_0,Y',\Pi'_0)$ such that $\Rdatum'^{deg}=\Rdatum^{deg}$. Let the above defined $\vet{k}'_{\a'}$ be the root labels for $R_0'$.  Then we define $\Psi:\Hgr(\Rdatum^{deg},\vet{k}) \to \Hgr(\Rdatum'^{deg},\vet{k'})$, by putting $s_\a \mapsto s_{w\phi(\a)}, x \mapsto w(x), \vet{k}_\a \mapsto \l_\a^{-1}\vet{k'}_{\phi(\a)}$. It is easy to see that $\Psi$ extends to an algebra isomorphism.


A first application is to take $\l \in \R^\times$ and to consider the scaling ${\vet k}_\a \to \l {\vet k}_\a$ for all $\a$.
We find that $\Hgr(\Rdatum^{deg},\vet{k}) \cong \Hgr(\Rdatum'^{deg},\l\vet{k})$ if $R_0'=\l R_0$ and $\Pi'_0=\l\Pi_0$ in $\real^*=X \otimes_\Z\R=X'\otimes_\Z $. However, since the isomorphism is the identity on $\C[W_0]$ and on $\t^*$, it follows that $\Hgr(\Rdatum^{deg},\vet{k}) \cong \Hgr(\Rdatum^{deg},\l\vet{k})$. Notice that it is not strictly necessary to assume that $\Pi'_0=\l^{-1}\Pi_0$, since a different choice of simple roots leads to isomorphic algebras. In particular, in a specialized algebra with root labels $k_\a \in \R$, we need only consider the $k_\a$ up to (simultaneous) scalar multiples.

A second application is to suppose that $R_0$ has type $B_n$, i.e. relative to the standard basis $\{e_i\}$ of $\real^*\subset\real^*_\Q=Q \otimes_\Z \R$, the simple roots can be chosen to be $\a_i=e_i-e_{i+1}$ for $i=1,\dots,n-1$ and $\a_n=e_n$. We choose $R'_0 \subset \real^*_Q$ of type $C_n$, with simple roots $\a_i'=\a_i$ for $i<n$ and $\a'_n=2e_n$. It follows that the graded Hecke algebra $\Hgr(B_n)$ associated to $(R_0,\real^*,\check{R}_0,\real,\Pi_0)$ 
and parameters determined by $k_1=k_{e_i-e_{i+1}}$ and $k_2=k_{e_n}$ is isomorphic to the graded Hecke algebra $\Hgr(C_n)$ associated to $(R'_0,\real^*,\check{R}'_0,\real,\Pi'_0)$ and parameters determined by $k_1'=k_{e_i-e_{i+1}}$ and $k'_2=k_{2e_n}$ where $k'_1=k_1$ and $k'_2=2 k_2$. In particular, the equal label case $k_1'=k_2'$ for $\Hgr(C_n)$ corresponds to the case $k_2=\half k_1$ for $\Hgr(B_n)$.
\end{rem}

\subsubsection{Tempered representations} The goal of this article is to make a conjecture on the representation theory of certain specializations of this algebra,
obtained by substituting the $\vet{k}_\a$ with $k_\a\in \R$. Thus, we obtain a map $k:R_0 \to \R: \a \to k_\a$, satisfying $k_\a=k_\b$ if $W_0\a=W_0\b$. We
denote such a specialized algebra by $\Hgr=\Hgr(\Rdatum^{deg},k)$ and refer
to it as the graded Hecke algebra as well. Let us review some well
known facts about the representation theory of $\Hgr$ in order to
fix notation. If $V$ is a finite-dimensional $\H$-module, then the
abelian algebra $S[\t^*]$ induces a weight space decomposition
\[ V = \bigoplus_{\l \in \t}V_\l,\]
where
\[ V_\l= \{ v \in V \mid \forall x \in \t^*,\ (x-\l(x))^kv=0 \mbox{ for some }k \in \N \}. \]
The $\l$ for which $V_\l \neq 0$ are called the weights of $V$. If
$V$ is irreducible, then by Dixmier's version of Schur's lemma,
the center $Z$ of $\Hgr$ acts by a character. Using the cross
relations one can show (\cite{lusztig}) that the center $Z$ of
$\Hgr$ is equal to $S[\t^*]^{W_0}$. It follows that any
irreducible representation is finite-dimensional, since $\Hgr$ is
finitely generated over $Z$. Let $\l \in \t$ be such that for all
$p \in Z$ and $v \in V$, we have $p \cdot v=p(\l)v$. Notice that
$\l$ is determined up to $W_0$-orbit only. We call $W_0\l$ (or any
of its elements) the central character of $V$.

We can now define the notion of temperedness. We define the positive Weyl chamber
\[ \t^{*+}=\{ \l \in \t^* \mid \langle \l,\ca_i\rangle \geq 0 \mbox{ for all } i
\in I \}, \] and its antidual
\[ \t_-=\{ \xi \in \t \mid \langle \l,\xi\rangle \leq 0 \mbox{ for all } \l \in
\t^{*+}\}. \] Then we define
\begin{defn} Let $V$ be a finite dimensional $\mathbb{H}$-module.

(i) $V$ is called tempered if every weight $\gamma$ of $V$
satisfies $Re(\gamma) \in \t_-$.

(ii) $V$ is called a discrete series representation if, for every
weight $\gamma$ of $V$, $Re(\gamma)$ lies in the interior of
$\t_-$.
\end{defn}

Notice that discrete series representations do not exist if $\mbox{rank}(R_0)<\mbox{rank}(X)$.

We will be interested in the following set.

\begin{defn}\label{defrcc}
If $(\pi,V)$ is an irreducible representation of $\Hgr$ with central character $W_0\l$, then we say that $(\pi,V)$ has real central character if $W_0\l \subset {\mathfrak a}$. We denote the set of equivalence classes of irreducible tempered representations of $\Hgr$ which have real central character by $\Hgrrcc$.
\end{defn}
\subsubsection{Residual subspaces} In our approach towards a classification of the tempered irreducible representations of the graded Hecke algebra, we use (following \cite{HO} and \cite{O1}) certain affine subspaces $L \subset \real$ which describe their central characters. These subspaces are defined as follows. Let $L \subset \real$ be an affine subspace. Define the parabolic root subsystem $R_L\subset R_0$ by
\[ R_L=\{\a \in R_0 \mid \a(L)={\rm constant}\}.\] Then we call $L$ a residual subspace if and only if
\begin{equation}
\label{res}
 |\{\a \in R_L \mid \a(L)=k_\a\}|=|\{\a \in R_L\mid \a(L)=0\}|+{\rm codim}(L).
\end{equation}

In particular, $\real$ itself is residual. It is clear that the
notion of residual subspace is $W_0$-invariant, since
$w(R_L)=R_{wL}$. The following important property, which reduces the classification of residual subspaces to that of residual points, was proved in \cite{HO}. If $L\subset \real$ is an affine subspace with $\mbox{codim}(L)=\mbox{rank}(R_L)$ then, putting
$\real_L=\check{R}_L \otimes_\Z\R$, we have $L=c_L+\real^L$, where $c_L=L
\cap \real_L$, and $\real^L=\real_L^\perp$.  We call $c_L$ the {\it center}
of $L$. Then $L$ is a residual subspace if and only if $\mbox{codim}(L)=\mbox{rank}(R_L)$ and its center $c_L$ is a residual point with respect to the graded Hecke algebra $\Hgr_L=\Hgr(\Rdatum_L^{deg},k_L)$, attached to the data $\Rdatum_L^{deg}=(R_L,\real_L^*,\check{R}_L,\real_L,\Pi_L)$ and labelling function $k_L$, the restriction of $k$ to $R_L$. Therefore, the determination of the residual subspaces boils down to the classification of residual points. This classification is done for all root systems in \cite{HO}. We remark that by Lemma 7.10 of \cite{O1} and induction on the rank of $R_0$, it is not hard to see that there are only finitely many residual subspaces.

The importance of the residual subspaces lies
in the following fact, that we will prove below using the affine
Hecke algebra. For a residual subspace $L=c_L+\real^L$, let $L^{\rm
temp}=c_L+i\mathfrak{a}^L \subset \t=\real_\C$ be the corresponding {\it tempered
form} of $L$, which we call a {\it tempered residual subspace}.

\medskip
The following theorem is the starting point for the combinatorics in the next sections. 

\begin{thm}\label{gradcc} The collection $\cup_L L^{\rm temp}$ of all tempered residual subspaces of ${\mathbb H}$ is equal to the set of central characters of irreducible tempered representations of ${\mathbb H}$. Moreover, the set of central characters of discrete series representations is equal to the set of $W_0$-orbits of residual points.
\end{thm}

Notice that in particular, since $k_\a \in \R$, all discrete series representations have real central character.

\subsection{The affine Hecke algebra} The proof of Theorem \ref{gradcc} relies on (1) theorems of Lusztig which establish important connections between the graded and the affine Hecke algebra, and (2) on its analogue  for the affine Hecke algebra, Theorem \ref{ODO}, due to Opdam and Delorme-Opdam.

\subsubsection{Definitions} Recall the root datum $\Rdatum=(R_0,X,\check{R_0},Y,\Pi_0)$. The elements of $X$ act on $X$ by translation, and it is easy to see that this action is normalized by the action of $W_0$. We define the extended affine Weyl group to be the product $W=W_0 \ltimes X$.  Clearly $W$ contains the subgroup $W^a$ generated by the affine reflections $s_a$, where $a=(\ca,k) \in R:=\check{R_0} \times \Z$, defined by \[ s_{(\ca,k)}(x)=x-(\langle x,\ca\rangle +k)\a. \]  We call $R$ the affine root system. Let $\check{\theta}$ be the highest coroot in $\check{R}_0$. Then $W^a=\langle s_a\mid a\in R \rangle= \langle s_i\mid i=0,\dots,n\rangle$ where $a_0=(-\check{\theta},1)$, $a_i=(\ca_i,0)$ for $i>0$ and we write $s_i=s_{a_i}$. Let $\Pi=\{a_i\mid i=0,\dots,n\}$. We call $W^a$ the affine Weyl group.
We have the decomposition $W=W^a \rtimes \Omega$, where $W^a=W_0 \ltimes Q$ is a Coxeter group and $\Omega\cong X/Q$.

Notice that the Dynkin diagram of the affine root system is the affine extension of the Dynkin diagram of the coroot system. In particular, if the root system $R_0$ is of type $B_n$, the set of affine roots $R=\check{R_0}\times \Z$ has type $C_n^{\rm aff}$. In this case, depending on the choice of $X$ one has either two or three orbits of simple roots under the action of $W$: if $\check{\theta} \in 2Y$ then there are three.

The affine Hecke algebra is a deformation of the group algebra of $W$. In order to define it we choose positive real numbers $q_{a}$ for every root in $R$, such that $q_a=q_{wa}$ for every $a\in R,w\in W$. Let $l$ be the length function of $W$; that is, the extension of the length function of the Coxeter group $W^a$ to $W$ by putting $l(\omega)=0$ for all $\omega \in \Omega$. Define, for a simple reflection $s_a$ of $W$, the number $q(s_a)=q_{a+1}$. In view of the assumptions on the $q_a$, there exists a unique $q:W \to \R_{>0}$ such that $q(ww')=q(w)q(w')$ if $l(ww')=l(w)+l(w')$, and $q(\omega)=1$ for all $\omega \in \Omega$.

\begin{defn}
We denote by $\H=\H(\Rdatum,q)$ the affine Hecke algebra associated to $(\Rdatum,q)$, i.e., the unique complex associative algebra with $\C$-basis $(T_w)_{w \in W}$ satisfying the following relations:
\begin{equation}
\left\{ \begin{array}{l} \mbox{If } l(ww')=l(w)+l(w') \mbox{ then } T_{ww'}=T_wT_{w'};
\\ \mbox{If }a \in \Pi, s=s_a, \mbox{ then } (T_s+1)(T_s-q(s))=0.
\end{array}\right.
\end{equation}
\end{defn}
It is well-known that such an object exists and is indeed unique. It admits a decomposition $\H=\H_0 \otimes \A$, where $\A$ is a commutative algebra which is isomorphic to $\C[X]$ through $x\mapsto \theta_x \in \A$. The cross relations  are known as the Bernstein-Zelevinsky-Lusztig relations. To write them down, it is convenient to define root labels for the possibly non reduced root system $R_{nr}$, where
\begin{equation}
R_{nr}=R_0 \cup \{ 2\a \mid \ca\in \check{R_0} \cap 2Y\}.
\end{equation}

We then label the coroots of $R_{nr}$. If $\ca\in \check{R}_{nr} \cap \check{R_0}$ then we put $q_{\ca}=q_{(\ca,0)}$ and if $\ca/2 \in \check{R}_{nr} \backslash \check{R_0}$ then we define
\begin{equation}
q_{\ca/2}=\frac{q_{(\ca,1)}}{q_{(\ca,0)}}.
\end{equation}
With this notation, we have for $s=s_\a, \a\in \Pi_0$ and $x \in X$:
\begin{equation}
\label{lusrel}
 \theta_xT_s-T_s\theta_{s(x)}=\left\{ \begin{array}{l@{\mbox{ if }}l} (q_{\ca}-1)\frac{\theta_x-\theta_{s(x)}}{1-\theta_{-\a}} & \ca \notin
 2Y \\ ((q_{\ca}-1)+\theta_{-\a}(q^{1/2}_{\ca}q^{1/2}_{\ca/2}-
 q^{1/2}_{\ca}q^{-1/2}_{\ca/2}))\frac{\theta_x-\theta_{s(x)}}{1-\theta_{-2\a}}
 & \ca \in 2Y \end{array} \right.
\end{equation}

 By an unpublished result of Bernstein (see \cite{lusztig}), the center $\cent$ of $\H$ can then be identified to be $\cent=\A^{W_0}$.

Let $T=\mbox{Hom}(X,\C^*)$, then $\mbox{Spec}(\cent)=T/W_0$. In an irreducible representation $(\pi,V)$ of $\H$, there exists a $t_\pi \in T$ such that any $z \in \cent$ acts by $\pi(z)=z(t_\pi)1_V$. We call $W_0t_\pi$, or by abuse of terminology, any of its elements, the central character of $(\pi,V)$.

The torus $T$ admits a polar decomposition $T=T_uT_{rs}=\mbox{Hom}(X,S^1)\mbox{Hom}(X,\R_{>0})$.
\begin{defn}
Let $(\pi,V)$ be a finite dimensional representation of $\H$.

(i) We call $V$ a tempered representation if all ${\mathcal
A}$-weights $t$ of $V$ satisfy $|t(x)| \leq 1$ for all $x \in
X^+=\{x \in X \mid \langle x, \ca \rangle\geq 0 \ \forall \
\a \in \Pi_0\}$.

(ii) We call $V$ a discrete series representation if all
${\mathcal A}$-weights of $V$ satisfy $|t(x)|<1$ for all $x\in
X^{++}=\{ x \in X \mid \langle x, \ca \rangle >0 \ \forall \
\a \in \Pi_0\}$.
\end{defn}

\begin{defn} Let $\Hrcc$ denote the set of equivalence classes of irreducible tempered representations of $\H$, whose central character lies in $T_{rs}$.
\end{defn}

\subsubsection{Residual cosets} By theorems of Opdam (cf. \cite{O1}) and Delorme-Opdam (cf. \cite{DO}), the set of central characters of irreducible tempered representations of $\H$ can be described in terms of the affine analog of the residual subspaces for $\Hgr$. Let $L \subset T$ be a coset of a subtorus $T^L \subset T$. Put $R_L=\{ \a \in R_0 \mid \a(T^L)=1\}$ and let $W_L$ be the corresponding parabolic subgroup of $R_0$. Then we define
\[ R_L^p=\{ \a \in R_L \mid \a(L)=-q_{\ca/2}^{1/2} \mbox{ or }\a(L)=q_{\ca/2}^{1/2}q_{\ca} \}\] and \[ R_L^z=\{ \a \in R_L \mid \a(L)=\pm 1 \}.\]
By definition, we call $L$ a residual coset if and only if
\[ |R_L^p|=|R_L^z|+\mbox{codim}(L).\]
A zero-dimensional residual coset is called a residual point. By Theorem 3.29 and Corollary 3.30 in \cite{O1}, a residual point is the common central character of a non-empty set of irreducible discrete series representations of $\H$. Conversely, the central character of an irreducible discrete series representation is always a residual point. 

For tempered representations in general, we need the notion of tempered form of a residual coset. Suppose that $L$ is a residual coset, and let $T_L$ be the subtorus of $T$ whose Lie algebra is spanned by $\check{R}_L$. Then we may choose $r_L$ in $L$ such that $L=r_LT^L$ and $r_L \in L \cap T_L$. We call $r_L$ the {\it center} of $L$. The tempered form of $L$ is defined to be $L^{temp}=r_LT^L_u$. Let $\L$ be the set of residual cosets (it is not hard to see that $\L$ is finite). We have
\begin{thm}\label{ODO}(Opdam, Delorme-Opdam)
\[ \cup_{L \in \L}L^{temp}=\{ \mbox{central characters of irreducible tempered
representations of }\H\}.\]
\end{thm}

\pf: The inclusion $\subset$ follows follows from \cite{O1}, Theorems 3.29 and 4.23. The other follows from \cite{DO} where it is shown that every irreducible tempered representation occurs in the support of the Plancherel measure for $\H$. \qed

For a residual coset $L$ with center $r_L$, it is not hard to see that in the polar decomposition $r_L=s_Lc_L \in T_uT_{rs}$, the real part $c_L$ is independent of the choice of $r_L$. We call a residual coset {\it real} if we can choose $s_L=1$.

Let $L=r_LT^L$ be a real residual coset, and let $c_L=\mbox{log}(r_L) \in \mbox{Lie}(T_{L,rs})$. Then $c_L+Lie(T^L_{rs})$ is a residual subspace for $\Hgr(\Rdatum^{deg},k)$ with root labels $k_\a$ as in \eqref{kq} below, and vice versa.

\subsubsection{Parabolic subalgebras} The proof of theorem \ref{gradcc} uses the parabolic nature of the classification of residual subspaces. Therefore we fix some more notation here. Let $R_L$ be a standard parabolic root subsystem of $R_0$ with simple roots $\Pi_L \subset \Pi_0$. Then $\mathcal{R}^L=(R_L,X,\check{R}_L,Y,\Pi_L)$ is a root datum.  Define also
\[ Y_L= Y \cap \R\check{R}_L \mbox{ and } X_L=X/(X \cap Y_L^\perp). \]
Then $\mathcal{R}_L=(R_L,X_L,\check{R_L},Y_L,\Pi_L)$ is a root datum satisfying $\mbox{rank}(R_L)=\mbox{rank}(X_L)$. On $R_{L,nr}:=\Q R_L \cap R_{nr}$ we define root labels $q_{L,\ca}=q_{\ca}^L$ by restricting $q$ from $R_{nr}$ to $R_{L,nr}$. We extend $q_L$, resp. $q^L$, to $W_0(R_L) \ltimes X_L$, resp. $W_0(R_L) \ltimes X$. 

Then we define the subalgebra $\H^L \subset \H$ to be $\H^L=\H(\Rdatum^L,q^L)$, the affine Hecke algebra associated to $(\mathcal{R}^L,q^L)$. We also define $\H_L=\H(\Rdatum_L,q_L)$.

\begin{rem}\label{temp}
For the graded Hecke algebra, recall the algebra $\Hgr_L=\Hgr(\Rdatum^{deg}_L,k_L)$ which is the analogue of $\H_L$. We also define $\Hgr^L=\Hgr(\Rdatum^{L,deg},k_L)$ where $\Rdatum^{L,deg}=(R_L,\real^*,\check{R}_L,\real,\Pi_L)$. Notice that $\Hgr^L=\Hgr_L \otimes S[\t^{L*}]$, where $\t^L$ is the orthogonal complement in $\t$ of $\t_L=Y_L \otimes \C$. In particular, the tempered representations of $\Hgr^L$ are obtained as $U \otimes C_\nu$ where $U$ is a tempered representation of $\Hgr_L$ and $C_\nu$ is the one-dimensional $\t^L$-module belonging to $\nu \in i\t^{L*}$ (see also \cite{evens}).

For the affine Hecke algebra, such a decomposition is not available since a sublattice $X_L$ in $X$ need not have a complement $X^L \subset X$ such that $X_L+X^L=X$.
\end{rem}

\subsubsection{Lusztig's reduction theorem}\label{bijections}

In \cite{lusztig}, Lusztig establishes the connection between the representation theory of the affine and graded Hecke algebra. Since these reductions play an essential role, we will explain them in some detail. It will be necessary to slightly adapt his constructions. Suppose the parameters of $\Hgr$ depend on those of $\H$ as follows (in view of \cite{O1}, (4.5)):
\begin{equation}
\label{kq}
k_\a=\begin{cases} \mbox{log}(q_{\ca}) & \mbox{ if } \ca \notin 2Y \\ \mbox{log}(q_{\ca}q_{\ca/2}^{1/2}) & \mbox{ if } \ca \in 2Y. \end{cases}
\end{equation}

Now fix a {\it real} central character $W_0t$ of an irreducible
representation of $\H$. Denote the corresponding maximal ideal of
$\mathcal{Z}$ by $I_t$, and let $\hat{\mathcal{Z}}_t$ be the
$I_t$-adic completion of $\mathcal{Z}$. Furthermore we put
\[ \hat{\H}_t=\H \otimes_{\mathcal{Z}}\hat{\mathcal{Z}}_t.\]
Denote the set of irreducible representations of $\H$ with central
character $W_0t$ by $\mbox{Irr}_t(\H)$. Then
\begin{equation}\label{Hiso1}
\mbox{Irr}_t(\H) \longleftrightarrow \mbox{Irr}(\hat{\H}_t),
\end{equation}
which follows from the fact that
\[ \H/I_t\H \cong \hat{\H}_t/\hat{I}_t\hat{\H}_t,\]
where $\hat{I}_t$ denotes the maximal ideal of
$\hat{\mathcal{Z}}_t$ corresponding to $I_t$.

Analogously, let $W_0\gamma$ be a $W_0$-orbit in $\t=Y \otimes \C$, then after
carrying out the analogous constructions for $\mathbb{H}$ we have
\begin{equation}\label{Hiso2}
 \mbox{Irr}_\gamma(\mathbb{H}) \longleftrightarrow
\mbox{Irr}(\hat{\mathbb{H}}_\ga).
\end{equation}

Let $\hat{\H}_{rcc}$ denote the set of equivalence classes of irreducible representations of $\H$ with real central character, and analogously for
$\hat{\mathbb{H}}_{rcc}$. We can now prove the following:
\begin{thm}\label{lus}(Lusztig)
There exists a natural bijection
\[ \{\mbox{irreducible representations in }\hat{\H}_{rcc}\}\longleftrightarrow \{\mbox{irreducible representations in } \hat{\mathbb{H}}_{rcc}\}.\]
\end{thm}

\pf: This is basically Lusztig's second reduction theorem (Theorem 9.3 in \cite{lusztig}). We explain how to adapt his construction,
since he works with the different assumptions that the root labels
$q_{\ca}$ are of the form $q_{\ca}=q^{n_\a}$ for some $q\in \C^*$ and $n_\a \in \N$,
and moreover that for all $t'\in W_0t$, one has $\a(t') \in \langle q\rangle$
(the group in $\C^*$ generated by $q$) if $\ca\notin 2Y$, and
$\a(t')\in\pm\langle q\rangle$ if $\ca\in 2Y$.

Fix a real central character $W_0t \subset T_{rs}$. Since the exponential map, restricted to $\mathfrak{a}=Y\otimes_\Z \R$, yields a $W_0$-equivariant isomorphism ${\rm exp}:\mathfrak{a} \to T_{rs}$, it gives rise to a bijection $W_0c \to W_0t$, where $c \in \mathfrak{a}$ is such that ${\rm exp}(c)=t$. In Lusztig's notation, this means that $t_0=1$. Our
assumption that all $q_{\ca}> 1$ then implies that Lemma 9.5
still holds. Hence, Theorem 9.3 still holds: the algebras
$\hat{\mathcal{Z}}_t$ and $\hat{Z}(\mathbb{H})_c$ are isomorphic, and
moreover \[ \hat{\H}_t \cong \hat{\mathbb{H}}_c. \] This isomorphism
is compatible with the $\hat{\mathcal{Z}}_t \cong
\hat{Z}(\mathbb{H})_c$-structures. We therefore find a natural
bijection between the irreducible representations of $\hat{\H}_t$
and the irreducible representations of $\hat{\mathbb{H}}_c$.
Combining with \eqref{Hiso1}, \eqref{Hiso2} and the fact that this
holds for any real central character $W_0t$, we find the desired
result.\qed

It is clear from the constructions that if the $\H$-module $V$ corresponds to the $\Hgr$-module $U$ under this bijection, then $V$ is a tempered (resp. discrete series) module if and only if $U$ is a tempered (resp. discrete series) module. Therefore we obtain in particular a bijection
\begin{equation}
\label{bijtemp}
\Hrcc \longleftrightarrow \Hgrrcc.
\end{equation}

\begin{lemma} \label{rcc}
Let $V$ be an irreducible discrete series representation of $\mathbb{H}$ with central character $\gamma$. Then $\gamma \in \mathfrak{a}$, i.e., $V$ has real central character.
\end{lemma}

\pf: We apply the analogue of Lusztig's first
reduction theorem (\cite{lusztig}, Theorem 8.6) to $\mathbb{H}$
(instead of $\H$). We adapt his definition of the root subsystem
$R_c$, and use the parabolic root system $R_\gamma$ of \cite{O1},
(4.15) instead. Explicitly,
$R_\gamma=R_0\cap\mathfrak{a}^*_\gamma$, where
$\mathfrak{a}^*_\gamma$ is the $\R$-span of $\{\a \in R_0 \mid
\a(\gamma) \in \{ 0,\pm k_\a\}\}$. We choose $\gamma$ in its
$W_0$-orbit such that $R_\gamma=R_P$, a standard parabolic
subsystem and put $c=W_P\gamma$. Let $\Gamma \subset W_0$ be the
group $\{w \in W_0 \mid w(c)=c \mbox{ and } w(R_P^+)=R_P^+\}$, and
$n=|W_0\gamma|/|W_P\gamma|$. Then, according to \cite{lusztig}
(see also \cite{O1}, Theorem 4.10),
\begin{equation}\label{st86} \hat{\mathbb{H}}_\ga \cong (\hat{\mathbb{H}}^P_\ga[\Gamma])_n, \end{equation}
where $\hat{\mathbb{H}}_\ga$ is the completed algebra as in
\eqref{Hiso2}, and the right hand side denotes the algebra of
$n\times n$-matrices with entries in $\hat{\mathbb{H}}^P_\ga[\Gamma]$.
The upshot is that $V={\rm
Ind}_{\mathbb{H}^P[\Gamma]}^{\mathbb{H}} U$, where $U$ is an
irreducible $\mathbb{H}^P[\Gamma]$-module. Furthermore, by
Lusztig's explicit description of the isomorphism \eqref{st86}, we
obtain the following formula for the set $\mathrm{Wt}(V)$ of
weights of $V$:
\begin{equation}\label{gewU} \mathrm{Wt}(V)=\bigcup_{w\in W^P}w\cdot \mathrm{Wt}(U),\end{equation}
where $W^P$ denotes the set of coset representatives of $W_0/W_P$ of minimal length.

Suppose that $\gamma \notin \mathfrak{a}$, then $R_\gamma(=R_P)
\neq R_0$. The assumption that $V$ is a discrete series module
implies that $\gamma' \in \mathrm{Wt}(U)$ satisfies $Re(\gamma')
\in \sum_{i \in I}\R_{<0}\ca_i$. Let $w_P$ (resp. $w^P$) be the longest
element of $W_P$ (resp. $W^P$), then  $Re(w^P\gamma') \notin \sum_{i \in
I}\R_{<0}\ca_i$ while on the
other hand by \eqref{gewU}, $w^P\gamma' \in \mathrm{Wt}(V)$. This
is a contradiction, so the Lemma follows.\qed

\subsubsection{Proof of Theorem \ref{gradcc}:}
 Let $L^{temp}=c_L+i\mathfrak{a}^L$ be a tempered
residual subspace of $\mathbb{H}$, and let $c_L+ix \in
L^{temp}$. Then $r_L={\rm exp}(c_L)$ is a residual
point for $\H_L$, and by \cite{O1} it is then the central
character of a non-empty set $\Delta_L$ of discrete series representations of
$\H_L$. Since $r_L$ is real, we can invoke \ref{lus} to see that $c_L$ is the central
character of a set of discrete series representations of
$\mathbb{H}_L$, which is in natural bijection with $\Delta_L$. It
follows that $c_L+ix$ is the central character of a tempered
representation of $\mathbb{H}^L$, and (by induction) also of
$\mathbb{H}$.

Conversely, let $V$ be an irreducible tempered representation of $\mathbb{H}$ with central character $W_0c$. Let the set of weights of $V$ be $\rm{Wt}(V)$. We then define, analogous to the analysis in \cite{DO} of (weak) constant terms of tempered representations of the affine Hecke algebra, the sets
\[ E_P^0(V)=\{ \l \in {\rm{Wt}}(V) \mid Re(\l) \in \R \check{R}_P\},\]
where $P\subset I$.
Now let $P$ be minimal such that $E_P^0(V)\neq \emptyset$. Then it is easy to see that
\[ V^P=\sum_{\l \in E_P^0(V)}V_\l \]
is a non-zero tempered $\mathbb{H}^P$-module, and minimality of $P$ implies that for all $\l \in {\rm{Wt}}(V^P)=E_P^0(V)$, we have $\l \in \sum_{\a \in \Pi_P}\R_{<0}\ca$.
Now take an irreducible subquotient $U$ of $V^P$, then ${\rm{Wt}}(U) \subset {\rm{Wt}}(V^P)$, so $U$ is an irreducible, tempered $\mathbb{H}^P=\mathbb{H}_P\otimes S[\t^{P*}]$-module, which implies (see Remark \ref{temp}) that $U=U_P \otimes \C_\nu$, where $U_P$ is an irreducible discrete series module of $\mathbb{H}_P$. Since $U$ is a tempered $\mathbb{H}^P$-module, $Re(\nu)|_{\t^{P*}}=0$, so $\nu \in i\mathfrak{a}^P$. By Lemma \ref{rcc}, $U_P$ has real central character $W_P\l$. We choose $c$ in its $W_0$-orbit such that $c=\l+\nu$. Finally, we apply Lusztig's theorem \ref{lus}, combined with \cite[Theorem 3.29]{O1}, to conclude that $\l$ is a residual point for $\mathbb{H}_P$. It follows that $c \in \l+i\mathfrak{a}^P$, i.e., indeed $c$ lies in a tempered residual subspace of $\mathbb{H}$.\qed

\begin{cor} The set $\Hgrrcc$ is finite.
\end{cor}

\pf: There are only finitely many equivalence classes of irreducible discrete series representations of $\Hgr$ whose central character is a prescribed residual point. Therefore, the number of residual points being finite, there are only finitely many equivalence classes of irreducible discrete series representations of $\Hgr$. Secondly, every irreducible tempered representation $V$ of $\Hgr$ occurs as a summand of the unitary induction of a representation $U \otimes \C_\nu$ of some $\Hgr^L$, where $U$ is a discrete series representation of $\Hgr_L$. If $V$ has real central character, then $\nu=0$. \qed


\section{(Generalized) Green functions}\label{greenfunctions}
In this section, we recall the description of $\Hgrrcc$ in the classical case where the root labels $k_\a$ are all equal. In this case, there are $|\hat{W}_0|$ such modules. They are in fact none other than the Springer modules, hence they have a natural grading as $W_0$-module. One can compute their $W_0$-strucure in every degree with Green functions. If the root system is of classical type, these Green functions are the solutions of a matrix equation which can be given entirely in terms of combinatorial objects. For types $B$ or $C$, we define the analogue of these objects for other choices of parameters, and  define Green functions for other parameter values of $\Hgr$. These functions are a special case of a more general notion of Green functions introduced by Shoji (cf. \cite{shoji}).

\subsection{Special and generic parameters}\label{spec}

One may also use \eqref{res} to define residual subspaces for the generic algebra $\Hgr(\Rdatum^{deg},\vet{k})$. Upon specialization $\vet{k}_\a \mapsto k_\a$, we distinguish between the following possibilities.



\begin{defn}
Let $k_\a \in \R$ be a choice of root labels of the root datum $\Rdatum$ and let $\Hgr$ be the associated graded Hecke algebra. For a subset $\Pi_P \subset \Pi_0$, let $\Hgr_P=\Hgr(\Rdatum_P^{deg},k_P)$ (see \ref{para}) and let the corresponding generic algebra $\tilde{\Hgr}_L(\Rdatum_L^{deg},\vet{k}_L)$ have residual points $\{W_Pc_{P,1}(\vet{k}),W_Pc_{P,2}(\vet{k}),\dots,W_Pc_{P,i_P}(\vet{k})\}$.  Then we call the parameters $k_\a$ {\it generic} if, for all $\Pi_P \subset \Pi_0$, the evaluation map $W_Pc_{P,i}(\vet{k}) \mapsto W_Pc_{P,i}(k)$ is a bijection onto the residual points of $\Hgr_P$. Otherwise, we call them {\it special}.
\end{defn}

\begin{ex} We give the special parameters for the root systems of type $B_n,C_n$. Let $V$ be an $n$-dimensional real vector space with orthonormal basis $e_1,\dots,e_n$ and consider $R_0=\{\pm e_i\pm e_j\}\cup \{\pm e_j\}$, the root system of type $B_n$. Let the set $\Pi_0$ of simple roots consist of  $\a_i=e_i-e_{i+1}$ for $i=1,\dots,n_1$ and $\a_n=e_n$. Choose root labels $k_1=k_{\pm e_i \pm e_j}$ and $k_2=k_{\pm e_i}$.

The parameters $k_1,k_2$ are shown to be special in \cite{HO} if
\begin{equation}
\label{special}
k_1k_2\prod_{j=1}^{2(n-1)}(2k_2-jk_1)(2k_2+jk_1)=0,
\end{equation}
and generic otherwise.

In particular, notice that in view of Remark \ref{scaling}, we may view these special values as the union of the integral relations $k_2=jk_1$ and $k_2'=jk_1'$ ($j=1,\dots,n-1$) of the graded Hecke algebras of type $B_n$ and $C_n$.

Since the choices $k_1=k_2$ resp. $k_1=2k_2$ give rise to a graded Hecke algebra with equal labels, we call any of these cases an {\it equal label case}.
\end{ex}

\subsection{The Springer correspondence}\label{groupcases}

In the equal label case,  $\hat{\H}$ has been determined by Kazhdan and Lusztig in \cite{KL}. They assume that $X=P$. In view of Lusztig's theorem and in particular the bijection \eqref{bijtemp}, we give their description for $\Hgr$ rather than for $\H$. Let $\H$ have labels $q_{\ca}=q>1$ for all $\a$, and put $k=\mbox{log}(q)$ (in accordance with \eqref{kq}).

\subsubsection{Springer correspondence}Let $\hat{G}$ be the complex reductive group whose root system is $(R_0,X,\check{R_0},Y,\Pi_0)$. 
For a unipotent $u \in \hat{G}$, let $\B_u$ be the variety of Borel subgroups of $\hat{G}$ containing $u$. Putting $\mbox{dim}(\B_u)=d_u$, the highest non-vanishing cohomology group of $H(\B_u)$ is $H^{2d_u}(\B_u)$. Springer has shown that $H(\B_u)$ is a $W_0$-module, and that the action of $W_0$ respects the grading.  The group $A(u)=C_{\hat{G}}(u)/Z_{\hat{G}}C_{\hat{G}}(u)^0$ also acts on $H(\B_u)$, such that the actions of $W_0$ and $A(u)$ commute. Let, for $\rho \in \hat{A}(u)$, $H(\B_u)_\rho$ be the $\rho$-isotypic part in $H(\B_u)$. Springer has shown that (i) if $H^{2d_u}(\B_u)_\rho \neq 0$, then it is an irreducible $W_0$-module (whose character we will denote by $\chi_{u,\rho}$); (ii) If $\chi_{u,\rho}=\chi_{u',\rho'}$ then $(u,\rho)$ and $(u',\rho')$ are conjugate under $\hat{G}$; (iii) for every $\chi \in \hat{W_0}$ there exists a pair $(u,\rho)$ such that $\chi=\chi_{u,\rho}$.

Given $u$, we denote by $\hat{A}(u)_0$ the set of irreducible characters $\rho \in \hat{A}(u)$ for which $H(\B_u)_\rho \neq 0$. In general, $\hat{A}(u)_0$ is strictly smaller then $\hat{A}(u)$, but never empty:  it is known to always contain the trivial representation. We write
\begin{equation}
\label{I0}
\I_0:=\{ (u,\rho) \mid u \in \hat{G} \mbox{ unipotent},\ \rho \in \hat{A}(u)_0\},
\end{equation}

where we identify conjugate pairs. Then $\I_0 \leftrightarrow \hat{W}_0$; this bijection is called the Springer correspondence. Given a unipotent class $C \subset G$, we denote the set of its Springer correspondents by
\[ \Sigma(C)=\{ \chi_{u,\rho} \mid \rho \in \hat{A}(u)_0, u \in C\} \subset \hat{W}_0. \]
\begin{thm}\cite{KL}\label{KLstelling}
$\Hgrrcc$ is naturally parametrized by $\I_0$ (and hence also by $\hat{W}_0$).
 Let $M_{u,\rho}$ be the $\Hgr$-module which corresponds to $(u,\rho) \in \I_0$. Then $M_{u,\rho}$ and $M_{u',\rho'}$ have the same central character if and only if $u=u'$. \end{thm}
We will also write $M_\chi$ for $M_{u,\rho}$ if, through the Springer correspondence, $\chi=\chi_{u,\rho}$.

\subsubsection{Bala--Carter bijection}\label{BC}It follows from Theorem \ref{gradcc} and Theorem \ref{KLstelling} that in the equal label cases, the central characters
of $\Hgrrcc$ are indexed both by the $W_0$-orbits $W_0c_L$ of the centers of the residual subspaces $L$ of
$\Hgr$, and by the unipotent conjugacy classes in $\hat{G}$. Therefore, there exists a bijection between these two sets, such that if $W_Lc_L$ corresponds to the unipotent class $C \subset \hat{G}$, then the modules in $\Hgrrcc$ with central character $W_0c_L$ are the modules $M_{u,\rho}$ such that $u \in C$. 

From \cite[Appendix]{O1} we recall the description of the resulting bijection. It is the map 
\begin{equation}
\label{u}
 u \mapsto \ga(u) \in \real 
\end{equation}
where $\ga(u)$ is determined from $u$ as follows: suppose that the weighted Dynkin diagram of $u$ is labelled $(x_1,\dots,x_n)$, then $\ga(u) \in \real$ is determined by the equations $\a_i(\ga(u))=\frac{k}{2}x_i$.

This map is such that if $u$ corresponds via the Bala-Carter classification (cf. \cite[Theorem 5.9.5]{carter}) to the pair $(L,P_L)$ where $L$ is a Levi subgroup of $\hat{G}$ and $P_L$ is a distinguished parabolic subgroup of the semisimple part of $L$, then $W_0\ga(u)=W_0c_M$ where $M$ is the residual subspace with root system $R_M=R_L$. In particular, we obtain a bijection between distinguished unipotent classes of $\hat{G}$ and $W_0$-orbits of residual points of $\Hgr$.

\subsubsection{Structure of homology modules}


Being a $\Hgr$-module, $M_{u,\rho}$ is also a $W_0$-module. Recall that $k>0$. According to \cite{cells4, cusploc2} (see also \cite{reeder2}) we have (letting $\e$ denote the sign representation of $W_0$):
\begin{equation}
\label{r}
M_{u,\rho}\mid_{W_0}=\e \otimes H(\B_u)_\rho.
\end{equation}
In particular the elements of $\Hgrrcc$, viewed as $W_0$-modules, are graded by homological degree.
We will see in Corollary \ref{mult1} that every element of $\Hgrrcc$ contains a multiplicity free $W_0$-type.

\subsection{Green functions} 
Let $q=p^l$ for a prime $p$. Let $k$ be an algebraically closed field of characteristic $p$. Let $G$ be a connected reductive algebraic group defined over $k$, with a split $\F_q$-structure, and corresponding Frobenius map $F$. Then $G^F$ is a finite group of Lie type. We assume that $G$ and $G^F$ have root system $(R_0,X,\check{R_0},Y,\Pi_0)$ of classical type.

\subsubsection{Definition} Let $C$ be a unipotent class
in $G$, then, since we assume that $G$ is split, $C$ is
$F$-stable. It is known (by \cite{shoji83}) that there exists a 
representative $u \in C^F$, unique up to $G^F$-conjugacy modulo
the center of $G$, such that $F$ acts trivially on the
component group $A(u)=C_G(u)/C_G(u)^0$, and the set of
$G^F$-conjugacy classes in $C^F$ is in bijective correspondence
with the set $A(u)/\sim$ of conjugacy classes of $A(u)$. We denote
by $u_a$ a representative of the $G^F$-conjugacy class in $C^F$
represented by $a \in A(u)/\sim$. Finally, let
$B_u$ denote the variety of Borel subgroups of $G$
containing $u$. Since $\hat{G}$ and $G$ have the same root datum, we have $H(B_u)=H(\B_u)$ (where for simplicity $u$ also denotes the corresponding unipotent element in $\hat{G}$ which has the same Jordan decomposition as $u \in G$). For $w \in W_0$ we define the Green function $Q_w$ as
\begin{equation} \label{geogreen} Q_w(g)=\sum_{m=0}^{d_u}  {\rm{Tr}}((w,a),H^{2m}(\mathcal{B}_u))q^m,
\end{equation}
if $g \in G^F_{\rm uni}$ is $G^F$-conjugate to $u_a$.

\subsubsection{The matrix equation}
In \cite{luscarV}, Lusztig finds an algorithm for computing the
Green functions as solutions of a matrix equation. We briefly review it here.
Let $\chi \in \hat{W}_0$, and let $\chi \leftrightarrow i \in \I_0$ under the Springer correspondence. Then one defines a function $Q_i=Q_\chi$, also called Green function, by
\[ Q_\chi=\frac{1}{|W_0|}\sum_{w\in W_0}\chi(w) Q_w.\]
Define for each $i=(u,\rho) \in {\mathcal I}_0$, a
$G^F$-invariant function $Y_i$ on $G^F_{\rm uni}$ by
\[ Y_i(g)= \left\{ \begin{array}{ll} \rho(a) & {\rm if}\ g \ {\rm is}\ G^F-{\rm conjugate\ to\ }u_a \\ 0 & {\rm if\ }g \notin C^F, \end{array} \right.  \]
where $u \in C$.
It has been shown by Lusztig (\cite{luscarV}) that there exist $\pi_{ji}$
such that
\begin{equation}
\label{defpi}
 Q_i = \sum_{j \in {\mathcal I}_0} \pi_{ji}Y_j.
\end{equation}
A priori, the $\pi_{ji}$ belong to $\bar{\mathbb{Q}}_l$, but
\begin{prop}\cite{lusgreen,shojigm}
There exist polynomials $\vet{\pi}_{ij}$ in $\Z[t]$, which are independent of $q$, such that $\pi_{ij}={\vet \pi}_{ij}(q)$.
\end{prop}

\pf: In view of \eqref{geogreen} and \eqref{defpi}, if $j=(u,\rho)$ and $Q_i=Q_\chi$ then
the entry $\pi_{ji}$ can be expressed as
\begin{equation} \label{pi}
\pi_{ji}=\sum_{m \geq 0}\langle H^{2m}(\mathcal{B}_u),\chi \otimes
\rho \rangle q^m,
\end{equation}
where $\langle H^{2m}(\mathcal{B}_u),\chi \otimes
\rho \rangle$ denotes the multiplicity of $\chi \otimes \rho$ in $H^{2m}(\B_u)$.
Moreover, it is known that the structure of
$H^{2m}(\mathcal{B}_u)$ as a $W_0 \times A(u)$-module is
independent of $p$. \qed

Given a unipotent $u \in G$, we denote the unipotent class in which it lies by $C_u$. We define a pre-ordering on $\I_0$ such that $i=(u,\rho) \leq i'=(u',\rho')$ if $C_u \subset \overline{C_{u'}}$. Let $\sim$ be the associated equivalence relation $i \sim i' \iff C_u=C_{u'}$, and refine the pre-ordering into a total one.

We can now describe the matrix equation. Let $t$ be an indeterminate. 
For a class function $f$ on $W_0$, we define $R(f)$ by
\[ R(f)=(t-1)^nP_0(t)\frac{1}{|W_0|}\sum_{w\in W_0}\frac{\e(w)f(w)}{\mbox{det}(t\cdot \mbox{id}_V-w)},\]
where $P_0(t)$ is the Poincar\'e polynomial, and $\e$ is the sign character of $W_0$. For an irreducible character $\chi \in \hat{W}_0$, $R(\chi)$ is the fake degree of $\chi$.
We define the matrix $\Omega$ to be the matrix whose entries are
\[ \omega_{A,B}=t^nR(\chi_A \otimes \chi_B \otimes \e).\]

The following theorem is the basis of our conjecture.
\begin{thm} \label{main}(\cite{luscarV})
Let $P=(p_{i,i})_{i,i \in \I_0}$ and $\Lambda=(\lambda_{i,j})_{i,j \in \I_0}$ be matrices of unknowns, subject to
\begin{equation}\label{greeneq}
\left\{ \begin{array}{l} \l_{i,j}=0  \mbox{ unless } i \sim j \\ p_{i,j}=0  \mbox{ unless either }i < j, i\nsim j, \mbox{ or } i=j \\ p_{i,i}=t^{d_u}, i=(u,\rho) \\ P \Lambda \ ^tP=\Omega.\end{array} \right.
\end{equation}
This equation has unique solution matrices. The matrix $P$ satisfies $P_{ij}(q)=\pi_{ij}$. In particular, $P_{ij} \in \Z[t]$.
\end{thm}

\subsection{Generalization}

For a graded Hecke algebra associated to a root system $R_0$ of type $B$, we will formulate the analogue of equation \eqref{greeneq} for other choices of its parameters $k_\a$. The particular choice of this generalization is based on the combinatorics in section \ref{combi} where we describe a combinatorial generalization of the Springer correspondence.

\subsubsection{Notation} We write $\P_{n,2}$ for the set of double partitions of $n$, that is, the set of pairs $(\xi,\eta)$ where $\xi$ and $\eta$ are partitions such that $|\xi|+|\eta|=n$. Let $\P_n$ be the set of partitions of $n$. We will usually write the parts of a partition in increasing order; parts of length zero are allowed.

Let $W_0$ be the Weyl group of type $B_n$ and $\hat{W_0}$ the set of its irreducible characters. It is well known that $\hat{W_0}$ is in bijection with $\P_{n,2}$. We adopt the explicit choice of this bijection, where the representation indexed by $(\xi,\eta)$ has a basis indexed by standard Young tableaux of shape $(\xi,\eta)$, as in \cite{hoefsmit}. In particular, the trivial representation is indexed by $(n,-)$ and the sign representation by $(-,1^n)$. We denote an ireducible representation of $W_0$ by $(\pi_A,V_A)$ for $A\in \P_{n,2}$, and its character by $\chi_A$ or sometimes simply by $A$. For the one-dimensional sign representation, we will sometimes also write $\e$ for $\chi_{(-,1^n)}$, or even for $V_{(-,1^n)}$.

\subsubsection{Symbols}

The Springer correspondence for classical groups of type $B,C$ has been described by Lusztig (in \cite{luscells}) in terms of certain objects which he calls symbols. These have been generalized by Malle (see \cite{malle}) and Shoji (see \cite{shoji}) in order to apply to complex reflection groups. We will use a certain type of these generalized symbols.

Fix $n$ and let $(m_1,m_2)$ be a pair of positive integers, both at least equal to $n$. Furthermore choose two non-negative integers $r,s$. Define $\Lambda_1=(0,r,2r,\dots,(m_1-1)r) \in \R^{m_1}$ and $\Lambda_2=(r+s,2r+2,\dots,(m_2-1)r+s) \in\R^{m_2}$. Let $(\xi,\eta) \in \P_{n,2}$, and let (by adding zeroes if necessary) $l(\xi)=m_1$ and $l(\eta)=m_2$. Then we denote by $\Lambda^{m_1,m_2}(\xi,\eta)$ the array with two rows, whose upper row consists of the entries of $\Lambda_1+\xi$ and whose lower row consists of the entries of $\Lambda_2+\eta$; the entries of $\xi$ and $\eta$ are not written on top of each other but interleave. If $m_1>m_2$, we start the array with the entry $\xi_1$, if $m_1<m_2$ with the entry $\eta_1+s$. If $m_1=m_2$, we obtain two types of symbols, written $\Lambda^+(\xi,\eta)$ resp. $\Lambda^-(\xi,\eta)$ if we start the symbol with $\xi_1$ resp. $\eta_1+s$.

For example, if  $n=5$, $(\xi,\eta)=(1^3,2), r=2, s=1$, then
\[\Lambda^{6,5}(\xi,\eta)= \left( \begin{array}{lllllllllll}0&&2&&4&&7&&9&&11 \\ &1&&3&&5&&7&&11& \end{array} \right), \]
while
\[\Lambda^{5,6}(\xi,\eta)= \left( \begin{array}{lllllllllll} &1&&3&&5&&7&&11&\\0&&2&&4&&7&&9&&11  \end{array} \right). \]

In fact we are not so much interested in these symbols, but rather in their equivalence classes generated by the shift operation $(m_1,m_2)\to (m_1+1,m_2+1)$. Putting $m=m_1-m_2$, we denote the equivalence class of $\Lambda^{m_1,m_2}(\xi,\eta)$ by $\bar{\Lambda}^{m}(\xi,\eta)$. We call this class, or any of its elements the symbol of $(\xi,\eta)$. We denote the set of symbols $\bar{\Lambda}^m(\xi,\eta)$ of $(\xi,\eta)\in\P_{n,2}$ by $\bar{Z}^{r,s}_n(m)$.

The symbols generate an equivalence relation on $\hat{W_0}$. For $A,B \in \P_{n,2}$, we say that $A \sim_m B$ if the $m$-symbols $\bar{\Lambda}^m(A)$ and $\bar{\Lambda}^m(B)$ have representatives which contain the same entries with the same multiplicities. Denote the equivalence class of $A$ by $[A]_m$. We transfer this equivalence relation to $\hat{W}_0$ by putting $\chi_A \sim_m \chi_B \iff A \sim_m B$.

The Springer correspondence for $SO_{2n+1}(k)$ and $Sp_{2n}(k)$ (if char($k$)$\neq$2) has been described by Lusztig (\cite{luscells}) in terms of symbols with $m=1$, and $(r,s)=(2,0)$ resp. $(2,1)$. In this description, the sets of Springer correspondents of the various unipotent classes are the equivalence classes under $\sim_m$ in $\hat{W}_0$.


\begin{defn}
Let $n \geq 2$ and $m \in \half\Z$. For $m \neq 0$, we define the $m$-symbol of $(\xi,\eta) \in \P_{n,2}$ to be $\bar{\Lambda}^m(\xi,\eta) \in \bar{Z}_n^{2,0}(m)$ if $0\neq m\in \Z$, or $\bar{\Lambda}^{m'}(\xi,\eta)\in \bar{Z}^{2,1}_n(m')$ where $m'=sgn(m)(|m|+\frac{1}{2})$ if $m \notin \Z$. 

If $m=0$, we define the $\pm 0$-symbol of $(\xi,\eta)$ to be $\bar{\Lambda}^\pm(\xi,\eta) \in \bar{Z}^{2,0}_n(0)$.
\end{defn}

In other words, we will consider for a given $m \in \half\Z$ only one choice of $(r,s)$, which depends on the integrality of $m$. Therefore we suppress $(r,s)$ from the notation and always write $\bar{\Lambda}^m(\xi,\eta)$ for the symbol of $(\xi,\eta)$.

\begin{defn}\label{e_m}
 Let $(\xi,\eta)$ be a double partition with $m$-symbol $\bar{\Lambda}^m(\xi,\eta)$. Denote the entry of $\xi_i$ in the $m$-symbol of $(\xi,\eta)$ by
$e_m(\xi_i)$, and analogously for $\eta$.
\end{defn}

Notice that $e_m(\xi_i),e_m(\eta_j)$ depend on the chosen representative of the length of the $m$-symbol. However we will be only interested in differences of the form $e_m(\xi_i)-e_m(\eta_j)$, which are independent of the choice of the symbol lengths.

\subsubsection{The $a$-function and truncated induction}
An important tool related to symbols is the $a$-function. Suppose that $r,s$ and $m_1,m_2$ have been fixed. For $(\xi,\eta) \in \P_{n,2}$, $m$-symbol $\Lambda^{m_1,m_2}(\xi,\eta)$ can be written as $\Lambda^{m_1,m_2}(-,-)+(\xi,\eta)$. Then we define
\[ a_{m_1,m_2}(\xi,\eta)=\sum_{x,y \in \Lambda^{m_1,m_2}(\xi,\eta)} \mbox{min}(x,y)- \sum_{x,y \in \Lambda^{m_1,m_2}(-,-)}\mbox{min}(x,y).\]

This function is actually invariant for the shift operation, and thus it induces a function on the set of $m$-symbols $\bar{Z}_n^{r,s}(m)$. We denote it by $a_m$. Clearly, $a_m$ is constant on similarity classes. We transport $a_m$ to $\hat{W}_0$ by putting $a_m(\chi)=a_m(\xi,\eta)$ if $\chi=\chi_{(\xi,\eta)}$. 

Originally, for $\chi \in \hat{W}_0$, the number $a(\chi)$ was defined to be the highest $s$ such that $u^s$ divides the generic degree $D_\chi(u)$ of $\chi$. For the root systems of type $B_n$ and $C_n$ this function coincides with the $a$-function for $(r,s)=(1,0)$ and $m=1$, as a special case of the definition above.

Given such $m$, we use the $a_m$-function to define an ordering $\succ_m$ on $\hat{W}_0$. First we choose a partial order on $\hat{W}_0$ by demanding that $\chi \succ_m \chi'$ implies that $a_m(\chi) \leq a_m(\chi')$, and that similarity classes form intervals. Then we arbitrarily refine this partial ordering into a total one.

We use the $a$-function to define a truncated induction that we will need later. Let $W' \subset W_0$ be a subgroup and let $\chi'$ be a character of $W'$. Let
\begin{equation}
\mbox{Ind}_{W'}^{W_0}(\chi')=\sum_{\chi \in \hat{W_0}}n_{\chi',\chi} \chi,
\end{equation}
and suppose that $\chi_0 \in \hat{W}_0$ is such that $n_{\chi',\chi_0}>0$ and $n_{\chi',\chi} >0$ implies $a_m(\chi)\leq a_m(\chi_0)$. Then we define
\begin{equation}
\trm_{W'}^{W_0}(\chi')=\sum_{\chi:a_m(\chi)=a_m(\chi_0)}n_{\chi',\chi}\chi.
\end{equation}


\subsubsection{Generalized Green functions}
In the equal label cases one considers the $m$-symbols for $m=\frac{1}{2},1$. It is known that $a_m(\chi)=\mbox{dim}(\B_u)$, if $\chi=\chi_{u,\rho}$. If $C \subset \bar{C'}$ then $a_m(\chi) \geq a_m(\chi')$ for all $\chi \in \Sigma(C), \chi' \in \Sigma(C')$. 

 Shoji has defined in \cite{shoji} a matrix of functions $P=(p_{A,B})_{A,B \in \P_{n,2}}$, which is the solution matrix of the following matrix equation. Let $t$ be a formal variable. Then we consider matrices $P^m,\Lambda^m$ of unknowns such that,
\begin{equation}\label{greengen}
\left\{ \begin{array}{l} \Lambda^m_{A,B}=0  \mbox{ unless } A \sim_m B \\ P^m_{A,B}=0  \mbox{ unless either }A > B, A \nsim_m B, \mbox{ or } A=B \\ P^m_{A,A}=t^{a_m(\chi)}, \chi=\chi_A \\ P^m \Lambda^m \ ^tP^m=\Omega.\end{array} \right.
\end{equation}
Shoji shows that this equation has unique solution matrices
$P^m,\Lambda^m$. Every entry can be shown to belong to $\Q(t)$ but
Shoji conjectures that actually $P_{A,B} \in \Z[t]$. Notice that
for $m=\frac{1}{2}$ resp. $m=1$, one recovers the (indeed
polynomial) Green functions of $Sp_{2n}(\F_q)$ resp.
$SO_{2n+1}(\F_q)$ for char($\F_q$)$\neq 2$: one has
$P^m_{A,B}=\pi_{ji}$ (for the appropriate $\pi_{ji}$) if $A
\leftrightarrow i, B \leftrightarrow j$. Therefore the $P^m_{A,B}$ are in general also called Green functions. On the other hand, they are a generalization of Kostka functions $K_{\l,\mu}(t)$ in the sense that they form the transition matrix between (appropriately defined, cf. \cite{shoji}) Schur functions and Hall-Littlewood polynomials.

Our objective is to use these Green functions to describe the $W_0$-structure of every module in $\Hgrrcc$. However, it is sufficient to determine $\Hgrrcc$ for a fixed quadrant of the parameter space $\{(k_1,k_2)\mid k_1,k_2 \in \R\}$:

\begin{rem}\label{pospar} Let $J \subset I$ be such that $W_0(R_J)=R_J$. Then the map 
\[f_J: s_i \mapsto \begin{cases} -s_i & \mbox{if }i \in J \\ s_i & \mbox{if }i \notin J \end{cases}\] extends to an automorphism of $\C[W_0]$. We also define $f_J(k_{\a_i})=\pm k_{\a_i}$ (with the sign equal to the sign of $f_J(s_i)$). Finally define, for $x \in X$, $f_J(x)=x$. Then it is easy to see that $f_J$ extends to an isomorphism of generic algebras $f_J:\Hgr(\mathcal{R}^{deg},\vet{k}) \to \Hgr(\mathcal{R}^{deg},f_J(\vet{k}))$. Specializing in both algebras $\vet{k}_\a \mapsto k_\a \in \R$, we notice that any tempered $\Hgr$-module $V$ is also a tempered $f_J(\Hgr)$-module. On the other hand, the $W_0$-type of $V$ as a $\Hgr$-module is equal to the $W_0$-type of $V$ as a $f_J(\Hgr)$-module, tensored with the one-dimensional $W_0$-representation $\C_J$, taking the simple reflections in $J$ to $-1$ and the others to $1$.

Therefore, it is sufficient to understand $\Hgrrcc$ for non-negative parameters in order to retrieve $\Hgrrcc$ for arbitrary parameters.
\end{rem}

From the above remark and remark \ref{scaling}, it follows that if $R_0$ has type $B_n$, we can restrict ourselves to determining $\Hgrrcc$ for every {\it ratio} $m$ of the parameters $k_2=mk_1$, and we can assume $k_1 \geq 0, k_2 \geq 0$.

In view of this remark, we will need to know later that for the Green functions we have
\begin{lemma}\label{mP}
Let $A=(\xi,\eta) \in \P_{n,2}$. If we write $\Phi(A)=A \otimes (n,-)= (\eta,\xi)$ (viewed as characters of $W_0$), then
\[ P^{-m}_{A,B}=P^{m}_{\Phi(A),\Phi(B)}. \]
\end{lemma}

\pf: This follows immediately from the definition of the orderings $\sim_{\pm m}$ and the matrix equation \eqref{greengen}. \qed

\section{Conjectures}\label{conjecture}
Recall that in the equal label cases, the set $\Hgrrcc$ is parametrized by $\hat{W}_0$. We transfer this parametrization to $\P_{n,2}$ by writing $M^m_A$ for $M_\chi$ if $\chi=\chi_A$, where $m=1$ if $k_1=k_2$ (the $B_n$-case) and $m=\frac{1}{2}$ if $k_2=\half k_1$ (the $C_n$-case).  As a $W_0$-module, $M^m_A$ is naturally graded, and we write $M_A^{m;l}$ for is degree-$l$ part. Fix $j=(u,\rho) \in \I_0$ and let
\[ H^{2l}(\B_u)_\rho\simeq\sum_{\chi \in \hat{W_0}}n^l_{j,\chi} V_\chi,\]
for $i=0,1,\dots,d_u$ (it is known that all odd cohomology vanishes).

Suppose that we write the Green polynomial $\pi_{ji}$ as
\begin{equation}
\label{expansion}
\pi_{ji}(q)=\sum_{l \geq 0}\pi_{ji}^l q^l.
\end{equation}
Then clearly, if $\chi \leftrightarrow i$, it follows that \[ \pi_{ji}^l=n^l_{j,\chi},\]
i.e., if we write $P^m_{B,A}$ as $P^m_{B,A}(q)=\sum_{l\geq 0}P^{m;l}_{B,A}q^l$, then
\[ M_A^{m;l}=\sum_B P^{m;l}_{B,A} V_B \otimes \e. \]
We can now state our conjecture.
\begin{conj}\label{conj}
Let $k_2=mk_1 \in \R$ with $m\in \frac{1}{2}{\mathbb \Z}$. Let $\Hgr$ be the graded Hecke algebra associated to a root system of type $B_n$, whose root labels are $k_1,k_2$. Then we conjecture that

\begin{itemize}
\item[(i)]{the Green functions $P^m_{A,B}$ do not depend on the chosen refinement of the pre-order defined by $a_m$ and $\sim_m$. They are polynomials of degree $\leq a_m(A)$, i.e. may be written as $P^m_{A,B}(q)=\sum_{l \geq 0}P^{m;l}_{A,B}q^l$; }
\item[(ii)]{We have a bijection \[ \Hgrrcc \longleftrightarrow \hat{W_0},\] written as $M^m_A \leftrightarrow A \in \P_{n,2}$, which is uniquely determined by requiring that $\chi_A$ (resp. $\chi_A \otimes \e$) occurs in $M^m_A$ if $k_1<0$ (resp. $k_1>0$). }
\item[(iii)]{The $W_0$-modules $M^m_A$ can be viewed as graded modules, via a grading inherited from the Green functions $P_{BA}^m$, as follows.  We write its degree-$l$ part as $M^{m;l}_A$. It is given by the formula \[ M^{m;l}_A\simeq\begin{cases} \sum_B P^{m;l}_{B,A}V_B \otimes \e & \mbox{ if } k_1>0 \\ \sum_B P^{m;l}_{B,A} V_B & \mbox{ if } k_1<0.\end{cases} \]}
\item[(iv)]{$M^m_A$ and $M^m_B$ have the same central character (see \ref{con2} below) if and only if $A \sim_m B$.}
\end{itemize}
\end{conj}

\begin{cor}
The bijection in (ii) is obtained by taking the top degrees of the
$M_A^m$, i.e. we have $M^{m;max}_A=\chi_A$ (resp. $\chi_A \otimes
\e$) if $k_1<0$ (resp. $k_1>0$).
\end{cor}
\pf: This follows from the defining relations of the Green
functions. In particular, the top non-vanishing degree of $M^m_A$
is $a_m(A)$.\qed

\begin{cor}\label{mult1} Every irreducible tempered representation of $\Hgr$ with real central character has a multiplicity free $W_0$-type.
\end{cor}
\pf: We show that, granted the Conjecture, the $W_0$-module $V_A$ occurs exactly once in $M_A^{m}$ if $k_1<0$. We have
\[ (V_A,{M_A^m}_{|W_0})=(V_A,\sum_{l \geq 0}P_{BA}^{m;l}V_B)=\sum_{l \geq 0}P_{AA}^{m;l}=1;\]
and similarly for $V_A \otimes \e$ if $k_1>0$.\qed

Although there may be several $W_0$-characters occurring with multiplicity one, there is thus a unique such character $\chi_V$ for every $V \in \Hgrrcc$ if we demand that $V \mapsto \chi_V$ is injective. Then $\chi_V$ forms the top degree of $V_{|W_0}$.

\subsubsection{Special case: $A=(-,1^n)$}\label{sign}In the equal label cases, the $W_0$-module $H(\B_u)$ is for $u=1$ isomorphic to the coinvariant algebra (with degrees doubled). The only Springer correspondent of $C=\{1\}$ is the sign representation, indexed by $(-,1^n)$. For arbitrary $m$, the conjecture shows that the module $M_{(-,1^n)}^m$ is isomorphic to the coinvariant algebra as well, and the $W_0$-grading has undergone a shift:
\begin{lemma}
Let $0 \leq m \in \half \Z$ . Then, with $\e=(-,1^n)$ the sign representation, we have $P^{m}_{A,\e}(t)=t^{n(m-1)}P^{1}_{A,\e}(t)$ if $m \in \Z$, and $P^{m}_{A,\e}(t)=t^{n(m-\half)}P^{\half}_{A,\e}(t)$ if $m \notin \Z$.
\end{lemma}

\pf: Suppose that $m\in\Z$. One checks that the $a_m$-value is maximal among $\P_{n,2}$. Therefore one can show, analogous to Lemma 7.2 in \cite{shoji} that $t^{-a_m(\e)}R(\e)P^m_{A,\e}=R(A)$ (where we identify $\hat{W}_0$ with $\P_{n,2}$ and $R$ is as in section 3). On the other hand one checks with an easy computation that $a_m(\e)=n^2+n(m-1)$. The proof for $m \notin \Z$ is analogous.\qed
\subsubsection{Consistency} Notice that in view of Remark \ref{pospar}, Conjecture \ref{conj} leads to two, potentially different expressions of the modules
$M^m_A$: we can either use the Green functions
for $m$ or for $-m$, and multiply with the appropriate character.

One checks easily that, in order for these two results to be the same, we need that $M_{\Phi(A)}^{m,(k_1,k_2)}\mid_{W_0}=\Phi(M_{A}^{-m,(-k_1,k_2)}\mid_{W_0})\otimes\e$. However, this follows from Lemma \ref{mP}. Therefore there is
no ambiguity.

If $m=0$ then the Conjecture implies that  $\Phi(M_A^{0;l})=M_{\Phi(A)}^{0;l}$ for all $l$, hence, if $A=\Phi(A)$, the $W_0$-module $M_A^0$ is
invariant for $\Phi$ in every graded part. The restriction to the graded Hecke algebra of type $D$ of the modules $M^0_A$ is discussed in \cite{proefschrift}.


Notice also that we do not necessarily assume that the parameters are special. Indeed, it is easy to see that the all equivalence classes in $\hat{W}_0$ under $\sim_m$ are singletons if $|m|>n-1$. This means that for such parameters, the members of $\Hgrrcc$ are separated by their central character. Thus, we recover the special parameters \eqref{special} computed in \cite{HO}.


\section{Combinatorial generalization of the Springer correspondence}\label{combi}

We have not yet motivated our definition of $m$-symbols for general $m \in \half\Z$. In this section, we will define combinatorics from which they arise naturally. Indeed, these $m$-symbols describe a combinatorial generalization of the Springer correspondence. We use the explicit knowledge of the central characters of $\Hgrrcc$ to define, for $k_2=mk_1$, a set of Springer correspondents $\Sigma_m(W_0c_L)$ attached to a central character $W_0c_L$ of $\Hgrrcc$. Then we define a set $\U_m(n)$ of partitions (replacing the unipotent conjugacy classes), such that $\U_m(n)$ is in bijection with the set of central characters of $\Hgrrcc$. We can then transfer the Springer correspondents of $W_0c_L$ to its associated ``unipotent class''. This leads us to reformulate our conjecture into the form \ref{finalconj}.


\subsubsection{}First we establish some notation. For $m \in
\half\Z_{\geq 0}$, let $\ceil{m}$ be the smallest integer at least
equal to $m$ and let $\floor{m}$ be the biggest integer at most
equal to $m$. We write $\l \vdash n$ to denote that $\l$ is a
partition of $n$.

\subsection{The $m$-tableau} 
If $R_0$ is of type $B_n$ then a standard parabolic root subsystem $R_L\subset R_0$ contains irreducible factors of type $A$ and at most one of type $B$. Therefore, we first review the classification of the residual points for these types. This has been done in \cite{HO}. First consider type $A_{n-1}$. Let $V$ be an $n$-dimensional real vactor space with root system $R_0=\{e_i-e_j;i,j=1,\dots,n\}$. Let $k_\a=k$ for all $\a$. Then the only residual point in $V$ has coordinates (up to Weyl group action)
\begin{equation}\label{acenter} (-\frac{n-1}{2},-\frac{n-3}{2},\dots,\frac{n-3}{2},\frac{n-1}{2})k.\end{equation}
Now consider a root system of type $B_n$ with labels $k_1,k_2$ as in section \ref{spec}. Then, according to \cite[Prop. 4.3]{HO}, the $W_0$-orbits of residual points in $V$ are in bijection with set $\P_n$ of partitions of $n$. Let $\l \in \P_n$. We consider its Young diagram $T_\l$. The content of a box $\square$ with coordinates $(i,j)$ is defined to be $c(\square)=i-j$. An example of a Young tableau where each box is filled with its content is given in Figure \ref{content} below.

\begin{figure}[htb]
\begin{center}
{\includegraphics[angle=0,scale=0.4]{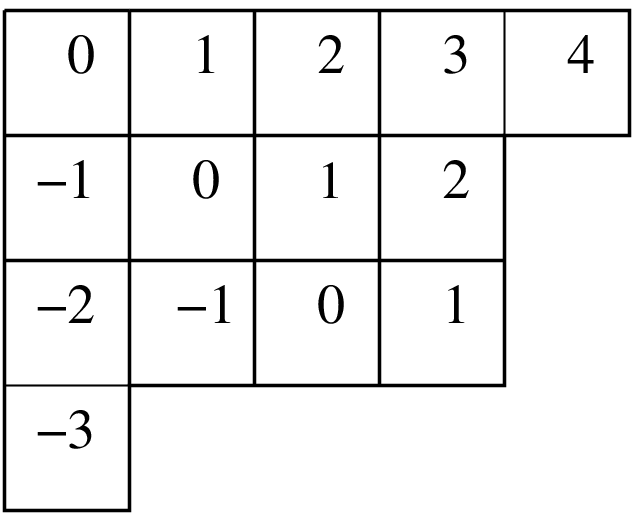}}
\end{center}
\caption{Content tableau of $(1445)$} \label{content}
\end{figure}

\begin{figure}[htb]
\begin{center}
{\includegraphics[angle=0,scale=0.4]{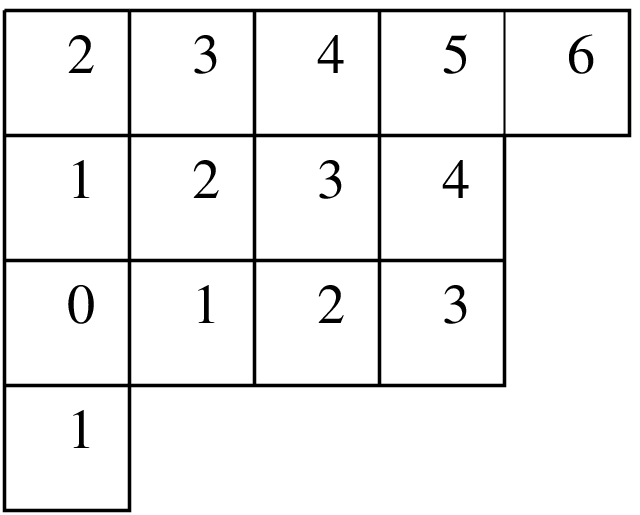}}
\end{center}
\caption{$T_2(\l)$ for $\l=(1445)$.}
\label{fig:qtableau}
\end{figure}

Let $c(\l,k_1,k_2)$ be the point in $V$ with coordinates $c(\square)k_1+k_2$ where $\square$ runs over the boxes of $T_\l$.
Then, by \cite{HO}, for generic parameters the points $c(\l,k_1,k_2)$ where $\l$ runs through $\P_n$ form a complete set of representatives for the $W_0$-orbits of residual points of $\Hgr$. At special parameters, one may have $c(\l,k_1,k_2)=c(\mu,k_1,k_2)$ even if $\l\neq \mu$, and the $c(\l,k_1,k_2)$  need not be residual any more. Let $k_2=mk_1$ with $m \in \halfZ$. Since $W_0$ acts by permutations and sign changes, we may write $c(\l,k,mk)$ as
 as
\begin{equation}\label{bnpunt}
\begin{array}{llll} (\underbrace{pk,\dots,pk,} & (p-1)k,\dots,k, & \underbrace{0,\dots,0}) & \in \R^n \\ M_p {\rm \ times} & & M_0 {\rm \ times}& \end{array}
\end{equation}
if $m$ is integer, or as
\begin{equation}\label{cnpunt}
\begin{array}{llll} (\underbrace{pk,\dots,pk,} & (p-1)k,\dots,k, & \underbrace{\frac{1}{2}k,\dots,\frac{1}{2}k}) & \in \R^n \\ M_p {\rm \ times} & & M_{\frac{1}{2}} {\rm \ times}& \end{array}
\end{equation}
if $m$ is half-integer. The conditions under which such a point remains residual differ between integer and half-integer $m$ and are investigated in \cite{HO}. We list the results here:
\begin{cond}\label{distvw} (cf. \cite{HO}) Consider $c(\l,k,mk)$ as in \eqref{bnpunt} or \eqref{cnpunt} for $k \neq 0$
and $m \in \half \Z$.
\begin{itemize}\item{If $m=0$ then the point \eqref{bnpunt} is
residual if and only if (i) $M_p=1$, (ii) $M_p \in
\{M_{p+1},M_{p+1}+1\}$ for $p>0$, and (iii) $M_0=\lfloor
\frac{M_1+1}{2}\rfloor$.}
\item{If $k_2=mk$, with $0 \neq m \in \Z$, then $p \geq |m|$ and the point
\eqref{bnpunt} is residual if (i) $M_p=1$, (ii) $M_l \in
\{M_{l+1},M_{l+1}+1\}$ for all $l \geq |m|$, (iii) $M_l \in
\{M_{l+1}-1,M_{l+1} \}$ for $l=1,2,\dots,|m|-1$, and finally (iv)
$M_0=\lfloor\frac{M_1}{2}\rfloor$};

\item{If $k_2=mk_1$ with $m \in \halfZ, m\notin \Z$, then $p \geq |m|$
and the point \eqref{cnpunt} is residual if (i) $M_p=1$, (ii) $M_l
\in \{M_{l+1},M_{l+1}+1\}$ for all $l \geq |m|$, (iii) $M_l \in
\{M_{l+1}-1,M_{l+1} \}$ for
$l=\frac{1}{2},\frac{3}{2},\dots,|m|-1$.}
\end{itemize}
\end{cond}
We will now translate this algebraic condition to a condition on
Young tableaux. To this end we define
\begin{defn} For $m \in \halfZ$ and $\l \vdash n$, let $T_m(\l)$ be the Young tableau of $\l$ where box $\square$ is filled with $|c(\square)+m|$. We call it the $m$-tableau of $\l$.
\end{defn}
An example is given in Figure \ref{fig:qtableau}. The
entries of $T_m(\l)$ arranged in
non-increasing order, multiplied with $k$, form the vector $c(\l,k,mk)$ in the form \eqref{bnpunt} or \eqref{cnpunt}.

\begin{defn} Let $\l \vdash n, m \in \halfZ$.
We define the {\it extremities} of $T_m(\l)$ to be the numbers in the
following list. If $m \in \Z$  (resp. $m \notin \Z$), the extremities are the entries of every last square of the rows of $T_m(\l)$ which lies on or above the zero-diagonal (resp. the upper $\half$-diagonal), and the entries of the every last square of the columns of $T_m(\l)$ which lies on or below the zero-diagonal (resp. the lower $\half$-diagonal). 
\end{defn}

Notice that it may happen that a number occurs twice as extremity. In particulat, zero can only occur twice as extremity.

On the other hand, let $\l$ be such that $c(\l,k,mk)$ is residual
and write it in the form \eqref{bnpunt} or \eqref{cnpunt}. Then
$c$ is clearly determined by the following numbers:
\begin{defn}\label{jumps}
 Suppose $c:=c(\l,k,mk)$ is a residual point written in the form \eqref{bnpunt} or \eqref{cnpunt}.
We define the \emph{jumps} of $c$ to be those $l \geq |m|$ such
that $M_l=M_{l+1}+1$, and those $l \in \{1,2,\dots,|m|-1\}$ such
that $M_l=M_{l+1}$. It will be convenient to have $\ceil{|m|}+2r$
jumps for some $r \in \Z_{\geq 0}$. Therefore, we add $0$ (for
integer $m$) or $-\frac{1}{2}$ (for half-integer $m$) as a jump if
necessary.

We write $J(c)$ for the set of jumps of $c$, and $J^+(c)$ for the
set of positive jumps of $c$. We write the jumps in increasing
order: $j_1<j_2\dots<j_{\ceil{|m|}+2r}$.
\end{defn}

For this definition to make sense, it remains to be seen that
there are enough positive jumps.
\begin{lemma}
\label{aantaljumps} Let $c=c(\l,k,mk)$ be a residual point. Then
$|J^+(c)|=M_1+|m|-1$ if $|m|\geq 1$ is integer, $|J^+(c)|= M_\frac{1}{2}+|m|-\frac{1}{2}$ if $m$ is half-integer, and
$|J^+(c)|=M_1$ if $m=0$.
\end{lemma}

\pf:  (i) Let $m=0$. Then $j \in J^+(c)$ gives rise to a sequence
$\{j,j-1,\dots,1\}$ in \eqref{bnpunt}.

(ii) If $m \in \Z$ then $|m| \leq j \in J(c)$ corresponds to a
subsequence $\{j,j-1,\dots,|m|\}k$ in \eqref{bnpunt}. Therefore
$M_{|m|}=|\{j \in J(c) \mid j \geq |m|\}|$. Since $l < m$ is a
jump if and only if $M_{l+1}=M_l$, it follows that $M_{|m|}-M_1=
|\{i\in\{1,2,\dots,m-1\}\mid i \notin J(c)\}|$, i.e. $|J^+(c)|=|\{j \in J(c)\mid j \geq |m|\}|+ | \{1 \leq j \leq
|m|-1 \mid j \in J(c)\}|=M_{|m|}+(|m|+1-M_{|m|}+M_1)=|m|+1+M_1$.

(iii) If $m \notin \Z$ then $|m| \leq j \in J(c)$ corresponds to a
subsequence $\{j,j-1,\dots,|m|\}k$ in \eqref{cnpunt}. Therefore
$M_{|m|}=|\{j \in J(c) \mid j \geq |m|\}|$. Since $l < m$ is a
jump if and only if $M_{l+1}=M_l$, it follows that
$M_{|m|}-M_{1/2}= |\{i\in\{1/2,3/2,\dots,m-1\}\mid i \notin
J(c)\}|$, i.e. $| J^+(c)|=|\{j \in J(c)\mid j \geq |m|\}|+
| \{1/2 \leq j \leq |m|-1 \mid j \in
J(c)\}|=M_{|m|}+(|m|+1/2-M_{|m|}+M_{1/2})=|m|+1/2+M_{1/2}$. \qed

Therefore we can indeed, by defining $0$ or $-\half$ to be a jump, arrange to have $|m|+2r$ or $|m|+2r+\frac{1}{2}$ jumps. Note that if $m \neq 0$, then $0$ (resp. $-\frac{1}{2}$) occurs as jump if
and only if $M_1$ (resp. $M_\frac{1}{2}$) is even.

\begin{lemma}\label{ml} Let $m \in \halfZ$ and $\l \vdash n$. Let $M_l$ be the multiplicity of $l$ in $T_m(\l)$. Then we have
\begin{itemize}
\item[(i)] $M_p\in \{M_{p+1},M_{p+1}+1,M_{p+1}+2\}$ for $p \geq |m|$;
\item[(ii)] $M_p \in \{M_{p+1}-1,M_{p+1},M_{p+1}+1\} $ for $p=1,2,\dots,|m|-1$ resp. $\half,\frac{3}{2},\dots,|m|-1$;
\item[(iii)] If $0 \neq m \in \Z$, then $M_0 \in \{ \lfloor\frac{M_1}{2}\rfloor,\lceil\frac{M_1}{2}\rceil\}$;
\item[(iv)] If $m=0$, then $M_0 \in \{ \lfloor\frac{M_1+1}{2} \rfloor
\}, \lceil \frac{M_1+1}{2} \rceil\}$.
\end{itemize}
\end{lemma}

\pf: Suppose that the Lemma holds for $m\geq 0$. Then it also
holds for $m<0$ since $T_{-m}(\l)=T_{m}(\l')'$. Therefore we need
only prove the case $m \geq 0$.

(i): This follows because $p$ and $p+1$ occur in two diagonals, an
upper and a lower one:

\begin{center}
{\includegraphics[angle=0,scale=0.4]{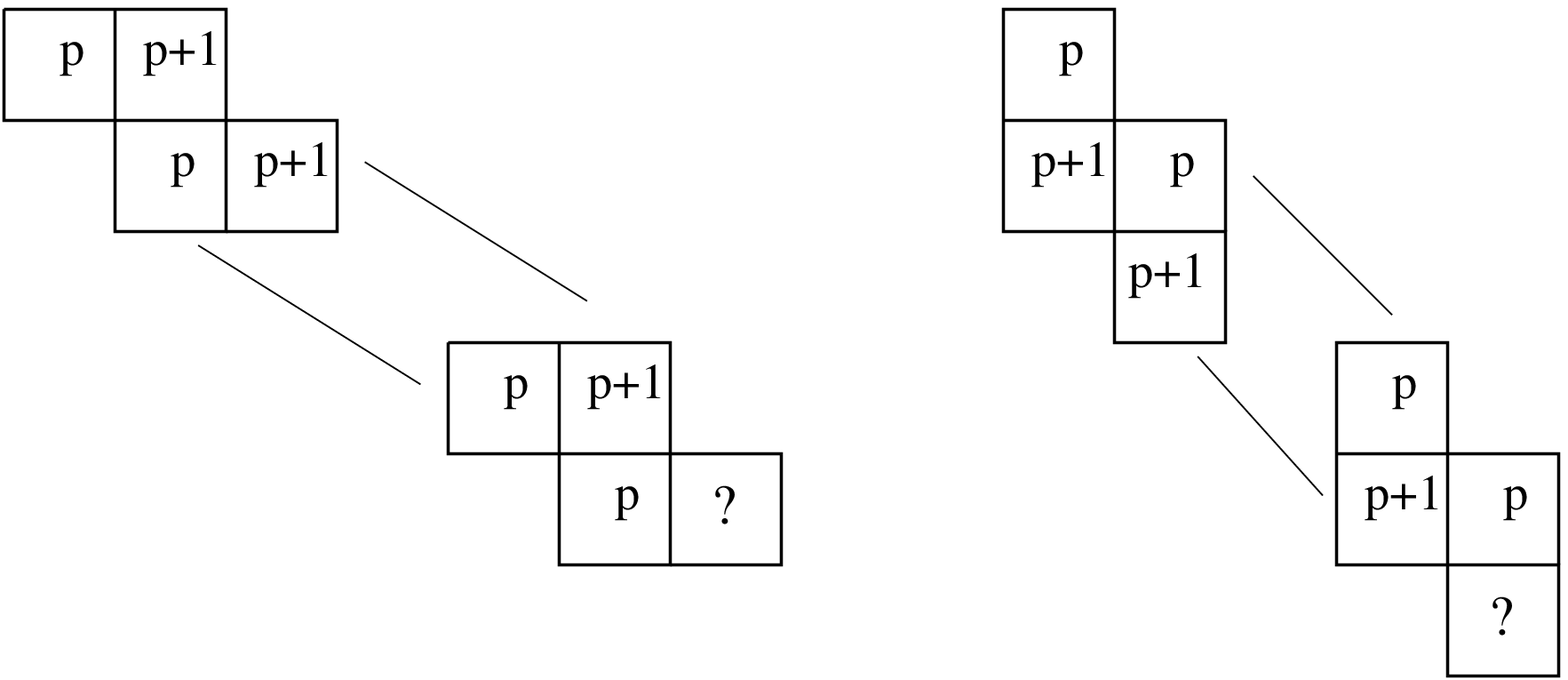}}
\end{center}

In the squares marked with $?$, we may or may not have a $p+1$.

(ii) Here, things change for the upper diagonal as we always have a $p+1$ above the first $p$:

\begin{center}
{\includegraphics[angle=0,scale=0.4]{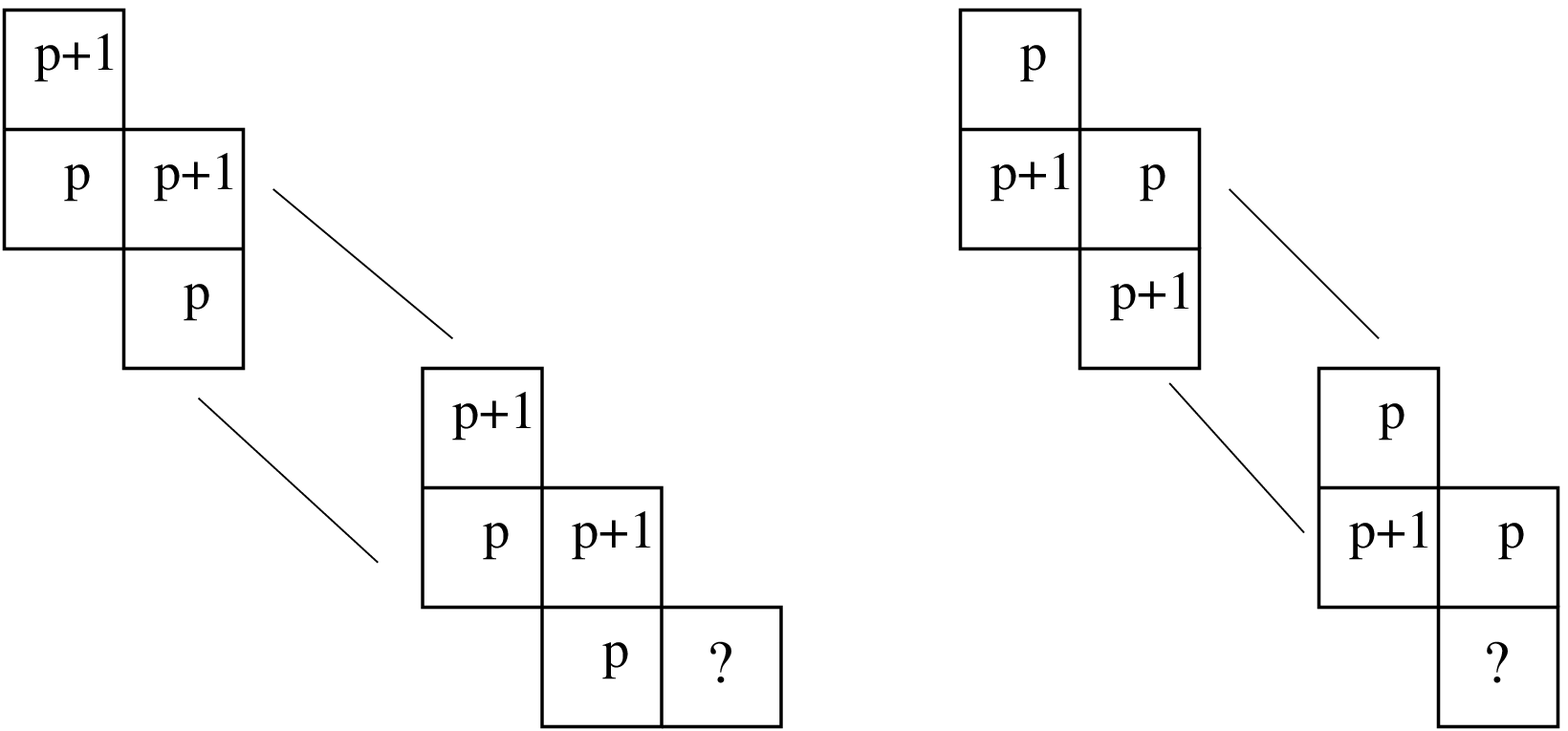}}
\end{center}

(iii) This follows from the fact that the zeroes occur as indicated:
\begin{center}
{\includegraphics[angle=0,scale=0.4]{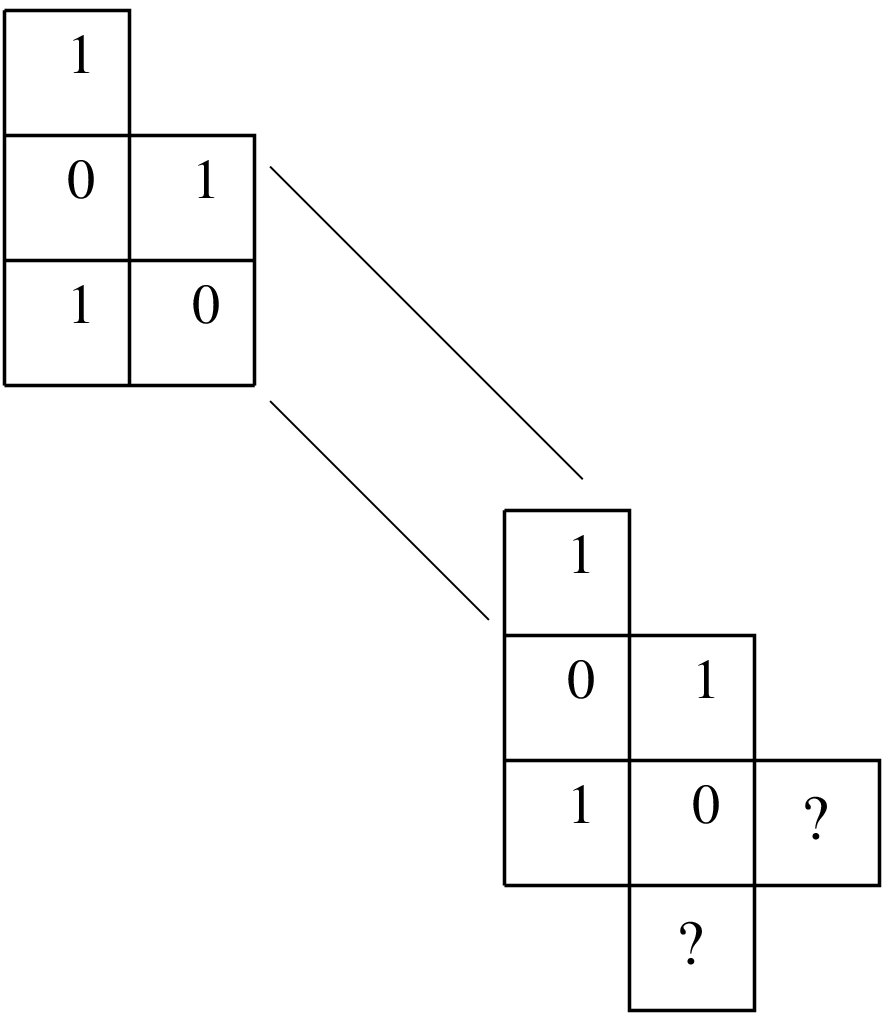}}
\end{center}
(iv) As in case (iii), but the top square containing a one is not
present now.\qed

\begin{cor}\label{extr} Let $\l \vdash n$ such that $c(\l,k,mk)$ is residual. Then for all $x$ occurring as entry in $T_m(\l)$ we have:
$x$ is a jump $\iff$ $x$ is an extremity.
\end{cor}

\pf: (i) Suppose $x>0$. Since $c(\l,k,mk)$ is residual, the
sequence covered by $T_m(\l)$ satisfies Condition \ref{distvw},
and therefore in the figures of the proof of parts (i) and (ii) of
Lemma \ref{ml}, at least one of the squares containing a ? are
occupied by $T_m(\l)$. But then we have that $x$ is a jump $\iff$
only one square containing a ? is occupied $\iff$ there is an
extremity square having entry $x$.

(ii) Suppose $x=0$. Then there is an extremity square with entry
zero $\iff$ only one square filled with ? in the proof of
Proposition \ref{ml} (iii) is occupied $\iff$ $m \neq 0$ and $M_1$
is even, or $m=0$ and $M_1$ is odd. \qed

\begin{cor}\label{distext} Let $\l \vdash n$ and $m\in \halfZ$. Then $c(\l,k,mk)$ is residual if and only if $T_m(\l)$ has distinct extremities.
\end{cor}

\pf: This follows directly from the Lemma and Corollary.\qed

\noindent It follows again that our notion of special parameters coincides with the one in \cite{HO}:
\begin{cor}\label{gendistinct} Let $m \in \R$. The points $c(\l,k,mk)$ are all residual and distinct if and only if $m \notin \pm\{0,\frac{1}{2},1,\dots,n-\frac{3}{2},n-1\}$.
\end{cor}

\pf: In view of the preceding Corollary, it suffices to check for
which $m$ every $T_m(\l)$ has distinct extremeties. Indeed, it is easy to see that if there are $\l,\mu \vdash n$ such that $c(\l,k,mk)$ and $c(\mu,k,mk)$ are both residual and $c(\l,k,mk)=c(\mu,k,mk)$, then there is also a $\nu \vdash n$ such that $c(\nu,k,mk)$ is not residual. Furthermore,
it suffices to check the Corollary for non-negative $m$ because of
symmetry. Let $m \in \R_{\geq 0}$ and write $m=l+\epsilon$, with
$l \in \mathbb{Z}_{\geq 0}$ and $0 \leq \epsilon <1$. Then
$T_m(\l)$ is drawn on box paper filled as follows:
\begin{center}
\includegraphics[scale=0.4]{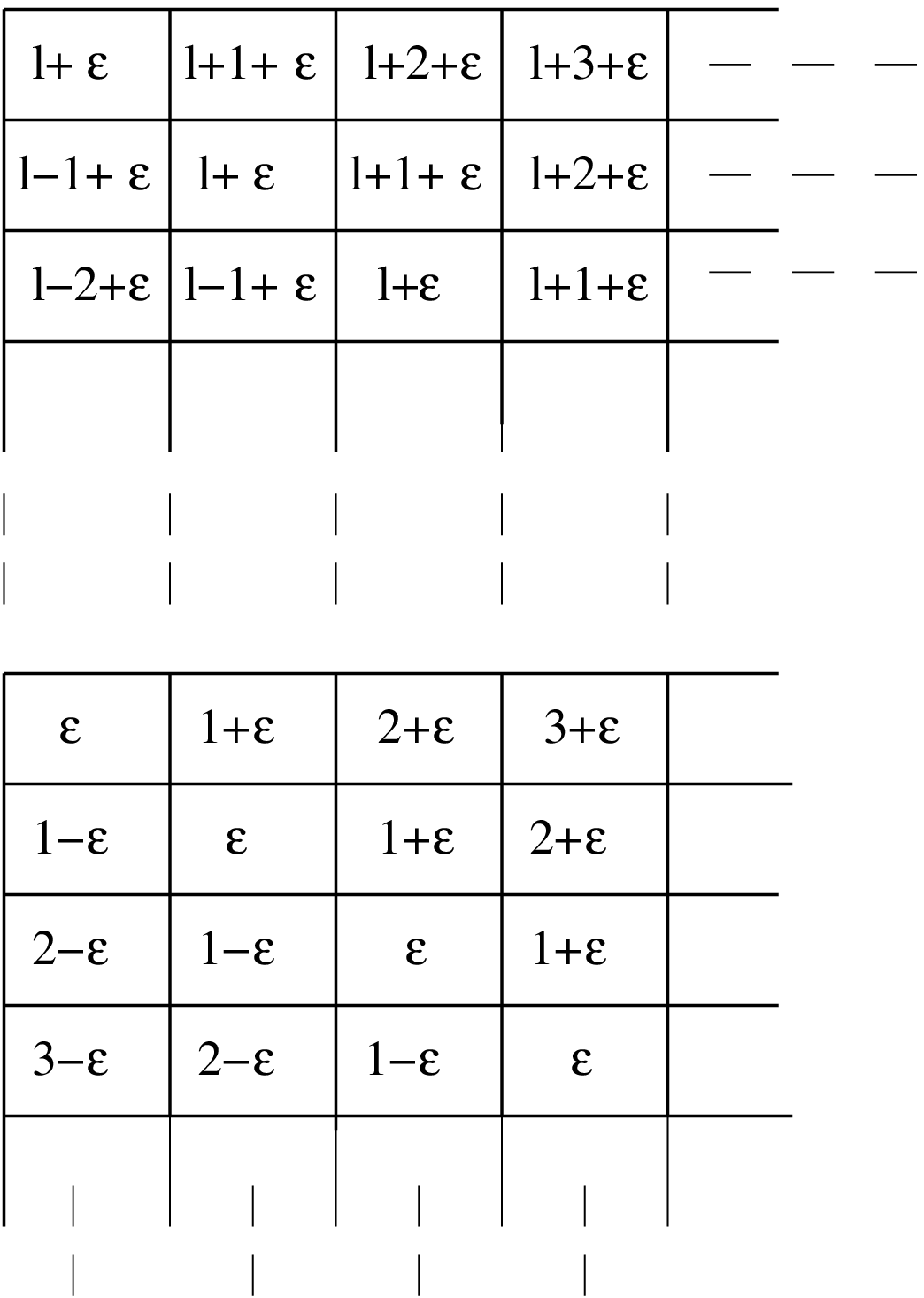}
\end{center}
There are no repetitions in
$l+f+\e,l+(f-1)+\e,\dots,1+\e,\e,1-\e,2-\e,\dots$, unless
$l+f+\e=d-\e$ for some $f,d$. This means that $\e\in\{0,\half\}$. Clearly, we can reduce to the case where the multiple extremity is either 0 or $\half$. But such an extremity can occur precisely for $m\in \pm\{0,\frac{1}{2},1,\dots,n-\frac{3}{2},n-1\}$.\qed

\subsection{Specialization of residual points} In this section we show that
 a generically residual
point specializes for all parameters values into the central character of a
tempered representation.

\begin{defn} Let $\Hgr$ have parameters $k_2=mk_1$.
Let $\L_m(n)$ denote the set of $W_0$-orbits of residual subspaces of $\Hgr$. We say that $W_0L \in \L_m(n)$ has type $(\l,\mu) \in
\P_{n,2}$, if $W_Lc_L=W_Lc(\mu,k,mk) \subset \t_L$ and $R_L=A_\l \times B_{|\mu|}$. Notice that for special parameters, $\mu$ need not be uniquely determined by $L$ since several residual points $c(\nu,k_1,k_2)$ may coincide with each other for $k_2=mk_1$.
We denote by $\centra=\{W_0c_L \mid W_0L \in \L_m(n) \}$ the set of central characters of $\Hgrrcc$.
\end{defn}

Sometimes we will abuse notation and write $L \in \L_m(n)$ instead of $W_0L \in \L_m(n)$, if it is clear that we are only interested in $L$ up to $W_0$-orbit.

\begin{prop}\label{altijdspec} Let $\l\vdash n$ and let $k_2=mk_1$ be special. Then there exists $W_0L \in \L_m(n)$ such that \[ W_0c(\l,k,mk)=W_0c_L \in \centra. \] Let $L$ have type $(\mu,\nu)$. Then $\mu$ consists of the lenghts of the hooks in $T_m(\l)$ ending on equal extremities, and $\nu$ is obtained by deleting these hooks from $\l$. In particular, the parts of $\mu$ are all different. If $m$ is integer they are odd, otherwise they are even. We write $sp_m(\l)=(\mu,\nu)$.
\end{prop}

\pf: This boils down to showing that the entries of $T_m(\l)$ can
be written as the union of the entries of a residual tableau
$T_m(\mu)$ for some partition $\mu$ whose Young diagram is
contained in that of $\l$, and a number of sequences of the form
$(l,l,l-1,l-1,\dots,1,1,0)$ (resp.
$(l,l,l-1,l-1,\dots,\frac{1}{2},\frac{1}{2})$) if $m$ is integer
(resp. not integer), with $l \geq 0$. We assume that $m$ is
integer, the other case being analogous. The entries of $T_m(\l)$
form the vector
\[ \begin{array}{llll} (\underbrace{p,\dots,p,} & (p-1),\dots,1, & \underbrace{0,\dots,0}) & \in \R^n \\ M_p {\rm \ times} & & M_0 {\rm \ times}& \end{array} \]
We want to prove the proposition by induction on $n$. If $n=1$, the assertion is trivial. Now suppose $n>1$, and the proposition holds for all partitions of weight less than $n$. Take the maximal $l$, for which $m_l$ and $m_{l+1}$ do not satisfy Condition \ref{distvw}. There are three cases:

(i) We have $l \geq m$, and $M_l=M_{l+1}+2$. Then both squares
with a $?$ in Lemma \ref{ml} are not occupied, and this means that
there are both a row and column ending on $l$. Together, they
contain the sequence $(l,l,l-1,l-1,\dots,1,1,0)$ if $m \in \Z$,
and $(l,l,l-1,l-1,\dots,\half,\half)$ otherwise. We remove this
row and column, and push the disconnected piece one block
northwest to reconnect it. This does not change the content of any
box, and we find the $m$-tableau of a smaller partition. The claim
follows by induction.

(ii) We have $l \in \{1,2,\dots,m-1\}$ and $M_l=M_{l+1}+1$. This means that both squares containing $?$ in part (ii) of lemma \ref{ml} are not occupied. Again, this means there is both a row and a column ending on $l$, and we can remove them as in the case (i).

(iii) If $m \in \Z, m \neq 0$ (resp. $m=0$) and
$M_0=\lceil\frac{M_1}{2}\rceil$ (resp.
$M_0=\ceil{\frac{M_1+1}{2}}$), then both squares containing a $?$
in case (iii) (resp. (iv)) of Lemma \ref{ml} are not occupied. We
can therefore remove the last square containing a zero.\qed

The following fact allows us to focus on the centers of the residual subspaces, i.e. on the central characters of $\Hgrrcc$, instead of on the residual subspaces.

\begin{prop}\label{bijektie}
For all parameters $k_1,k_2$, the map
\[ \L_m(n) \to \centra:W_0L \to W_0c_L\]
is a bijection.
\end{prop}
\pf: This fact is known to hold for the equal label case (for any root system). Therefore, for generic $k_1,k_2$, the statement holds due this fact and Corollary \ref{gendistinct}. Next let $k_2=mk_1$ with $m \in \half \Z$. If $k_1=0$, then all orbits $W_0c(\l,0,k_2)$ coincide with each other into $W_0(k_2,\dots,k_2)$, which is not residual. A generically residual point in $\Hgr_L$ with $R_L=A_\l \times B_l$ and $l \neq 0$ therefore does not specialize into a residual point. If $R_L=A_\l$ then its center specializes into $c_L=0$, which is also not residual. Thus, $\L_m(n)$ is empty and $\centra$ as well.

Finally suppose that $k_2=mk_1 \neq 0$ with $m\in \half\Z$. Put $k=k_1$. We suppose for simplicity that $m \in \Z$, the other case being analogous. Suppose that $W_0L_1,W_0L_2$ are such that $W_0c_{L_1}=W_0c_{L_2}=:W_0c$. We choose $c$ in its $W_0$-orbit to take the form \eqref{bnpunt}. Let $L_i$ have type $(\l^{(i)},\mu^{(i)})$. We may assume that $\l^{(1)}$ and $\l^{(2)}$ have no parts in common. Therefore, we may also suppose that all the entries of $c_{L_i}$ are in $\Z_{\geq 0}k$. Indeed, otherwise we would find factors $A_{\l^{(i)}_j-1}$ with odd $\l^{(i)}_j$, but for type $A$, the required bijection is known to hold, hence we would find common parts in $\l^{(1)}$ and $\l^{(2)}$. 

The coordinates of the center $c_{L_i}$ of the residual subspace $L_i$ of type $(\l^{(i)},\mu^{(i)})$ can be written as the union of the lists $(\frac{\l^{(i)}_i-1}{2},\frac{\l^{(i)}_i-1}{2},\frac{\l^{(i)}_i-3}{2},\frac{\l^{(i)}_i-3}{2},\dots,1,1,0)$ and the coordinates of $c(\mu^{(i)},k,mk)$.

Suppose that for $p \in \Z_{\geq 0}$, the multiple $pk$ occurs $M_p$ times in $c_{L_1}=c_{L_2}$. Recall that $c_{L_i}$ is a residual point if and only if Condition \ref{distvw} holds for $M=(M_0,M_1,\dots,M_p)$. 
In this case, we are done. Otherwise, since $W_0c \in \centra$ it follows from Condition \ref{distvw} and the fact that every non-zero entry occurs twice in a series corresponding to an $A$-factor, that there is exactly one way to decompose $M$ into the coordinates of a residual point, together with series of the form $(l,l,\dots,1,1,0)$ corresponding to type $A$ factors. Therefore, if $W_0c_{L_1}=W_0c_{L_2}$ then $\l^{(1)}=\l^{(2)}$ and $c(\mu^{(1)},k,mk)=c(\mu^{(2)},k,mk)$. This implies that $W_0L_1=W_0L_2$. \qed

Remark this bijection only holds on the level of $W_0$-orbits, and not on the level of residual subspaces and their centers. Indeed, in general one may have $1 \neq w \in W_0$ such that $w(c_L)=c_L$ but $w(L) \neq L$.
 
We expect that for generic parameters, the representations in
$\Hgrrcc$ are separated by their central characters, and that
there are $|\hat{W}_0|$ such. 
Therefore we remark
\begin{lemma}\label{generiekW0} Let $k_1,k_2$ be generic and put $m=k_2/k_1$. Then $|\centra|=|\hat{W}_0|$.
\end{lemma}
\pf: Associate to $(\l,\mu) \in\P_{n,2}$ the residual subspace of
type $(\l,\mu)$. Since the parameters are generic, these are all distinct. Now use Proposition \ref{bijektie}. \qed

It remains to be shown that generically, for every
residual point $c$, there is exactly one irreducible discrete
series representation of $\Hgr$ with central character $W_0c$. For
regular $c$, this is not hard (see \cite{proefschrift}), but for
singular $c$, it remains an open problem.

\subsection{Springer correspondence, combinatorially} We recall the explicit combinatorial description of the Springer correspondence in the equal label cases.
\subsubsection{Unipotent classes}\label{unipotentclasses} The set of unipotent conjugacy classes of $SO_{2n+1}(\C)$ is parametrized by the elementary divisors partitions
\[ \U_1(n)=\{ \l=(1^{r_1}2^{r_2}\dots) \vdash 2n+1\mid  r_i \mbox{ is even if }i \mbox{ is even}\}.\]
Similarly for $Sp_{2n}(\C)$, we have the set
\[ \U_\half(n)=\{ \l=(1^{r_1}2^{r_2}\dots) \vdash 2n\mid  r_i \mbox{ is even if }i \mbox{ is odd}\}.\]
The distinguished unipotent classes among these are the ones parametrized by partitions which consist of distinct parts.
\subsubsection{The map $\phi$}\label{phi} Lusztig has defined maps, which we denote by $\phi_m$, $ \U_m(n) \to \P_{n,2}$ such that $\Sigma_m(C_\l)=[\phi_m(\l)]_m$. These have the following form. Let $\l \in \U_m(n)$. If $m=\half$ we make sure that $\l$ has an even number of parts by calling the first part zero if necessary. Then we define a partition $\l^*$ by putting $\l^*_i=\l_i+(i-1)$. We define two partitions $\xi^*$ and $\eta^*$ by letting the odd parts of $\l^*$ be $2\xi_1^*+1<2\xi^*_2+1<\dots$ and by letting the even parts of $\l^*$ be $2\eta_1^*<2\eta_2^*<\dots$. Finally we define $\xi_i=\xi^*_i-(i-1)$ and $\eta_i=\eta_i^*-(i-1)$. Then by definition $\phi_m(\l)=(\xi,\eta)$.

\subsubsection{Bala--Carter map}\label{Balacarter}
We now describe the bijection $\U_m(n) \leftrightarrow \L_m(n)$
discussed in \ref{BC}, by giving explicitly the map $f_m^{BC}: \U_m(n) \to \L_m(n)$ implementing it (BC standing for Bala-Carter). It is defined as follows (cf. \cite{jantzen}): 

(1) If $C_\l$ is distinguished, i.e. $\l$ consists of
distinct parts, then let $c$ be the residual point associated to
$C_\l$ as described in \ref{groupcases}. Put $W_0c=f^{BC}_m(L)$.

(2) For $\l \in \U_m(n)$, let $l_1<l_2< \dots$ be
the parts that occur an odd number of times. Let their sum be
$n_0$. If we remove each part $l_i$ from $\l$ we obtain a
partition in which all parts have even multiplicity. Remove each
second part and let the remaining partition be $d_1 < d_2 <
\dots$. Then $C_\l$ corresponds to the residual subspace $L \in
\L_m(n)$ such that $R_L=A_{d_1-1} \times A_{d_2-2} \times \dots
\times B_l$, and whose center $c_L$ is the residual point in
$\t_L$ which corresponds to the distinguished unipotent class in
$SO_{n_0}(\C)$ (resp. $Sp_{n_0}(\C)$) with partition $(l_i)$, as decribed in (1). We put $W_0c_L:=f^{BC}_m(\l)$.

  We remark here that it follows 
from the definition of the weighted Dynkin diagram of $C_\l$ (cf. 
\cite[p. 395]{carter}) that if $C_\l$ is distinguished, the jumps of $\ga$ (as in (1)) are equal to $\frac{\l_i-1}{2}$.

Our generalization starts from the remark that we may transfer the Springer correspondents  from the unipotent classes $\U_m(n)$ to the central characters $\centra$ by defining
\begin{equation}\label{doel} \Sigma_m(f_m^{BC}(\l))=[\phi_m(\l)]_m.\end{equation}

\subsection{Generalized definition of Springer correspondents} \label{sigmarespunt} 
In this section we define a set $\Sigma_m(W_0c_L)$, where $m \in \half\Z$ and $W_0c_L \in \centra$. 
Let $k_2=mk_1$ with $m \in \halfZ$. Let $\l \vdash n$ and suppose that $c=c(\l,k,mk)$ is a residual point with jumps $j_1<j_2<\dots j_{\ceil{|m|}+2r}$
. Then we define, for non-negative $m$:
\begin{equation}
(\xi(c),\eta(c))=((j_1,j_3,\dots,j_{2r+1},j_{2r+2}-1,j_{2r+3}-2,\dots,j_{2r+m}-(m-1)),
\end{equation}
\[ (j_2+1,j_4+1,\dots,j_{2r}+1)) \]
if $m$ is integer, and
\begin{equation}
\label{sigmarespunt2}
(\xi(c),\eta(c))=((j_1+\half,j_3+\half,\dots,j_{2r+1}+\half,j_{2r+2}-\half,j_{2r+3}-\frac{3}{2},\dots,j_{2r+\ceil{m}}-(m-1)),
\end{equation}
\[ (j_2+\half,j_4+\half,\dots,j_{2r}+\half)) \]
if $m$ is not integer.

If $m < 0$, we put $(\xi(\tilde{c}),\eta(\tilde{c}))=(\eta(c),\xi(c))$ where $\tilde{c}$ is a residual point for $k_2=mk_1$ and $c$ is the residual point for $k_2=-mk_1$ with the same jumps as $\tilde{c}$. 
For the classical values, we have
\begin{lemma}
Let $m\in\{\half,1\}$ and suppose that $c$ is a residual point corresponding to the distinguished unipotent conjugacy class $C_\l$. Then $(\ga(c),\eta(c))=\phi_m(\l)$.
\end{lemma}
\pf: This boils down to showing that $\l_i=2j_i+1$, which was
already observed above.  \qed

\begin{defn} (i) We define the set of Springer correspondents of the residual point $W_0c \in \centra$ to be
\[ \Sigma_m(W_0c)=[(\ga(c),\eta(c))]_m.\]
(ii) Let $W_0L \in \L_m(n)$ with $R_L=A_\l \times B_l$. Then we define the set of Springer correspondents of $W_0c_L$ as
\begin{equation}
\Sigma_m(W_0c_L)=\bigcup_{(\a,\beta) \in \Sigma_m(W_Lc_L)}{\rm{tr}}_m-{\rm{Ind}}_{W_0(R_L)}^{W_0}(\rm{triv}_\l\otimes (\a,\beta)).
\end{equation}
\end{defn}

Note that, in the notation as above, $\Sigma_m(W_0c)=\Sigma_{-m}(\tilde{c}) \otimes (-,n)$ (where we identify 2-partitions with irreducible
characters, i.e., we switch the factors of the
members of $\Sigma_{-m}(\tilde{c})$ to obtain $\Sigma_m(W_0c)$. For
example, if $m=-1$ and the jumps of $c$ are ${0,1,2}$ then we have
$\Sigma_{-1}(c)=\{(2,02),(0,22),(4,00)\}$.

\subsubsection{Objective} In the equal label cases, the representations in $\Hgrrcc$ with central character $W_0c_L=f_m^{BC}(\l)$ are in bijection with the set $\Sigma_m(C_\l)$ of Springer correspondents of $C_\l$. We will show that the number of centers of generically residual cosets which specialize into $W_0c_L$ is equal to the numbers of Springer correspondents of $C_\l$. Therefore we define the sets, for $m \in \halfZ$ and $L \in \L_m(n)$,
\[ C_m(W_0c_L) =\]\[ \{ W_0c_{L'}\mid  L' \mbox{ generically residual subspace s.t. }  W_0c_{L'}=W_0c_L \mbox{ if }k_2=mk_1 \}.\]
Notice that, by \ref{generiekW0} and \ref{altijdspec}, we have
\[ \bigcup_{L \in \L_m(n)}|C_m(W_0c_L)|=|\hat{W_0}|.\]
We will show that $\Sigma_m(W_0c_L)$ is a similarity class for $\sim_m$
in $\hat{W}_0$ which is in bijection with $C_m(W_0c_L)$ (natural up
to at most involution), and that in the equal label cases we retrieve the already existing Springer correspondence, i.e., we have $\Sigma_m(W_0c_L)=\Sigma_m(C)$ if $C=C_\l$ and $W_0c_L=f^{BC}_m(\l)$.

\subsection{Splitting and joining} First we show that for a residual point $c$, we  have a canonical bijection $\Sigma_m(W_0c) \leftrightarrow C_m(W_0c)$. Let $\P_m(W_0c)=\{\l \vdash n \mid W_0c(\l,k,mk)=W_0c\}$. Of course, we may identify $\P_m(W_0c)$ with $C_m(W_0c)$.
We now define a map $\S_m:\P_m(W_0c)\to \P_{n,2}$ which implements the bijection.
\begin{defn} \label{defsplit}
 Suppose $\l\vdash n$. We divide $T_m(\l)$ into blocks as follows. Locate
the maximal entry of $T_m(\l)$, and draw a block from its square
up to the top left square. This is either a horizontal or a
vertical box. Subsequently, look for the maximal remaining entry
and enclose it in vertical or horizontal box, such that the formed
boxes for the diagram of a partition. This induces a splitting of
$T_m(\l)$ into horizontal and vertical boxes. A $1 \times 1$-block is considered as horizontal if it lies above the zero (resp. on or above the upper $\half$)-diagonal, and as vertical if it lies below the zero-diagonal (resp. on or below
the lower $\half$)-diagonal. Let $\xi$, resp. $\eta$,
be the partition whose parts are the lengths of the horizontal, resp. the vertical,
blocks. Then we define $\S_m(\l)=(\xi,\eta)$.
Examples are given in Figure \ref{fig:demosplits} below.
\begin{figure}[htb]
\begin{center}
{\includegraphics[angle=0,scale=0.4]{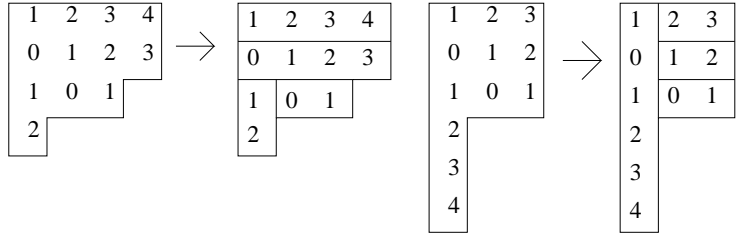}}
\end{center}
\caption{Examples: $\S_1(1344)=(244,2)$ and
$\S_1(1^33^3)=(222,6)$} \label{fig:demosplits}
\end{figure}

\end{defn}

\begin{lemma} The procedure $\S_m$ is well-defined on $\l \vdash n$ if and only if $c(\l,k,mk)$ is residual.
\end{lemma}
\pf: In view of Corollary \ref{distext} it remains to check that
for $m \in \Z$ we do not get $1\times 1$-boxes containing a zero. Such a
box can only arise as the last box under $\S_m$. But then in the
situation of Lemma \ref{ml}(iii)-(iv), the squares marked with ?
do not belong to $T_m(\l)$, which implies that $c(\l,k,mk)$ is not residual.
\qed

\begin{figure}[htb]
\begin{center}
{\includegraphics[angle=0,scale=0.4]{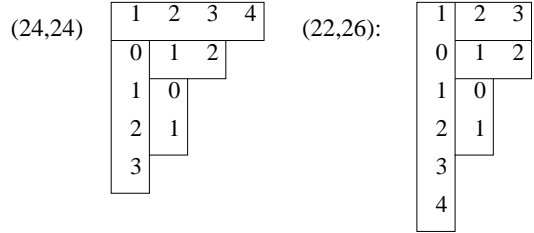}}
\end{center}
\caption{Examples: $\J_1(024,24)=(12234)$,
$\J_1(022,26)=(112233)$} \label{fig:demovoeg}
\end{figure}

We also define the inverse map $\J_m$.
\begin{defn} \label{defjoin}
Let $(\xi,\eta)\in \P_{n,2}$. We place horizontal blocks of length
$\xi_i$ and vertical blocks of length $\eta_j$ such that we obtain
the $m$-tableau of a partition, in the following way. Each time,
we look for the block which can reach the maximal number in the
$m$-tableau. The blocks are nested in the sense that after every
step the placed blocks form the diagram of a partition. We repeat
this process until we have placed all blocks. The end result is
the $m$-tableau $T_m(\l)$ of the partition $\l:=\J_m(\xi,\eta)$.
Two examples for $m=1$ are provided in Figure \ref{fig:demovoeg}.
\end{defn}

This map is also well-defined:\begin{lemma}
For arbitrary parameter values $k_2=mk_1$, $\J_m$ is well defined on
$\Sigma_m(W_0c)$ for every residual point $c$.
\end{lemma}

\pf: (i) For generic parameters, this follows analogously to
Corollary \ref{gendistinct}.

(ii) Now suppose that the parameters are special. First we observe that it follows from the definition of $\Sigma_m(W_0c)$ that all the entries in the $m$-symbols of $\Sigma_m(W_0c)$ are mutually distinct.

(ii-a) First we assume $m \geq 0$ is integer.
Suppose $J(c)=\{j_1,\dots,j_{m+2r}\}$. Then $\J_m$ is not
well-defined if after placing $h$ horizontal and $v$ vertical
blocks, both the horizontal block of length $\xi_{r+m-h}$ and the
vertical block of length $\eta_{r-v}$ can reach the same number
$x$ in the $m$-tableau. After $h$ horizontal and $v$ vertical
blocks, the next block will start on a square containing
$|m-h+v|$. Suppose first that $m \leq h-v$. Then the square under
consideration lies above the zero diagonal and hence we find that
we have $\xi_{r+m-h}=x-(m-h+v)+1$ and $\eta_{r-v}=x+(m-h+v)+1$,
i.e., that $\xi_{r+m-h}=\eta_{r-v}-2(m-h+l)$. But then
$e_m(\xi_{r+m-h})=e_m(\eta_{r-v})$, which is a contradiction.  The case $m<h-v$ is treated along the same lines.

(ii-b) Now suppose that $0<m\in\half\Z, m\notin \Z$. Suppose  $J(c)=\{j_1,\dots,\allowbreak j_{\ceil{m}+2r}$. Again, $\J_m$
will only be not well defined if after placing $h$ horizontal and
$v$ vertical blocks, both the horizontal block of length
$\xi_{r+\frac{1}{2}+m-h}$ and the vertical block of length
$\eta_{r-v}$ can reach the same number $x$ in the $m$-tableau.
The first placed of these blocks starts on a square containing $|m-h+v|$.
Suppose for the moment that $m> h-v$, then this square lies on or
above the upper $\frac{1}{2}$-diagonal, and we find that the block
lengths $\xi_{r+m+\frac{1}{2}-h}$ and $\eta_{r-v}$ satisfy
$\xi_{r+m+\frac{1}{2}-h}=x-(m-h+v)+1$ and
$\eta_{r-v}=x+(m-h+v)+1$. But again, this means that $e_m(\xi_{r+m+\frac{1}{2}-h})=e_m(\eta_{r-v})$.
The case where $m<h-v$ can be treated in the same way.

(iii) Finally we remark that the case $m<0$ is analogous.\qed

\subsection{Confluence of residual points} We want to show that $\S_m(\P_m(W_0c))=\Sigma_m(W_0c)$. First we show
the existence of a canonical ``balanced'' $m$-tableau of a
partition $\l(c) \in \P_m(W_0c)$.

\begin{lemma}\label{startpart}
\label{begintableau} Let $(\xi(c),\eta(c))$ be the double partition
described in section \ref{sigmarespunt}. Then $\l(c)=\J_m(\xi(c),\eta(c))\in \P_m(W_0c)$.
\end{lemma}

\pf: (i) Let $m\geq 0$ and let $\tilde{c} \in \L_{-m}$ be the residual point for which $J(\tilde{c})=J(c)$.  Then, if $c=c(\l,k,mk)$, we have $\tilde{c}=c(\l',k,-mk)$ and hence the extremeties of $T_{-m}(\l)$ and of $T_m(\l')$ are the same. Since $(\xi(\tilde{c}),\eta(\tilde{c}))=(\eta(c),\xi(c))$, one has $\P_{-m}(W_0\tilde{c})=\P_m(W_0c)'$. It is easy to see that in general
$\J_m(\xi,\eta)=\J_{-m}(\eta,\xi)'$. Therefore, if $\l(c) \in \P_m(W_0c)$, we have $\J_{-m}(\xi(\tilde{c}),\eta(\tilde{c}))=\J_m(\xi(c),\eta(c))'=\l(c)'\in \P_{-m}(W_0\tilde{c})$. Hence, it is sufficient to consider $m\geq 0$.

\begin{figure}[htb]
\begin{center}
{\includegraphics[angle=0,scale=0.4]{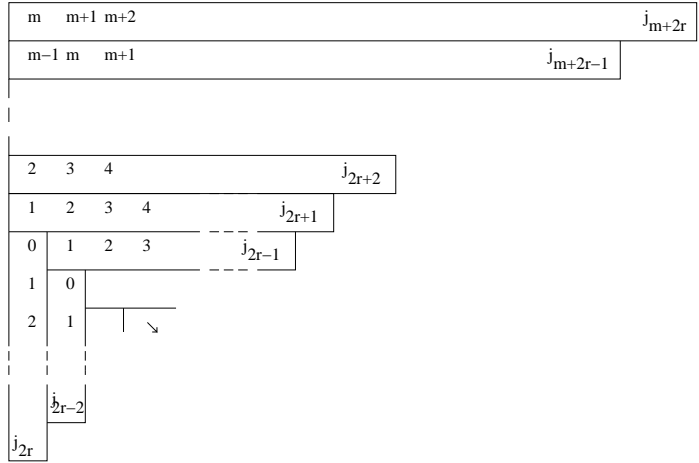}}
\end{center}
\caption{Construction of $\l(c)=\J_m(\xi(c),\eta(c))$}
\label{fig:begin}
\end{figure}

(ii) Suppose $m\in\Z$ and let $c$ have jumps $j_1< j_2< \dots< j_{m+2r}$. The lengths
of the parts of $(\xi(c),\eta(c))$ are $j_1<j_2+1\leq \dots\leq
j_{2r-1}<j_{2r}+1\leq j_{2r+1}\leq j_{2r+2}-1 \leq \dots \leq
j_{2r+m}-(m-1)$. This means that if we carry out the map $\J_m$
(see Figure \ref{fig:begin}), we first place a horizontal block of
length $j_{2r+m}-(m-1)$, which contains
$(m,m+1,m+2,\dots,j_{m+2r})$. Next will be a horizontal block of
length $j_{m+2r-1}-(m-2)$ containing
$(m-1,m,m+1,\dots,j_{2r+m-1})$. If $r=0$ the procedure stops after
having placed at most $m-1$ horizontal blocks. Suppose now that
$r>0$. Then, after having placed the $m-1$ horizontal blocks, we
still have to place $r+1$ horizontal blocks of length
$j_1,j_3,\dots,j_{2r+1}$ and $r$ vertical blocks of length
$j_2+1,j_4+1,\dots,j_{2r}+1$ on a 1-tableau. We first get a
horizontal block of length $j_{2r+1}$. This is clear since we can
reach $j_{2r+1}$ with a horizontal block and only
$j_{2r}-1<j_{2r+1}$ with a vertical block. The $(m+1)^{st}$ block
will have length $j_{2r}+1$, which is placed vertically and
contains $(0,1,\dots, j_{2r})$. For the $(m+2)^{nd}$ block, we
find after a similar consideration that we get a horizontal block
of length $j_{2r-1}$, containing $(1,2,\dots,j_{m+2})$. Notice
that as $j_{2r-1}+1<j_{2r+1}$, the intermediate diagram we have
obtained is still that of a partition. In the same way the
$(m+3)^{rd}$ block is vertical of length $j_{2r-2}+1$, with
content $(0,1,2,\dots,j_{2r-2})$. We continue in this way to place
all blocks. In the end, we have a diagram covering a sequence with
jumps $j_1,j_2,\dots,j_{m+2r}$, together with a zero for every
second block containing a one, which is the right amount.

(ii) Now suppose that $m$ is half-integer. Then the idea of the
proof for integer $m$ carries over without modification: if $r=0$
we need at most $m-\frac{1}{2}$ horizontal blocks, and if $r>0$ we
find the same alternating pattern of horizontal and vertical
blocks after the first $m-\frac{1}{2}$ horizontal ones. \qed
This gives us a starting point for the 1-1 correspondence, because
from Figure \ref{fig:begin} it is clear that also
$\S_m(\l(c))=(\xi(c),\eta(c))$.

Let $\l \in \P_m(W_0c)$. From $\l$ we can construct
another partition $\mu$ in $\P_m(W_0c)$ by cutting off a piece of
$\l$ and reattaching it elsewhere. Recall that the extremeties of
$\l$ are jumps for $c$. Therefore, consider the case where it is
possible to cut from $\l$ a block which contains the sequence
$(j_{p-1}+1,j_{p-1}+2,\dots,j_p)$ in $T_m(\l)$, and that it can be reattached elsewhere. Then there is a
unique other position where it may be reattached: to the
extremity of $T_m(\l)$ containing $j_{p-1}$. Let the new partition
thus obtained be $\mu$. Since $\mu$ is uniquely determined by $\l$
and $p$, we write $\mu=\flip_p(\l)$.

We show that these flips do not change the similarity class of the split tableaux. Denote, for a block $\xi_i$ which appears in $\S_m(\l)$, by $j(\xi_i)$ the extremity on which block $\xi_i$ ends.
\begin{prop}\label{symbol} Let $\l \in \P_m(W_0c)$. If $\flip_p(\l)$ exists, then $\S_m(\l) \sim_m \S_m(\flip_p(\l))$.
\end{prop}

\pf: Denote $\S_m(\l)=(\xi,\eta)$ and $\S_m(\flip_p(\l))=(\xi',\eta')$.

\begin{figure}[htb]
\begin{center}
{\includegraphics[angle=0,scale=0.4]{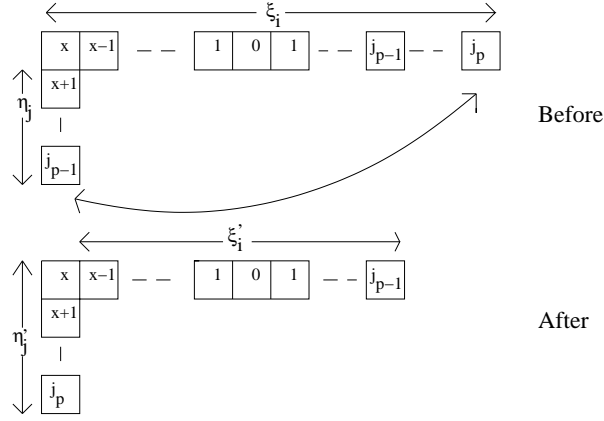}}
\end{center}
\caption{Example of flip}
\label{fig:flip}
\end{figure}

(i) Suppose $m$ is integer. We consider the flip of a part of
length $l=j_p-j_{p-1}$ from $\xi_i$ to $\eta_j$. Since
$j(\xi_i)=j_p$ and $j(\eta_j)=j_{p-1}$, it is clear that $\xi_i$
and $\eta_j$ are neighbours in the sense that block $\eta_j$ is
selected by $\S_m$ immediately after block $\xi_i$. Since there
are $m+r-i$ horizontal and $r-j$ vertical blocks above $\eta_j$,
it follows that $\xi_i$ starts on a square containing
$x:=|m-(m+r-i)-(r-j)|=|i-j|$. Suppose for the moment that $i\leq
j$, the case $i>j$ being similar. Then the starting square of
$\xi_i$ lies to the left of the zero-diagonal, and the flip is
visualized in Figure \ref{fig:flip}.

We see that we have $\xi_i=x+j_p+1, \eta_j=j_{p-1}-x$ and $\xi_i'=x+j_{p-1}, \eta_j'=j_p-x+1$. In order for the $m$-symbol to be preserved, we need to show that $\xi_i+2(i-1)=\eta_j'+2(j-1)$. But this is equivalent to $x+j_p+1+2(i-1)=j_p-x+1+2(j-1) \iff 2x=2(j-i)$ which is the case. The case where $j<i$ is treated with a similar calculation, which we omit here.

(ii) Suppose that $m$ is half-integer. The proof carries over directly.
 \qed

Next we show that using these flips, we can construct any partition in $\P_m(W_0c)$, starting from $\l(c)$.
\begin{prop}\label{allpart}
Let $\l \in \P_m(W_0c)$. Then there exists a sequence $i_1,i_2,\dots,i_f$ such that $\l=\prod_{a=1}^f \flip_{i_a}(\l(c))$.
\end{prop}

\pf:  Obviously it suffices to prove the case where $m\geq 0$. We only treat the case where $m$ is integer, the case where
$m$ is half-integer being analogous. Let
$(\xi,\eta)=\S_m(\l)$. Let the boxes of $(\xi,\eta)$ emerge from
the splitting procedure in the order $b_1,b_2,\dots,b_{m+2r}$.
Notice that it may happen that some boxes are empty, i.e., there
are less than $m+2r$ boxes. The non-empty boxes all end on a jump. We
put $j(b_i)$ to be the jump on which box $b_i$ ends, where $j(b_i):=0$ if $b_i$ is empty. The diagram $\l(c)$ is the
unique diagram for which
\[ j(\xi_1) \leq j(\eta_1) \leq j(\xi_2) \leq \dots \leq j(\xi_r) \leq j(\eta_r) \leq j(\xi_{r+1}) \leq j(\xi_{r+2}) \leq \dots \leq j(\xi_{r+m}). \]
If we look at the corresponding series for $\l \neq \l(c)$, this
series of inequalities will be violated in a number of places.
However, we always have $j(\xi_1) \leq j(\xi_2) \leq \dots \leq
j(\xi_{m+r})$ and $j(\eta_1) \leq j(\eta_2) \leq \dots \leq
j(\eta_r)$.  Now take the first violated inequality, for
simplicity we assume that
\begin{equation}\label{ass}  j(\xi_1)\leq j(\eta_1)\leq j(\xi_2) \leq \dots \leq j(\xi_{i-1})\leq j(\eta_{i-1}) \leq j(\xi_i) \nleq j(\eta_i), \end{equation}
the case where the first violated inequality is of the kind
$j(\eta_i) \nleq j(\xi_{i+1})$ admitting a similar consideration.

This means that the part $\xi_i$ is too big. We now show that $\flip_i(\l)$ exists. Say $j(\xi_i)=j_p$,
then $\xi_i$ ends on a series of squares containing
$(j_{p-1}+1,j_{p-1}+2,\dots,j_p)$, and we can cut off these
squares while still maintaining the graph of a partition if
$j(\xi_{i-1})\neq j_{p-1}$. But this inequality holds, since if
$j(\xi_{i-1})= j_{p-1}$, we would have $j_{p-1} \leq
j(\eta_{i-1})\leq j_p$, which is impossible. So we cut off the
rectangle. Since $j_{p-1}$ is a jump, it follows from Lemma
\ref{extr} that there are now two extremity squares containing
$j_{p-1}$, unless $j_{p-1}$ does not appear in $T_m(\l)$. But if that were true, then $j_p$ is the smallest jump appearing as extremity, giving a contradiction to \eqref{ass}.
Thus, after cutting off the rectangle, we have two extremities equal to $j_{p-1}$,  one in the former box $\xi_i$ and one other. This
latter one can not belong to a horizontal box, since in that case
again we would have $j(\xi_l)=j_{p-1}$ for some $l<i$, and we
would find an earlier violated inequality. So we attach the
rectangle vertically to this square, which is the last square of
$\eta_l$ for some $l$. If $l<i$, then
$j_{p-1}=j(\eta_l)<j(\eta_i)=j_t$, so $p-1<t$. On the other hand,
we also have $j_p=j(\xi_i)>j(\eta_i)=j_t$ so $p>t$. This is a
contradiction and so we see that $l \geq i$. After the flip, we
still have the diagram of a partition if $j_s=j(\eta_{l+1})>j_p$,
i.e., if $s>p$. But before the flip we had
$j_s=j(\eta_{l+1})>j(\eta_l)=j_{p-1}$, so $s>p-1$. Since $s \neq
p$, it follows that indeed we still have the diagram of a
partition after this flip, and moreover the inequality we were
looking at, is less or no longer violated by $\flip_i(\l)$. We keep applying this
procedure until we reach $\l(c)$. \qed

In $\Sigma_m(W_0c)$, we can do something similar. First we define
\begin{defn}Consider
a $m$-tableau $T_m(\l)$ such that $\S_{m}(\l)$ has blocks
$b_i$. We denote by $e_m(b_i)$ the entry of block $b_i$ in the
$m$-symbol of $\S_m(\l)$.
\end{defn}

Notice that $e_m$ depends on the chosen lengths of the rows of the
symbol. However, we will only be interested in the difference of
two entries, which is independent of this choice.

It now turns out that the $m$-symbol of
$(\xi,\eta)\in \Sigma_m(W_0c)$ determines the order in which the
parts of $\xi$ and $\eta$ are placed when applying $\J_m$:

\begin{prop} \label{volgorde}
If we apply $\J_m$ to $(\xi,\eta)$, then the blocks are laid down
in the order of decreasing entries in the $m$-symbol, i.e., block
$b_i$ is placed before block $b_j$ if $e_m(b_i)>e_m(b_j)$.
\end{prop}

\pf: Again, to simplify notation, we only treat the case where $m
\geq 0 $ is integer, the other cases where $m$ is half-integer or
negative being similar. We use induction on the number of blocks
already placed. Suppose this is zero, then we either place
$\xi_{r+m}$ or $\eta_{r}$, and $\xi_{r+m}$ will be selected if and
only if it can reach a higher number than $\eta_r$, which is the
case if and only if $m+\xi_{r+m}-1 > \eta_r-m-1 \iff
\xi_{r+m}+2(r+m-1) > \eta_r+2(r-1)$, which is what we wanted.

Now suppose that we have placed already $k$ horizontal blocks
$\xi_{r+m},\xi_{r+m-1},\dots,$ $\xi_{r+m-k+1}$ and $l$ vertical
blocks $\eta_r,\eta_{r-1},\dots,\eta_{r-l+1}$. Then the next block
starts on a square containing $|m-k+l|$. First we assume that $m
\geq k-l$, i.e., the square under consideration is above the
zero-diagonal. Then the next block to be placed will either be
$\xi_{r+m-k}$ or $\eta_{r-l}$, and it will be $\xi_{r+m-k}$ if and
only if
\begin{eqnarray*} &&m-k+l+\xi_{r+m-k}-1  >  \eta_{r-l}-m+k-l-1 \\ &\iff& \xi_{m+m-k}+2(m-k) > \eta_{r-l}-2l \\ &\iff& \xi_{r+m-k}+2(r+m-k-1) > \eta_{r-l}-2(r-l-1), \end{eqnarray*}

Second, assume that $m<k-l$, then the starting square for the next block is below the zero-diagonal. In this case, the next block will be $\xi_{r+m-k}$ if and only if \begin{eqnarray*} &&\xi_{r+m-k}-(-m+k-l)-1  >  \eta_{r-l}-m+k-l-1 \\ &\iff& \xi_{r+m-k}+2(m-k) > \eta_{r-l}-2l \\ &\iff& \xi_{r+m-k}+2(r+m-k-1) > \eta_{r-l}-2(r-l-1). \end{eqnarray*}

This proves the induction step and therefore the claim.\qed

Then a flip in the $m$-symbol this corresponds to the following:
\begin{lemma}
In the situation of Figure \ref{fig:flip} we have $e_m(\xi_i)>e_m(\eta_i)$ and
$e_m(\eta_j)<e_m(\xi_i)<e_m(\eta_{j+1})$. Therefore the flip
interchanges $e_m(\xi_i)$ with the unique entry on the other row
of the $m$-symbol with which it can be interchanged.
\end{lemma}

\pf: In the notation of the proof of Proposition \ref{symbol},
$\xi_i=j-i+j_p+1$ and $\eta_j=j_{p-1}-j+i$, therefore
$e_m(\xi_i)-e_m(\eta_j)=\xi_i+2(i-1)-\eta_j-2(j-1) \geq 2$. For
$\eta_{j+1}$ we know that $\eta_{j+1} \geq j_{p+1}-(x+1)+1$, and
therefore $e_m(\eta_{j+1})-e_m(\xi_i) \geq j_{p+1}-j_p+1 \geq 2$.
Finally, since $j \geq i$, $\xi_i \geq \eta_j \geq \eta_i$ and so
also $e_m(\xi_i) \geq e_m(\eta_i)$. \qed

On the other hand, given any $m$-symbol of a partition in
$\Sigma_m(W_0c)$, then we can change it into the similar symbol with
the same, but increasing entries, i.e. where
$e_m(\xi_1)<e_m(\eta_1)<e_m(\xi_2)<\dots<e_m(\eta_{r-1})<e_m(\xi_r)<e_m(\xi_{r+1})<\dots<e_m(\xi_{m+r})$.
We can for example achieve this by finding the first violated
inequality, e.g. say $e_m(\xi_i)>e_m(\eta_i)$, and then
interchanging $e_m(\xi_i)$ with $e_m(\eta_j)$ for the unique $j$
such that $e_m(\eta_j)<e_m(\xi_i)<e_m(\eta_{j+1})$. We denote the
permutation where $e_m(\xi_i)$ is moved to the bottom row, where
it takes the place of entry $e_m(\eta_j)$, by $\flip^{\downarrow
i}(\xi,\eta)=\flip^{\uparrow j}(\xi,\eta)$.

On the level of partitions, such a permutation of the $m$-symbol
always corresponds to a flip:

\begin{lemma} \label{flipflip}
Let $(\xi,\eta) \in \Sigma_m(W_0c)$ and let $\l=\J_m(\xi,\eta)$.
Suppose
$e_m(\xi_1)<e_m(\eta_1)<e_m(\xi_2)<\dots<e_m(\eta_{i-1})<e_m(\xi_i)$
but $e_m(\xi_i)>e_m(\eta_i)$. Then $\J_m(\mathcal{F}^{\downarrow
i}(\xi,\eta))={\mathcal F}_p(\l)$ if the box of $\xi_i$ ends on
jump $j_p$ in $T_m(\l)$.
\end{lemma}

\pf: Let $\eta_j$ be such that
$e_m(\eta_j)<e_m(\xi_i)<e_m(\eta_{j+1})$. Then clearly $j \geq i$.
From Proposition \ref{volgorde} it follows that $\eta_j$ is placed
after $\xi_i$. Since $e_m(\xi_{i-1})<e_m(\eta_j)$ it follows that
$\eta_j$ is placed immediately after $\xi_i$. Say block $\xi_i$
ends on $x=j_p$ and block $\eta_j$ ends on $y<x$, then we can cut
a rectangle containing $y+1,y+2,\dots,x$ off $\xi_i$ and attach it
to $\eta_j$ provided $j(\xi_{i-1})\neq y$ and $j(\eta_{j+1}) \neq
x$. But this would contradict the fact that $\J_m$ is well defined
on $\Sigma_m(W_0c)$. As we have seen in the proof of Proposition
\ref{allpart}, this flip of $\xi_i$ to $\eta_j$ corresponds with
the interchange of $e_m(\xi_i)$ and $e_m(\eta_j)$ in the
$m$-symbol. \qed

Finally, we can show
\begin{thm}\label{puntenstelling}
The maps $\S_m$ and $\J_m$ implement a bijection  $\P_m(W_0c) \leftrightarrow \Sigma_m(W_0c)$.
\end{thm}

\pf: Let $\l \in \P_m(W_0c)$, then it follows from \ref{startpart} and \ref{allpart} that $\S_m(\l) \in \Sigma_m(W_0c)$. Since it is easy to see that $\S_m$ is injective it remains to show that the image of $\S_m$ is all of $\Sigma_m(W_0c)$.
  Take $(\xi,\eta) \in
\Sigma_m(W_0c)$, then there is a sequence of permutations ${\mathcal
F}_{\downarrow l_j}$ such that \[ (\xi,\eta)=\prod_j{\mathcal
F}^{\downarrow l_j}(\xi(c),\eta(c)).\] Therefore, by Lemma
\ref{flipflip}, there is a sequence of jumps $t_j$ such that
\[ (\xi,\eta)=\prod_j{\mathcal
F}^{\downarrow l_j}(\xi(c),\eta(c))=\prod_j {\mathcal F}^{\downarrow
l_j}(\S_m(\l(c)))=\S_m(\prod_j{\mathcal F}_{t_j}(\l))\in
\S_m(\P_m(W_0c)).\] This proves the claim. \qed

\subsection{Properties of truncated induction} Recall the truncated induction $\trm$. We will prove its
transitivity. We first consider its analogue on $\P_n$.
Consider a partition $\l \vdash n$, then the analogue of the symbol is just $\l$ itself. In this section we will write
(as is more common) the parts of all partitions in decreasing
order. Then the $a$-function reduces to the well-known
\begin{equation}
\label{n} n(\l)=\sum_i(i-1)\l_i=\sum_i{\l'_i \choose 2},
\end{equation}
where $\l'$ denotes the partition conjugate to $\l$.

In this case, the partial ordering defined by the $a$-function
thus refines the dominance ordering, in which by definition $\l
\geq \mu$ if
\begin{equation}
\label{dominance} \l_1+\l_2+\dots+\l_i \geq
\mu_1+\mu_2+\dots+\mu_i {\rm \ for\ all\ } i \in \{1,2\dots,n\}.
\end{equation}

\begin{lemma}
\label{n=dom} For all $\l,\mu \vdash n, \l \geq \mu \Rightarrow
n(\l)\leq n(\mu)$.
\end{lemma}

\pf: $\l \geq \mu \iff \l' \leq \mu'\Rightarrow \sum_i {\l'_i
\choose 2} \leq \sum_i {\mu'_i \choose 2} \iff n(\l) \leq
n(\mu)$.\qed

Now let us see how the $a$-function behaves under induction.
Recall the Littlewood-Richard\-son rule (cf. \cite{Mac} for a
proof). For a partition $\l \vdash n$ we denote by $V_\l$ the
corresponding irreducible $S_n$-module. If $\mu \vdash k$ and $\nu
\vdash l$, then

\begin{equation}
\label{LW} {\rm Ind}_{S_k \times S_l}^{S_{k+l}}(V_\mu \otimes
V_\nu)=\sum_{\l\vdash n} c_{\mu\nu}^\l V_\l.
\end{equation}

In this formula $c_{\mu\nu}^\l$ is the number of ways the Young
tableau of $\mu$ can be made into the Young tableau of $\l$ by
means of a strict $\nu$-expansion. This means the following. We
add $|\nu|$ squares to the tableau of $\mu$ to obtain the tableau
of $\l$, by first adding $\nu_1$ squares filled with a 1, then
adding $\nu_2$ squares filled with a 2, etc, until finally we add
$\nu_n$ squares filled with $n$. After each of these steps, we
must still have the diagram of a partition, and no two squares in
the same column may have the same entry. It now remains to explain
when such an expansion is called strict. Suppose we read the added
squares from right to left, and from top to bottom. We then have a
sequence of numbers $a_1a_2\dots a_l$. Now the expansion is called
strict if for any intermediate sequence $a_1a_2\dots a_i$ ($1\leq
i \leq l$), any of the occurring integers $1\leq k \leq n$ occurs
at least as many times as the next integer $k+1$.

The consequence is that

\begin{lemma}
\label{cup} Let $\mu\vdash k$ and $\nu \vdash l$. Then
\[ {\rm tr-Ind}_{S_k \times S_l}^{S_{k+l}}(V_\mu \otimes V_\nu)=V_{\mu \cup \nu}.\]
\end{lemma}

\pf: First we check that $\l =\mu \cup \nu$ indeed occurs in the
induction, by constructing its Young tableau as a strict
$\nu$-expansion of the Young tableau of $\mu$. This can be seen as
follows. Fill the part $\nu_i$ in $\l$ with $i's$. Then we obtain
a series of $l=l(\nu)$ horizontal bars in $\l$, filled with 1 up
to $l$. Now move the contents of the filled squares as far down as
possible. Then the unfilled squares make up a tableau of shape
$\mu$ in $\l$, so that we obtain $\l$ as a $\nu$-expansion of
$\mu$. This expansion is strict since if we have a square filled
with $i$, then $i-1$ occurs more often than $i$ in the preceding
boxes of the expansion, by construction.

Thus indeed $\l:=\mu \cup \nu$ occurs in \eqref{LW}. Now we show
that its $a$-value is maximal. This follows from \ref{n=dom},
since
\[ n(\l)=n(\mu \cup \nu)=\sum_i { (\mu \cup \nu)'_i \choose 2} =\sum_i { \mu'_i+\nu'_i \choose 2},\]
but for any $\kappa\neq\l$ occurring in \eqref{LW}, $\kappa' <
\mu' \cup \nu'$. \qed

\begin{cor}
\label{ntrans} The truncated induction is transitive in the sense
that for all $\l \vdash k, \mu \vdash l, \nu \vdash m$, we have
\[ {\rm tr-Ind}_{S_k\times S_{l+m}}^{S_{k+l+m}}(V_\l \otimes {\rm tr-Ind}_{S_l \times S_m}^{S_{l+m}}(V_\mu \otimes V_\nu))={\rm tr-Ind}_{S_k \times S_l \times S_m}^{S_{k+l+m}}(V_\l \otimes V_\mu \otimes V_\nu). \]
\end{cor}

\pf: The argument of the previous Lemma generalizes immediately to
any number of partitions, which implies that both sides are equal
to $V_{\l \cup \mu \cup \nu}$.\qed

We want to prove that this also holds for ${\rm tr}_m{\rm -Ind}$.
Let us first see how the $a$-value of a partition is related to
$a_m(\xi,\eta)$.

\begin{lemma}
\label{kappa}
 Let $\l$ and $\mu$ be such that $n(\l)\leq n(\mu)$.
Then $a_m(\l,\kappa)\leq a_m(\mu,\kappa)$ and $a_m(\kappa,\l) \leq
a_m(\kappa,\mu)$ for all $m\in\half \Z$ and all partitions
$\kappa$.
\end{lemma}

\pf: The function $a_m(\xi,\eta)$ is nothing but the $a$-function
on the partition whose parts are the entries of the $m$-symbol of
$(\xi,\eta)$. The lemma therefore follows from the fact that
$n(\l\cup\kappa) \leq n(\mu\cup\kappa)$ for all $\kappa$, if
$n(\l)\leq n(\mu)$.\qed

\subsubsection{Transitivity of $\trm$}We now return to ${\rm tr}_m-{\rm Ind}$. We check what happens in
case $\l$ has only one part:
\begin{prop}
\label{1deel} Let $k_2=mk_1$ with $k_1 \neq 0, m \in \half\Z$.
Then for every 2-partition $(\xi,\eta)$ of $n$,
\[ {\rm tr}_m{\rm -Ind}_{W_0(B_l) \times S_t}^{W_0(B_n)}( (\xi,\eta) \otimes (t))=\sum_i(\xi,\eta) \cup (t_{i,1},t_{i,2}),\]
where $i$ runs over an indexing set containing one or two
elements, i.e., for one or two decompositions $t_{i,1}+t_{i,2}=t$.
\end{prop}

\pf: Suppose that in a constituent $\chi_A$ of maximal $a$-value,
$t_{1,1}$ squares are added to $\xi$ and $t_{1,2}$ to $\eta$. Then
by Lemma \ref{kappa}, $A=(\xi \cup t_{1,1},\eta\cup t_{1,2})$, and
$\chi_A$ occurs with multiplicity one in the induction. The
entries $e_m(t_{1,j})$ corresponding to $t_{1,1},t_{1,2}$ in the
$m$-symbol of $(\xi\cup t_{1,1},\eta\cup t_{1,2})$ are equal or
have a difference of one. Suppose we are in the latter case and
$e_m(t_{1,1})=e_m(t_{1,2})+1$. Then, {\it if} $\chi_B$ with $B=(\xi \cup
(t_{1,1}-1), \eta \cup (t_{1,2}+1))$ occurs in the induction, we
have $a_m(\chi_A)=a_m(\chi_B)$. Put $t_{2,1}=t_{1,1}-1$ and
$t_{2,2}=t_{1,2}+1$. \qed

Now we investigate the transitive behaviour of ${\rm tr}_m{\rm
-Ind}$:
\begin{lemma}
\label{trans}
 Let $(\xi,\eta)\in\P_{l,2}$, and $a,b \in
 {\mathbb N}$. Suppose that $(\xi \cup a_1, \eta \cup a_2)$ occurs in ${\rm tr}_m{\rm -Ind}_{S_a \times W_0(B_l)}^{W_0(B_n)}((a) \otimes (\xi,\eta))$ and that  $(\xi \cup b_1, \eta \cup b_2)$ occurs in ${\rm tr}_m{\rm -Ind}_{S_b \times W_0(B_l)}^{W_0(B_n)}\allowbreak ((b) \otimes (\xi,\eta))$. Then $(\xi,\eta) \cup (a_1,a_2) \cup (b_1,b_2)$ occurs in
\begin{equation}
\label{ab} {\rm tr}_m{\rm -Ind}_{S_a \times S_b \times
W_0(B_l)}^{W_0(B_n)}((a)\otimes (b) \otimes (\xi,\eta)).
\end{equation}
Conversely, any constituent of \eqref{ab} is of this form.
\end{lemma}

\pf: Suppose for convenience that $a \geq b$. Remark that ${\rm
Ind}_{S_a \times S_b}^{W_0({B_{a+b}})}((a)\otimes(b))$ consists of
2-partitions $(x_1y_1,x_2y_2)$, where $x_i\geq y_i$, $x_1+x_2 \geq
a$ and $y_1+y_2\leq b$. It is not hard to see that $a\geq b$
implies $a_i \geq b_i$ and therefore $(\xi,\eta) \cup
(a_1b_1,a_2b_2)$ occurs in the induction \eqref{ab}. It 
remains to see that this is indeed a constituent of maximal
$a$-value. It is clear from \ref{cup} that \[ {\rm
tr-Ind}_{S_{|\xi|}\times S_{a_1} \times
S_{b_1}}^{S_{|\xi|+a_1+b_1}}(\xi \otimes (a_1) \otimes (b_1))=\xi
\cup a_1\cup b_1\] and similarly for ${\rm tr-Ind}(\eta \otimes
(a_2)\otimes (b_2))$. By Lemma \ref{kappa}, we only need to check
now that the $a_m$-value does not increase by moving an induced
square from $\xi \cup a_1 \cup b_1$ to $\eta \cup a_2 \cup b_2$, in the set of 2-partitions occurring in the induction.
Suppose we move a square of $b_1$, then we either obtain $(\xi\cup
a_1 \cup (b_1-1), \eta \cup a_2 \cup (b_2+1))$ or $(\xi \cup a_1
\cup (b_1-1), \eta \cup (a_2+1) \cup b_2)$. However, the first
possibility would imply that $a_m(\xi \cup (b_1-1), \eta \cup
(b_2+1))>a_m(\xi \cup b_1,\eta \cup b_2)$ which is a
contradiction. Now consider the second case. It is clear that $e_m(a_1) \in \{e_m(a_2)-1,e_m(a_2),e_m(a_2)+1\}$ and likewise for the $b_i$. Since $b_1 \leq a_1$ it follows that $e_m(b_1)+2 \leq e_m(a_1)$, hence also $e_m(b_1)<e_m(a_1)$. Therefore the second possibility also does not increase the $a_m$-value.
 Similarly we exclude the possibility that $a_m(\xi \cup (a_1-1) \cup b_2,\eta
\cup (a_2+1) \cup (b_2))>a_m(\xi\cup a_1 \cup b_1,\eta \cup a_2
\cup b_2)$, since this would imply that $a(\xi\cup
(a_1-1),\eta\cup (a_2+1))>a(\xi\cup a_1,\eta \cup a_2)$.  The
representation $(\xi \cup (a_1-1) \cup b_1,\eta\cup a_2 \cup
(b_2+1))$ does not occur in the induction, since $b_1+b_2+1>b$.
Since we can reach any constituent of the induction by a series of
these moves of one square, and the $a_m$-value increases strictly
along a shortest path towards the 2-partition with maximal
$a_m$-value, we see that indeed the $a_m$-value of $(\xi,\eta)
\cup (a_1,a_2) \cup (b_1,b_2)$ is maximal.

Conversely, suppose that $(\xi',\eta')$ is a constituent of
\eqref{ab}. Then it is of the form $(\xi', \eta')=(\xi,\eta) \cup
(x_1y_1,x_2y_2)$ for some $x_i,y_i$ as above. Maximality of the
$a$-value implies that $x_1+x_2=a$, and so $y_1+y_2=b$. But then,
keeping $(x_1,y_1)$ fixed, it is clear that
$(\xi,\eta)\cup(x_2,y_2)$ has maximal $a_m$-value among the
constituents of ${\rm Ind}_{S_a \times W_0(B_l)}^{W_0(B_{l+a})}
((a) \otimes (\xi,\eta))$, and likewise for $(\xi,\eta) \cup
(x_1,y_1)$. \qed

\begin{cor}
\label{atrans} Let $k_2=mk_1\neq 0$ be such that $m \in \half \Z$,
and consider $(\xi,\eta)\in \P_{l,2}$. Let $a+b+l=n$. Then

\[ \begin{array}{l}  {\rm tr}_m{\rm -Ind}_{S_a \times W_0(B_{b+l})}^{W_0(B_n)}((a) \otimes {\rm tr}_m{\rm -Ind}_{S_b \times W_0(B_l)}^{W_0(B_{b+l})}( (b) \otimes(\xi,\eta)))= \end{array}\]\[ \begin{array}{r} {\rm tr}_m{\rm -Ind}_{S_a\times S_b \times W_0(B_l)}^{W_0(B_n)}((a)\otimes (b) \otimes (\xi,\eta)).\end{array}   \]
\end{cor}
\pf: By the above Lemma, both sides are equal to $\sum_{i,j}
(\xi,\eta) \cup (a_{i,1},a_{i,2}) \cup (b_{j,1},b_{j,2})$ for
certain $a_{i,1}+a_{i,2}=a, b_{j,1}+b_{j,2}=b$. \qed

By repeated application of this corollary, we get:
\begin{prop} Let $k_2=mk_1\neq 0$ and $m \in \half \Z$. Let $(\xi,\eta)\in \P_{l,2}$ and $\l
\vdash n-l$. Then
\begin{equation}
\label{amax}
 {\rm tr}_m{\rm -Ind}_{S_\l \times W_0(B_l)}^{W_0(B_n)}(\rm{triv}_\l \otimes (\xi,\eta))=\sum_\mu (\xi \cup \mu, \eta \cup (\l-\mu)),
\end{equation}
where $\mu$ ranges over a set of partitions, whose Young diagrams are
contained in the one of $\l$, which depends on $(\xi,\eta)$ and $m$.
\end{prop}

\subsection{Confluence of residual subspaces} In this section we show that $C_m(W_0c_L) \leftrightarrow \Sigma_m(W_0c_L)$ and that $\Sigma_m(W_0c_L)$ is a full similarity class. 

\subsubsection{} First we take a look at the case where a residual point is no longer residual at a given special value, but coincides with the center of a higher-dimensional residual subspace:

\begin{prop}
\label{tellenkorrespondenten} Let $k_1 \neq 0, m \in\half\Z$, and suppose that $c(\mu,k_1,mk_1)=c_L$ for $L\in \L_m(n)$ of type $A_\kappa \times (B_l,\nu)$. Then
the set $\Sigma_m(W_0c_L)$ of Springer correspondents of $W_0c_L$, i.e., all
the irreducible representations of $W_0$ in
\[ \sum_{ \{ \nu' | \S_m(\nu') \sim_m \S_m(\nu)\}} {\rm tr}_m{\rm -Ind}_{S_\kappa \times W_0(B_l)}^{W_0(B_n)}({\rm triv}_\kappa \otimes \S_m(\nu')), \]
satisfies \[ \Sigma_m(W_0c_L)=[\S_{m+\e}(\mu)]_m=[\S_{m-\e}(\mu)]_m.\]
Moreover, we have a bijection
\[ \Sigma_m(W_0c_L) \longleftrightarrow C_m(W_0c_L).\]
\end{prop}

\pf: We will prove the Proposition by induction. Therefore we introduce, for $i=0,1,\dots,l=l(\kappa)$:
\[ \kappa^{(i)}=(\kappa_1,\dots,\kappa_i); \mu^{(i)} \mbox{ such that }sp_m(\mu^{(i)})=(\kappa^{(i)},\nu),\]
i.e., $\mu^{(i)}$ is the partition obtained by removing the strips of $\kappa_{i+1},\dots,\kappa_l$ from $T_m(\mu)$. Put $L^{(i)}$ to be the residual coset of type $A_{\kappa^{(i)}}\times (B_l,\nu)$. 

We will use induction on $i$. Suppose that $i=0$. Then the statement reduces to Theorem \ref{puntenstelling}.

Therefore we take $i+1$ and suppose that the Proposition holds for $i$. That is, (putting $|\kappa^{(i)}|=n_i$)
\[ \Sigma_m(W_0c_{L^{(i)}})= \sum_{ \{ \nu' | \S_m(\nu') \sim_m \S_m(\nu)\}} {\rm tr}_m{\rm -Ind}_{S_{\kappa^{(i)}} \times W_0(B_l)}^{W_0(B_{n_i+l})}({\rm triv}_{\kappa^{(i)}} \otimes \S_m(\nu')), \]
satisfies the bijection
\[ \Sigma_m(W_0c_{L^{(i)}}) \leftrightarrow C_m(W_0c_{L^{(i)}}),\]
and
\[ \Sigma_m(W_0c_{L^{(i)}})=[\S_{m+\e}(\mu^{(i)})]_m=[\S_{m-\e}(\mu^{(i)})]_m.\]

Suppose that $\kappa_{i+1}=t$ and pick $(\xi,\eta) \in \Sigma_m(W_0c_{L^{(i)}})$.
We will show that
\begin{equation}\label{i+1} \trm_{S_t\times W_0(B_{n_i+l})}^{W_0(B_{n_{i+1}+l})}((t)\otimes (\xi,\eta))\end{equation} consists of two elements, which belong to $[\S_m(\mu^{(i+1)})]_m$.

In view of the following implication
\[ (\xi,\eta)\sim_m (\xi',\eta') \Rightarrow \trm((t)\otimes (\xi,\eta)) \sim_m \trm((t)\otimes (\xi',\eta')),\]
we may assume that $(\xi,\eta) \in \S_{m\pm\e}(\mu^{(i)})$. Doing so, we will calculate explicitly both \eqref{i+1} and $\S_{m\pm\e}(\mu^{(i+1)})$ and show that they agree.
Recall from \ref{altijdspec} that a strip in $T_m(\mu)$ corresponding to an $A$-factor is found inside $T_m(\mu)$ as a hook-shaped strip, whose extremities contain $\lfloor\frac{t}{2}\rfloor$ (recall also that $t$ is necessarily odd). For $m\pm\e$ the splitting procedure $\S_{m\pm\e}$ is well-defined; the square in the corner of the hook belongs to a horizontal box in one splitting, and to a vertical one in the other. 


We are going to construct the constituents $(\xi',\eta')$ whose symbol has maximal $a$-value in the induced representation \eqref{i+1}.

(i) Recall that the Littlewood-Richardson rule implies
that the diagrams of $\xi'$ resp $\eta'$ are obtained from those
of $\xi$ resp. $\eta$ by means of strict expansions. In this case,
this means that there are no two new squares in the same column. Suppose that $l(\xi)=m+r$, $l(\eta)=r$. Then we write $\xi=(0=\xi_0\leq \xi_1 \leq \dots \leq \xi_{r+m})$, $\eta=(0=\eta_0 \leq \eta_1 \leq \dots \leq \eta_r)$, and $\xi'=(\xi'_0 \leq \dots \leq \xi'_{m+r})$, $\eta'=(\eta'_0 \leq \dots \leq \eta'_r)$. We have
\begin{equation}
\label{induktievw}
\xi_i \leq \xi'_i \leq \xi_{i+1}, \eta_i\leq \eta'_i \leq \eta_{i+1},
\end{equation}
and  $(\xi,\eta)$ has
$m$-symbol
\begin{tiny}
\begin{equation}
\label{symb} \left( \begin{array}{lllllllllll} 0& & \xi_1+2 &
&\dots
&\xi_r+2r && \xi_{r+1}+2(r+1)& &\dots&\xi_{r+m}+2(r+m)\\
    &0  &         &\eta_1+2  &          & \dots&\eta_r+2r&& &  \end{array}
    \right)
\end{equation}
\end{tiny}
which we have to expand in such a way that we maximize the $a$-value.

(ii) Now we describe the splittings $\S_{m\pm\e}(\mu^{(i+1)})$. Suppose for convenience that $(\xi,\eta)=\S_{m+\e}(\mu^{(i)})$. Recall that by $j(b)$ we denote the extremity on which the block $b$ ends. Since the
strip of length $t$ fits into $T_m(\mu^{(i)})$, 
there exist $\xi_{a+1},\xi_a,\eta_{b+1},\eta_b$ such that $j(\xi_{a+1})>\floor{\frac{t}{2}},j(\eta_{b+1})>\floor{\frac{t}{2}}$, while $j(\xi_a)<\floor{\frac{t}{2}},j(\eta_a)<\floor{\frac{t}{2}}$. 
Several of these blocks may be empty,
in which case we consider these equalities to hold, since this can
only happen for the inner- or outermost blocks.

The $A$-factor is thus inserted in $T_m(\mu^{(i)})$ in the place of $(\xi_a,\eta_b)$. Since there are $m+r$ horizontal blocks and $r$ vertical ones, this means that the corner of the strip is positioned on a square containing $|a-b|$, as seen before. Therefore, $\S_{m\pm\e}(\mu^{(i+1)})=\S_{m\pm\e}(\mu^{(i)}) \cup (\xi'_a,\eta'_b)$, where $(\xi'_a,\eta'_b)$ is given by 
\begin{equation}
\label{verandering1}
 (i)-\left\{ \begin{array}{r} \xi'_a=\lceil \frac{t}{2} \rceil-(a-b)  \\ \eta'_b= \lfloor \frac{t}{2} \rfloor+(a-b)  \end{array} \right. {\rm for \ } m+\e {\rm \ or\ } (ii)-\left\{ \begin{array}{r} \xi'_a=\lfloor \frac{t}{2} \rfloor -(a-b)\\ \eta'_b= \lceil \frac{t}{2} \rceil+(a-b)  \end{array} \right. {\rm \ for \ } m-\e.
\end{equation}

(iii) Now we show that we also find \eqref{verandering1} when we calculate \eqref{i+1}. By \ref{1deel},a constituent $(\xi',\eta')$ with maximal $a$-value is of the form $(\xi',\eta')=(\xi,\eta)\cup(x,y)$ with $x+y=t$. It remains to show that $x$ and $y$ are as in \eqref{verandering1}.


We have to add $t$ squares to $(\xi,\eta)$, such that $a_m$ is maximized. Consider for a moment all $m$-symbols of $(\a,\beta)$ occurring in $\mbox{Ind}((t)\otimes (\xi,\eta))$ to have lengths $(n+m,m)$. Suppose that $(\a,\beta)$ occurs in the induction. Let the $m$-symbol of $(\a,\beta)$ (resp. $(\xi,\eta)$) have entries $x'_1\geq x'_2\geq \dots \geq x'_{2n+m}$ (resp. $x_1\geq x_2 \geq \dots \geq x_{2n+m}$). Let $s_{(\a,\beta)}$ be such that $x'_i=x_i$ for $i=1,2,\dots,s_{(\a,\beta)}$ and $x'_{s_{(\a,\beta)}+1} \neq x_{s_{(\a,\beta)}+1}$. Then, by \ref{n=dom}, if $(\xi',\eta')$ occurs in \eqref{i+1} then $s_{(\xi',\eta')}$ is maximal among $\{ s_{(\a,\beta)} \mid (\a,\beta) \mbox{ occurs in } \mbox{Ind}((t)\otimes (\xi,\eta)\}$; i.e., the truncated induction changes the parts of $(\xi,\eta)$ whose entries are ``as small as possible''.

Since we write the parts of $\xi$ and $\eta$ in increasing order, we therefore locate a pair $(\xi_i,\eta_j)$ such that $i+j$ is maximal among the pairs $\{(i,j)\mid (\xi_i,\eta_j) \neq (\xi'_i,\eta'_j), a_m(\xi',\eta') \mbox{ is maximal}\}$.

For example, if $t>\xi_{m+r}+\eta_r$, then we can not perform the induction without enlarging either $\xi_{m+r}$ or $\eta_r$ and we find $(\xi_i,\eta_j)=(\xi_{r+m},\eta_r)$. In general, there are two analogous situations possible. 

(iii-a) Suppose we have a sequence
\begin{equation}
\label{mog1}
 e_m(\eta_{f+1}) > e_m(\xi_k) > e_m(\xi_{k-1}) \dots > e_m(\xi_l) \geq e_m(\eta_f)
\end{equation}
which is as long as possible, i.e. $e_m(\xi_{k+1})>e_m(\eta_{f+1})$ and $e_m(\xi_{l-1})<e_m(\eta_f)$, in which \[ \xi_k+\eta_{f+1} \geq t > \xi_l+\eta_f.\] Some of these parts may not exist, in which case we delete the relation in which they occur. We check if, and how, we need to change parts $(\xi_i,\eta_j)$ into $(\xi'_i,\eta'_j)$. First we consider $\eta_f$ and one of the parts $\xi_l,\dots,\xi_k$. We do not need to change $\xi_k$ if $t-(\xi_k+\eta_f) < e_m(\xi_k)-e_m(\eta_f)=\xi_k-\eta_f+2(k-f)$. In that case we check if we have to change $\xi_{k-1}$, etc. This means that we have to look for the smallest $p\in\{0,1,\dots,k-l\}$ such that
\begin{equation}
\label{i}
 t \geq  2\xi_{k-p}+2(k-p-f).
\end{equation}
In fact, since $t$ is odd, this will be a strict inequality. Suppose for a moment that such $p$ exists. Then we calculate how the parts of $\xi_{k-p}$ and $\eta_f$ should be changed, in order to maximize the $a$-value. This means that we distribute the $t-\xi_{k-p}-\eta_f$ squares of the strip, which need to be added to $\xi_{k-p}$ and $\eta_f$, in such a way that we increase the entries $e_m(\xi_{k-p})$ and $e_m(\eta_f)$ minimally. This produces the equations $\xi'_{k-p}=\xi_{k-p}+x, \eta'_f=\eta_f+y, x+y=t-\xi_{k-p}-\eta_f$ and $e_m(\xi'_{k-p})=e_m(\eta'_f)$. One checks that these equations lead to the following:
\begin{equation}
\label{verandering2}
 (i)-\left\{ \begin{array}{r} \xi'_{k-p}=\lfloor \frac{t}{2} \rfloor-(k-p-f)  \\ \eta'_{f}= \lceil \frac{t}{2} \rceil+(k-p-f)  \end{array} \right.  {\rm \ or\ } (ii)-\left\{ \begin{array}{r} \xi'_{k-p}=\lceil \frac{t}{2} \rceil -(k-p-f)\\ \eta'_f= \lfloor \frac{t}{2} \rfloor+(k-p-f).  \end{array} \right.
\end{equation}

It remains to show that $a=k-p$, $b=f$, for then we recover
\eqref{verandering1}. This can be shown as follows: if
$p>0$, then \eqref{i} and minimality of $p$ imply that $\xi_{k-p}
\leq \lfloor\frac{t}{2}\rfloor -(k-p-f)$ and $\xi_{k-p+1} \geq
\lfloor \frac{t}{2} \rfloor -(k-p-f)$. Also $e_m(\eta_f)<
e_m(\xi_{k-p})$ implies $\eta_f<\lfloor\frac{t}{2}\rfloor+(k-p-f)$
and $e_m(\eta_{f+1})> e_m(\xi_{k-p+1})$ implies $\eta_{f+1} >
\lceil\frac{t}{2}\rceil +(k-p-f)$. Thus by \eqref{induktievw} both
\eqref{verandering2}(i)-(ii) are allowed, unless
$\xi_{k-p+1}=\lfloor\frac{t}{2}\rfloor-(k-p-f)$.
If $\xi_{k-p+1}\neq 0$, then it starts on a square
containing $|k-p-f+1|$, hence ends on
$\lfloor\frac{t}{2}\rfloor$, which is impossible. Therefore
$\xi_{k-p+1}$ ends on more than $\lfloor\frac{t}{2}\rfloor$ and
$\xi_{k-p}$ on less. On the other hand, if $\xi_{k-p+1}=0$, then it follows that the block of length $\eta_f$ starts on more than$\floor{\frac{t}{2}}$, from which we see that a strip of length $t$ can not be fitted into the $m$-tableau $T_m(\mu^{(i)})$, which is a contradiction. Again because of \ref{volgorde}, we see from \eqref{mog1} that if $\eta_{f+1}\neq 0$, then $j(\eta_{f+1})>\lfloor\frac{t}{2}\rfloor$, and $j(\eta_f)<\lfloor\frac{t}{2}\rfloor$.
This means that, as we wanted to show, $a=k-p$ and $b=f$.
The upshot is that we find two 2-partitions with maximal $a$-value, {\it because} the
strip of length $t$ fits into the diagram. The
new entries in the $m$-symbol are consecutive. 

Notice that $\xi_{k-p+1}\geq \lfloor \frac{t}{2} \rfloor -(k-p-f)$ implies
that $e_m(\xi_{k-p+1})\geq \lceil \frac{t}{2}\rceil+(k-p+f)+1$.
In equality holds, then the corresponding block of length $\xi_{k-p+1}$ ends on
$\lfloor\frac{t}{2}\rfloor$, and the strip of length $t$ can not be fitted
into $T_m(\mu^{(i)})$. In that situation only \eqref{verandering2}-(i) is possible, since
\eqref{verandering2}-(ii) would yield two consecutive entries in
the top row of the $m$-symbol, which is impossible. Then the
entries of $\xi'_{k-p}, \eta'_f$ and $\xi'_{k-p+1}$ form three
consecutive numbers, i.e., an interval in the sense of Lusztig
(see \cite{carter}) is formed.

We still have to check the case $p=0$. Then $t>2\xi_k+2(k-f)$, so
$\xi_k\leq\lfloor\frac{t}{2}\rfloor-(k-f)$, and $t\leq
\xi_k+\eta_{f+1}$ then implies
$\eta_{f+1}\geq\lceil\frac{t}{2}\rceil+(k-f)$. Since
$e_m(\eta_f)\leq e_m(\xi_k)$, we also find
$\eta_f\leq\lfloor\frac{t}{2}\rfloor+(k-f)$ and finally
$e_m(\xi_{k+1})\geq e_m(\eta_{f+1})$ implies $\xi_{k+1}\geq
\lceil\frac{t}{2}\rceil-(k-f)$. This shows that indeed both
 \eqref{verandering2}(i)-(ii) are allowed. As in the case $p>0$,
here also $a=k-p$ and $b=f$, which can be shown similarly.

The last possibility in situation \eqref{mog1} is that there is no
$p\in\{0,1,\dots,k-l\}$ satisfying \eqref{i}. Then we have $t\leq
2\xi_l+2(l-f)$, so $\xi_l \geq \lceil\frac{t}{2}\rceil-(l-f)$. In
this case we do not have to change $\xi_l$, but we have to change
$\eta_f$, since $t> \xi_l+\eta_f$. From $e_m(\eta_f)>e_m(\xi_{l-1})$ it follows that $\xi_{l-1}\leq \ceil{\frac{t}{2}}-(l-f-1)$. This means (as it did in \eqref{i}) that we have to change the pair
$(\xi_{l-1},\eta_f)$. Calculating the extension which maximizes
the $a$-value yields \eqref{verandering1} with $a=l-1,b=f$. Since
$\eta_f$ resp. $\xi_l$ start on a square containing $|l-f-1|$
resp. $|l-f|$, the blocks of length $\eta_f$ resp. $\xi_l$ end on
at least, resp. at most $\lfloor\frac{t}{2}\rfloor$. Therefore
indeed $a=l-1,b=f$, and both extensions are allowed since there is
no block ending on $\lfloor\frac{t}{2}\rfloor$.

(iii-b) The other possibility is that we have a sequence of the form
\begin{equation}
\label{mog2}
e_m(\xi_{f+1})>e_m(\eta_k)>\dots e_m(\eta_l)> e_m(\xi_f)
\end{equation}
such that $e_m(\eta_{k+1})>e_m(\xi_{f+1}), e_m(\eta_{l-1})<e_m(\xi_f)$, and \[\xi_{f+1}+\eta_k\geq t > \eta_l+\xi_f.\] In this situation, we can do computations analogous to the ones of situation \eqref{mog1}. We first check if there is $p\in \{0,1,\dots,k-l\}$ such that
\begin{equation}
\label{i2}
t \geq 2\eta_{k-p}+2(k-p-f)
\end{equation}
If such $p$ exists, take the smallest such. Then we find that $a=f,b=k-p$ and that both
\eqref{verandering1}(i)-(ii) are allowed, where we need to know
that there is no strip ending on $\lfloor\frac{t}{2}\rfloor$ if
$i>0$. Otherwise, we have to change $(\xi_f,\eta_{l-1})$.  Again
the $a$-value is always maximized by \eqref{verandering1}.

(iv) Suppose that in (iii) we have changed the parts $\xi_{i}$ and $\eta_{j}$. We check when it is necessary to increase the entries of $\xi_{i-1}$ and $\eta_{j-1}$. Since we have placed $t-\xi_{i}-\eta_{j}$ squares, there now remain $\xi_{i}+\eta_{j}$ to be placed. It is easy to see that \eqref{induktievw} implies that all $\xi'_k=\xi_{k+1}$ $ (k=1,\dots,i-1)$ and $\eta'_k=\eta_{k+1}$ $ (k=1,\dots,j-1)$, and finally $\xi'_0=\xi_1,\eta'_0=\eta_1$. This is also what we found in (ii). We conclude that indeed $(\xi',\eta')=(\xi,\eta) \cup (x,y)$ with $x,y$ as in \eqref{verandering1}.

(v) Suppose that $(\xi,\eta), (\tilde{\xi},\tilde{\eta})$ are two different members of $\Sigma_m(W_0c_{L^{(i)}})$,  that $(\xi',\eta')$ appears in $\trm((t)\otimes (\xi,\eta))$, and similarly for $(\tilde{\xi}',\tilde{\eta}')$. Then it is easy to see that $(\xi',\eta') \neq (\tilde{\xi}',\tilde{\eta}')$.
Therefore we find that $|\Sigma_m(W_0c_{L^{(i+1)}})|=2|\Sigma_m(W_0c_{L^{(i)}})|$.
On the other hand, it is also clear that
\[ | C_m(W_0c_{L^{(i+1)}})|=2| C_m(W_0c_{L^{(i)}})|,\]
thus the required bijection between $C_m(W_0c_{L^{(i+1)}})$ and $\Sigma_m(W_0c_{L^{(i+1)}})$ follows. Moreover, we have already seen that $\Sigma_m(W_0c_{L^{(i+1)}}) \subset [\S_{m\pm\e}(\mu^{i+1})]_m$, and it is not hard to see that $|[\S_{m\pm\e}(\mu^{(i+1)})]_m|=2|[\S_{m\pm\e}(\mu^{(i)})]_m|$, thus it follows from the induction hypothesis that indeed $\Sigma_m(W_0c_{L^{(i+1)}})= [\S_{m\pm\e}(\mu^{(i+1)})]_m$. This proves the induction step, and therefore also the Proposition if $m$ is integer.

(vi) If $m$ is half-integer, we can apply the same type of
reasoning. We briefly state the results in this case. Now $t$ must
be even. The strip corresponding to the $A$-factor contains the
entries
$(\frac{t-1}{2},\frac{t-3}{2},\dots,\frac{1}{2},\frac{1}{2},\dots,\frac{t-3}{2},\frac{t-1}{2})$.
Again, we consider $(\xi,\eta) \in \S_{m\pm\e}(\mu^{(i)})$. Suppose that blocks $\xi_{a+1}$ and $\eta_{b+1}$ end on more than $\frac{t-1}{2}$, and $\xi_a$ and
$\eta_b$ end on less. Then the strip of the $A$-factor will
account for the entries $(\xi'_a,\eta'_b)$, taking the place of
$(\xi_a,\eta_b)$. Since the corner square of the strip contains
this time the value $|a-b-\frac{1}{2}|$, we find the following
analogue of \eqref{verandering1}:
\begin{equation}
\label{verandering3}
 \left\{ \begin{array}{r} \xi'_a= \frac{t}{2}-(a-b)+1  \\ \eta'_b= \frac{t}{2}+(a-b)-1  \end{array} \right. {\rm for \ } m+\e {\rm \ or\ } \left\{ \begin{array}{r} \xi'_a= \frac{t}{2} -(a-b)\\ \eta'_b= \frac{t}{2}+(a-b)  \end{array} \right. {\rm \ for \ } m-\e.
\end{equation}
This time the entries in the $m$-symbol are $e_m(\xi_i)=\xi_i+2i$ and $e_m(\eta_i)=\eta_i+2i+1$, which means that if we are in situation \eqref{mog1}, this time we need to find the smallest $p \in \{0,1,\dots,k-l\}$ such that
\begin{equation}
\label{i3}
t \geq 2\xi_{k-p}+2(k-p-f)-1.
\end{equation}
In situation \eqref{mog2}, we have to find the smallest $p$ for which
\begin{equation}
\label{i4}
t \geq 2\eta_{k-p}+2(k-p-f)+1.
\end{equation}
Then it turns out, that in the same way as for integer $m$, we find the same
$(\xi',\eta')$ as obtained in \eqref{verandering3}, where the fact that both 2-partitions occur in the induction is due to the fact that the strip of length $t$ can be inserted into $T_m(\nu)$.\qed

\subsubsection{}\label{algemeen} Now we can consider the general situation. Let $L \in \L_m(n)$ be a
residual coset. Then we choose $M \in \L_m(n)$ of minimal dimension among those $M$ such that $W_0c_M  \in C_m(W_0c_L)$. Suppose $M$ has type $(\rho,\mu)$, and $sp_m(\mu)=(\kappa,\nu)$. Putting $\a=\rho \cup \kappa$, it
follows that $L$ is of type $(\a,\nu)$. It follows
from the minimality of $\mbox{dim}(M)$  that a strip
of length $\rho_i$ can not be incorporated into $T_m(\nu)$. This
means that:
\begin{lemma} Let $ m \in \half \Z$. A strip of length $t$ can not be inserted into $T_m(\mu)$ if and only if one of the following holds:
\begin{itemize}
\item[(i)]{$m$ is integer and $t$ is even, or $m$ is not integer and $t$ is odd.}
\item[(ii)]{One of the extremities of $T_m(\mu)$ is equal to $\floor{\frac{t}{2}}$;}
\item[(iii-a)]{$m >0$, $l(\mu)<m$ and $\lfloor\frac{t}{2}\rfloor<m-l(\mu)$;}
\item[(iii-b)]{$m <0$, $l(\mu')<m$ and $\lfloor\frac{t}{2}\rfloor<m-l(\mu')$;}

\end{itemize}
\end{lemma}

\pf: This is straightforward.\qed

Now suppose that we consider a
residual subspace of type $A_{t-1} \times (B_s,\mu)$. Then, as in
Proposition \ref{tellenkorrespondenten}, we look for the
constituent in ${\rm Ind}((t)\otimes\S_m(\mu))$ with maximal
$a$-value. Then we observe that such constituent is unique
precisely when there is no confluence:

\begin{prop}
\label{geenkonf}
Let $\mu, \rho$ be partitions, and let $(\xi,\eta) \in \Sigma_m(W_0c_L)$ where $L \in \L_m(n)$ has type $sp_m(\mu)=(\kappa,\nu)$. Suppose that strips of length $\rho_i$ can not be inserted into $T_m(\nu)$. Then ${\rm tr}_m-{\rm Ind}(\rm{triv}_\rho \otimes (\xi,\eta))=(\xi',\eta')$ for some $(\xi',\eta') \in \P_{n,2}$ with $n=|\mu|+|\rho|$.
\end{prop}

\pf: We adopt the notation of Proposition \ref{tellenkorrespondenten} and its proof. We suppose that $m \in \Z_{\geq 0}$, the other cases being analogous. First we remark that if a strip of length $t$ can not be inserted into $T_m(\nu)$, then it can also not be inserted into $T_m(\mu)$. Let $t$ be odd, such that $\floor{\frac{t}{2}}<m-l(\nu)$. If $\kappa=\emptyset$ then $\floor{\frac{t}{2}}<m-l(\mu)$, else there is an extremity of $\mu$ equal to $\floor{\frac{t}{2}}$.

Suppose that $\kappa=\emptyset$. In this case, it is clear that $a_m$ is maximal for $(\mu,\rho)$. The reader may verify that this is also the 2-partition obtained in the proof of \ref{tellenkorrespondenten}. Therefore if $\kappa=\emptyset$, we are done.

Next, suppose that $\kappa\neq\emptyset$. We can now go through the proof of \ref{tellenkorrespondenten}. The only difference is that there may be pairs of equal entries in the $m$-symbol of $(\xi,\eta)$. Therefore, in part (iii), the sequence \eqref{mog1} now reads
\begin{equation}
\label{mog3}
 e_m(\eta_{f+1}) \geq e_m(\xi_k) > e_m(\xi_{k-1}) \dots > e_m(\xi_l) \geq e_m(\eta_f)
\end{equation}
which is as long as possible, i.e. $e_m(\xi_{k+1})>e_m(\eta_{f+1})$ and $e_m(\xi_{l-1})<e_m(\eta_f)$, in which \begin{equation}\label{vw1} \xi_k+\eta_{f+1} \geq t > \xi_l+\eta_f.\end{equation} The other possibility is that one finds a sequence
\begin{equation}
\label{mog4}
e_m(\xi_{f+1})\geq e_m(\eta_k)>\dots e_m(\eta_l)\geq e_m(\xi_f)
\end{equation}
such that $e_m(\eta_{k+1})>e_m(\xi_{f+1}), e_m(\eta_{l-1})<e_m(\xi_f)$, and \begin{equation}\label{vw2}\xi_{f+1}+\eta_k\geq t > \eta_l+\xi_f.\end{equation}

Notice that there is a unique sequence \eqref{mog3} satisfying \eqref{vw1} or else a unique sequence \eqref{mog4} satisfying \eqref{vw2}.

In each case, one obtains two potential Springer correspondents given by \eqref{verandering1} for some $a,b$. However,  if the strip of length $t$ can not be fitted into the diagram, then only one of them is a 2-partition, and v.v.. For even $t$, the two possibilities for $(\xi'_i,\eta'_j)$ of \ref{tellenkorrespondenten} reduce to one. Therefore in this case as well, we are done. \qed

Now we show that indeed we find the desired result:

\begin{thm}\label{tellen2}
 Let $L \in \L_m(n)$ be a residual subspace. Then the set $\Sigma_m(W_0c_L)$ forms a similarity class in $\hat{W}_0/\sim_m$, and there is a bijection
\[ C_m(W_0c_L) \longleftrightarrow \Sigma_m(W_0c_L).\]

\end{thm}

\pf: Let $L$ have type $A_\a \times (B_l,\nu)$ and let
$\rho,\kappa,\nu$ be as described in the introduction to this
subsection. If $\rho=\emptyset$ then we have already treated this
case in Proposition \ref{tellenkorrespondenten}. Suppose therefore
that we have $\rho \neq \emptyset$. Let $(\xi,\eta)$ be a Springer correspondent of the residual subspace of type $(\kappa,\mu)$. By  Proposition \ref{geenkonf} above, $(\xi,\eta)$ leads (by inducing $\mbox{triv}_\rho \otimes (\xi,\eta)$) to exactly one Springer correspondent $(\xi',\eta')$ of $W_0c_L$. It remains to check that the thus obtained $(\xi',\eta')$ are mutually different and in the same similarity class. This is however not difficult to see. 

The fact that $\Sigma_m(W_0c_L)$ is a full similarity class in bijection with $C_m(W_0c_L)$ follows from
the corresponding fact in Proposition \ref{tellenkorrespondenten}. Indeed, let $M$ be the residual subspace of type $(\kappa,\nu)$, then $\Sigma_m(W_0c_M)$ is in bijection with $C_m(W_0c_M)$ by \ref{tellenkorrespondenten}. But $|\Sigma_m(W_0c_L)|=|\Sigma_m(W_0c_M)|$ by Proposition \ref{geenkonf}, and $|C_m(W_0c_L)|=| C_m(W_0c_M)|$ as well by definition of $M$ and $L$.
\qed
\subsubsection{Deformed symbols} Thus far, the $m$-symbols only apply to parameters $k_2=mk_1$ with $k_1 \neq 0$ and $m \in \half\Z$. However, we can use them as well to apply to any choice of parameters with $k_1 \neq 0$ if we deform them slightly. Choose $m \in \half \Z$ and let $\e >0$ be very small. Let $(\xi,\eta) \in \P_{n,2}$. Then we define the $m+\e$-symbol (resp. the $m-\e$-symbol) of $(\xi,\eta)$ to be the $m$-symbol of $(\xi,\eta-\e)$ (resp. the $m$-symbol of $(\xi,\eta+\e)$. By $\eta\pm\e$ we mean the partition whose parts are $\eta_i\pm\e$.
Let $a_{m\pm\e}(\xi,\eta)=a_m(\xi,\eta\mp\e)$ be corresponding $a$-value, with corresponding truncated induction $\mbox{tr}_{m\pm\e}-\mbox{Ind}$. Then it is not hard to show that 
\[ \Sigma_{m\pm\e}(L)=\mbox{tr}_{m\pm\e}-\mbox{Ind}_{S_\l \times W_0(B_l)}^{W_0(B_n)}(\mbox{triv}_\l \otimes \S_{m\pm\e}(\mu))\]
consists of exactly one element $(\xi,\eta) \in \Sigma_m(W_0c_L)$, and moreover that we obtain bijections 
\[ \bigcup_{W_0c_{L'} \in C_m(W_0c_L)}\Sigma_{m+\e}(W_0c_{L'}) \longleftrightarrow \Sigma_m(W_0c_L) \longleftrightarrow  \bigcup_{W_0c_{L'} \in C_m(W_0c_L)}\Sigma_{m+\e}(W_0c_{L'}) .\]
Finally, one checks that we can now define $\Sigma_x(L)$ for any $x \in \R$, since the Springer correspondence is constant between two consecutive $m \in \half\Z$:
\begin{prop}
Let $m,m' \in \half\Z_{\geq 0}$ with $m-m'=\half$. Let $L$ be a generically residual subspace of type $(\l,\mu)$. Then
\[ \Sigma_{m-\e}(W_0c_L)=\Sigma_{m'+\e}(W_0c_L). \]
Therefore, for $m>x>m'$, we put $\Sigma_x(W_0c_L)=\Sigma_{m-\e}(W_0c_L)=\Sigma_{m'+\e}(W_0c_L)$.
\end{prop}

\pf: This boils down to going through the proof of \ref{tellenkorrespondenten} and checking that for both $m$ and $m'$ the same of the two possibilities in \eqref{verandering2} or \eqref{verandering3} is obtained by $\mbox{tr}_{m-\e}-\mbox{Ind}$ and $\mbox{tr}_{m'+\e}-\mbox{Ind}$. For details, see \cite{proefschrift}.\qed

\subsection{Unipotent classes} In the classical cases, we have the bijection $\centra \leftrightarrow \U_m(n)$, such that the Springer correspondents of $W_0c_L$ are those of the corresponding unipotent class. We will show that we can define a set $\U_m(n)$ with the analogous properties for every $m \in \halfZ$. Furthermore, by also defining the generalizations $f^{BC}_m$ and $\phi_m$, we will show that our generalized Springer correspondence reduces to the classical one in case $m \in \{\half,1\}$.

Recall that in this case, for a residual point $c$ with jumps
$j_i$ the partition $\l$ of the corresponding unipotent class
satisfies $\l=(2j_i+1)$. We therefore calculate what the weight of
the corresponding partition for general $m$ is:

\begin{lemma}
\label{2j+1}
Let $k_2=mk_1$ with $m \in \halfZ$ and let $c=c(\mu,k,mk)$ be residual with jumps $\{j_i\}$. Define $\l$ to be the partition with parts $2j_i+1$. If $m$ is integer, then $|\l|=2n+m^2$, otherwise $|\l|=2n+m^2-\frac{1}{4}$.
\end{lemma}

\pf: Clearly it is sufficient to consider $m \geq 0$.

(i) Suppose that $m$ is integer. First consider $\mu_1 \leq \mu_2
\leq \dots \leq \mu_p$ with $p \leq m-1$. Then
$\S_m(\mu)=(\mu,-)$, and the jumps of $c$ are
$\{m+\mu_p-1,m+\mu_{p-1}-2,\dots,m+\mu_1-p,
m-p-1,m-p-2,\dots,1,0\}$. Therefore, indeed
\begin{eqnarray*} \sum_{i}2j_i+1 &=& \sum_{i=1}^l (2(m+\mu_i-i)+1) + \sum_{i=0}^{m-p-1} (2i+1) \\ &=& 2mp+2n -2\frac{p(p+1)}{2}-2\frac{(m-p-1)(m-l)}{2}+(m-p) \\ &=& m^2+2n \end{eqnarray*}

Now consider the general situation, where we may (and will) assume
that $T_m(\mu)$ is in standard position, i.e., that $\S_m(\mu)=(\xi(c),\eta(c))$ as in \ref{startpart}.
Hence, $\S_m(\mu)$ consists of $m-1$ horizontal blocks
$(\xi_{r+2},\dots,\xi_{r+m})=(\mu_{p-m+2},\dots,\mu_p)$, followed by a
series of alternating horizontal and vertical blocks. Then
$(\xi_{tr},\eta)$, with
$\xi_{tr}=(\xi_1,\xi_2,\dots,\xi_{r+1})$, is the $m$-tableau
of a residual point for $m=1$. Therefore we know that $\sum_{i=1}^{2r+1}
(2j_i+1)= 2(\sum_{i<p-m+2}\mu_i)+1$. Thus
\begin{eqnarray*} \sum_{i}(2j_i+1)&=& \sum_{i=1}^{2r+1}(2j_i+1)+\sum_{i=2}^{m}(2(i+\mu_{p-m+i}-1)+1)\\&=&2\sum_{i<p-m+2}\mu_i+1+2\sum_{i=2}^m\mu_{p-m+i}+m^2-1\\&=&m^2+2n.\end{eqnarray*}

(ii) If $m$ is half-integer, we find that if $l(\mu)=p \leq
m-\frac{1}{2}$ that the jumps are
$\{m+\mu_p-1,m+\mu_{p-1}-2,\dots,m+\mu_1-p,m-p-1,m-p-2,\dots,\frac{3}{2},\frac{1}{2}
\}$, which after a similar computation yields $\sum_{i}2j_i+1 =
m^2+2n-\frac{1}{4}$. The general case is treated analogous to
integer $m$ as well.\qed

Therefore we define
\begin{defn} Let $m \in \halfZ$.
If $m$ is integer we define
\[ {\mathcal U}_m(n):=\{ \l=(1^{r_1}2^{r_2}\dots) \vdash 2n+m^2 \mid r_i {\rm \ is\ even\ if\ }i {\rm \ is\ even\ and\ } \sum_{i {\rm \ odd}} (r_i {\rm \ mod\ }2) \geq |m| \}. \]
That is, we consider partitions of $2n+m^2$ in which even parts
have even multiplicity and which have at least $|m|$ odd parts with odd multiplicity.

For half-integer $m$ we define \begin{flushleft} ${\mathcal
U}_m(n):=\{ \l=(1^{r_1}2^{r_2}\dots) \vdash 2n+m^2-\frac{1}{4}
\mid $ \end{flushleft}  \begin{flushright}$ r_i {\rm \ is\ even\
if\ }i {\rm \ is\ odd\ and\ } \sum_{i {\rm \ even}} (r_i {\rm \
mod\ }2) \geq |m|-\frac{1}{2} \}.$ \end{flushright}

That is, we consider partitions of $2n+m^2-\frac{1}{4}$ in which
odd parts have even multiplicity, and which have at least
$|m|-\frac{1}{2}$ even parts with odd multiplicity.
\end{defn}
Notice that in the equal label cases $m=\half,1$ we indeed recover the
old $\mathcal{U}_1(n)$. 

We now define a map $f^{\mbox{\tiny{BC}}}_m:{\mathcal U}_m(n) \to \centra$, associating a central character of $\Hgrrcc$  to $\l \in {\mathcal U}_m(n)$.

\begin{defn} Let $k_2=mk_1 \neq 0$ with $m \in \half \Z$ and suppose $\l \in {\mathcal U}_m(n)$.

(i) If $m$ is integer, let $l_1<l_2<\dots<l_s$ be the parts that
occur an odd number of times in $\l$, and put $n_0:=\sum_jl_j$.
All $l_j$ are odd and $s \geq m$. Notice that $n_0 \equiv m(2)$.
Let $l=\frac{n_0-m^2}{2}$. By removing each $l_j$ once from $\l$,
we obtain a partition where every part occurs an even number of
times. If we remove each second part, we find a partition $d_1
\leq d_2 \leq \dots \leq d_r$. The associated residual subspace
$L$ is the one with \[ R_L= A_{d_1-1} \times A_{d_2-1} \times
\dots A_{d_r-1} \times B_{l},\] and the residual point in $B_{l}$
has jumps $\lfloor \frac{l_i}{2} \rfloor$. Then
$W_0c_L:=f^{\mbox{\tiny{BC}}}_m(\l)$.

(ii) If $m$ is half-integer, let again be $l_1<\dots< l_s$ be the
parts of $\l$ that occur an odd number of times. Now all $l_i$ are
even, and $s \geq m-\frac{1}{2}$. Let $n_0=\sum_jl_j$ and put
$l=\frac{n_0-m^2+1/4}{2}$. By removing each $l_j$ once from $\l$
we obtain a partition where each part has even multiplicity.
Remove each second part to obtain a partition $d_1\leq d_2 \leq
\dots \leq d_r$, then $L=f^{\mbox{\tiny{BC}}}_m(\l)$ has root
system \[ R_L= A_{d_1-1}\times A_{d_2-1} \times \dots \times
A_{d_r-1} \times B_{l},\] and the residual point in $B_{l}$ has
jumps $\frac{l_j}{2}$. Then $W_0c_L:=f^{\mbox{\tiny{BC}}}_m(\l)$.
\end{defn}

\begin{lemma}
The map $f^{BC}_m$ is a bijection, i.e. the set ${\mathcal U}_m(n)$ parametrizes the set of central characters of $\Hgrrcc$.
\end{lemma}

\pf: Consider a residual subspace $L$ with $R_L$ of type $A_\kappa
\times B_l$. The residual points in $B_l$ are characterized by
their jumps $j_i$. From \ref{aantaljumps} we know that the parts
$2j_i+1$ form a partition of at least $m$ (resp. at least
$m-\frac{1}{2}$) positive and distinct odd (resp. even) parts. By
adding two parts $\kappa_i$ for every $A$-factor, we obtain indeed
a partition $\l \in {\mathcal U}_m(n)$. Clearly
$f^{\mbox{\tiny{BC}}}_m(\l)=L$. Therefore $f^{\mbox{\tiny{BC}}}_m$
is a bijection between ${\mathcal U}_m(n)$ and the set of $W_0$-orbits
of residual subspaces. Finally we use the bijection $W_0L \to W_0c_L$. \qed

We next define the generalization $\phi_m$, associating to $\l \in \U_m(n)$ an equivalence class in $\hat{W}_0$. Later we show that this class is $\Sigma_m(f^{BC}_m(\l))$.

\begin{defn}\label{defphi}Let $m \in \half\Z_{\geq 0}$ and consider $\l \in {\mathcal U}_m(n)$. Its parts are arranged in
increasing order. If $m \notin \Z$ then we ensure that
$l(\l)=|m|-\half+2r$ for some $r \in \Z_{\geq 0}$ by putting $\l_1=0$ if necessary.

(i) For integer resp. half-integer $m$, replace in the last $|m|$
resp. $|m|-\frac{1}{2}$ parts of $\l$ a pair of even resp. odd
entries $(x,x)$ by $(x+1,x-1)$ to obtain an $n$-composition $\mu$.

(ii) Define $\mu^*_i=\mu_i+(i-1)$ for $i=1,2,\dots,2r$ and
$\mu^*_{2r+i}=\mu_{2n+i}+2r$ for $i=1,2,\dots,|m|$ if $m$ is
integer; and $\mu^*_i=\mu_i+(i-1)$ ($i=1,2,\dots,2r)$;
$\mu^*_{2r+i}=\mu_i+(2r-1)$ for $i=1,2,\dots,|m|-\frac{1}{2}$ if
$m$ is half-integer.

(iii) Form the 2-composition $(\xi^*,\eta^*)$ where $\eta^*$
contains the $\frac{\mu_i}{2}$ for even $\mu_i$, and $\xi^*$
contains the $\lfloor \frac{\mu_i}{2} \rfloor$ for odd $\mu_i$ (in
the order they appear in $\mu^*$).

(iv) Define the 2-composition $(\xi,\eta)$ by
$\xi_i=\xi^*_i-(i-1)$ and $\eta_i=\eta^*_i-(i-1)$. If $m$ is
half-integer, readjust the lengths of $\xi$ and $\eta$ to
$l(\xi)=l(\eta)+|m|+\frac{1}{2}$.

(v) Finally $\phi_m(\l)=[(\xi,\eta)]_m=\{ (\a,\beta) \in \P_{n,2} \mid (\a,\beta) \sim_m (\xi,\eta) \}$.

(vi) Put $\phi_{-m}(\l)=[(\eta,\xi)]_{-m}$ if
$\phi_m(\l)=[(\xi,\eta)]_m$.
\end{defn}

Notice that for $m \in \{\half,1\}$, we recover the map $\phi_m$ of section \ref{phi}. As in this case, the $m$-symbol of $(\xi,\eta)$ (of part (iv) in Definition \ref{defphi}) is
always increasing. This follows easily from the corresponding, known,
fact in the equal label cases. However, it may happen that
$(\xi,\eta) \notin \P_{n,2}$ since $\xi$ need not be a partition.

Finally, notice also that if $m=0$, we have $(\xi,\eta) \sim_0 (\eta,\xi)$,
so the definition is unambiguous.

We want to show that $\phi_m$ establishes a bijection between $\U_m(n)$ and $\hat{W}_0/\sim_m$.
Therefore we define the candidate inverse map $\psi_m=\phi_m^{-1}$:

\begin{defn} Let $k_2=mk_1$ and $m\in \halfZ_{\geq 0}$. For an equivalence class $\Sigma \in \hat{W}_0/\sim_m$, define $\l = \psi_m(\Sigma)$ as follows. Let $m \geq 0$.

(i) Let $(\xi,\eta)$ be the unique 2-{\rm composition}
$(\xi,\eta)$ such that $[(\xi,\eta)]_m \cap \hat{W}_0= \Sigma$ and
$(\xi,\eta)$ has increasing $m$-symbol. Let $l(\xi)=m+n$,
$l(\eta)=n$ if $m$ is integer or $l(\xi)=m-\half+n$, $l(\eta)=n$
if $m$ is half-integer. Notice that for $m \notin \Z$ this is
different from the symbol lengths!

(ii) Let $\xi_i^*=\xi_i+(i-1)$ and $\eta_i^*=\eta_i+(i-1)$. Note that
$\xi_i^*$ still is in general a composition rather then a
partition.

(iii) Let $\mu^*$ be the composition which consists of the parts
$2\xi_i+1$ and $2\eta_i$. We order the first $2n$ parts in
increasing order, and add the last $m$ (resp. $m-\half$) parts in
the order they appeared in $\xi^*$.

(iv) Now we construct from $\mu^*$ the $n$-composition $\mu$. If
$m$ is integer we define $\mu_i:=\mu_i^*-(i-1)$ for
$i=1,2,\dots,2n$ and $\mu_{2n+i}=\mu_{2n+i}^*-2n$ for
$i=1,2,\dots,|m|$. If $m$ is not integer then we define
$\mu_i=\mu_i^*-(i-1)$ for $i=1,2,\dots,2n$ and
$\mu_{2n+i}=\mu^*_{2n+i}-(2n-1)$ for $i=1,2,\dots,m-\half$.

(v) Finally we remove all zeroes from $\mu$, and replace any
entries in $\mu$ of the form $(x+1,x-1)$ by $(x,x)$. The end
result is $\l:=\psi_m(\Sigma)$.

(vi) Let $m \geq 0$ and let $\Sigma$ be an equivalence class in
$\hat{W}_0$ for $\sim_{-m}$. We define
$\psi_{-m}(\Sigma)=\psi_m(\Sigma \otimes (-,n))$.

\end{defn}

Notice that since $\sim_0$-equivalence classes are invariant under
multiplication with $(-,n)$, the definition is unambiguous for
$m=0$. In the equal label cases, we recover the inverse of $\phi_m$ of section \ref{phi}.

\begin{lemma}
The map $\psi_m$ is invariant under the shift operation for
symbols, and therefore really is a map on
$\bar{Z}_n^{2,0}(n+\ceil{m},n)$.
\end{lemma}

\pf: This is clear, since if we add a zero to both $\a$ and $\beta$, we add a zero to $\xi$ and $\eta$, which means we add two zeroes to $\l$. \qed

To prove that $\phi_m$ and $\psi_m$ realize a bijection between $\U_m(n)$ and $\hat{W}_0/\sim_m$, it is sufficient to show that $\psi_m(\Sigma) \in \U_m(n)$, since both $\phi_m$ and $\psi_m$ are easily seen to be injective. To show this, we first prove a Lemma.
\begin{lemma}
\label{psilemma} Let $(\a,\beta)$ be a 2-partition of $n$ and let
$m \in \half \Z_{\geq 0}$.

(i) Suppose that in the $m$-symbol of $(\a,\beta)$, one has
$e_m(\a_{n+i+1})>e_m(\beta_n) \geq e_m(\a_{n+i})$ for some $1 \leq
i \leq m$. Let $\a_{tr}=(\a_1,\dots,\a_{n+i})$, and suppose that
$\psi_m(\a,\beta)=\l$. If $m \in \Z$ then
$\psi_i(\a_{tr},\beta)=\l_{tr}=(\l_1,\l_2,\dots,\l_{2n+i})$.
Otherwise,
$\psi_{i-1/2}(\a_{tr},\beta)=\l_{tr}=(\l_1,\l_2,\dots,\l_{2n+i})$.

 (ii-a) Suppose that $m \in \Z$ and $e_m(\beta_n) \geq e_m(\a_{n+m})$.
Let $\a_{tr}=\a-\a_{n+m}$ and $\beta_{tr}=\beta-\beta_n$. Suppose
$\psi_m(\a,\beta)=\l$. Then
$\psi_m(\a_{tr},\beta_{tr})=\l_{tr}=(\l_1,\dots,\l_{2(n-1)+m})$.

(ii-b) Suppose that $m \notin \Z$ and $e_m(\beta_n) \geq
e_m(\a_{n+m+\half})$. Let $\a_{tr}=\a-\a_{n+m+\half}$ and
$\beta_{tr}=\beta-\beta_n$. Suppose $\psi_m(\a,\beta)=\l$. Then
$\psi_m(\a_{tr},\beta_{tr})=\l_{tr}=(\l_1,\dots,\l_{2(n-1)+m-\half})$.

\end{lemma}

\pf: (i) This is clear, since if we denote the corresponding
2-compositions with increasing $m$-symbols by
$(\xi,\eta)\sim_m(\a,\beta)$ and
$(\xi_{tr},\eta)\sim_i(\a_{tr},\beta)$ (resp. $\sim_{i-\half}$),
then $\xi=\xi_{tr}\cup (\xi_{n+i+1},\dots,\xi_{n+m})$.

(ii-a) Suppose that $(\xi,\eta) \sim_m (\a,\beta)$ and
$(\xi_{tr},\eta_{tr}) \sim_m (\a_{tr},\beta_{tr})$ have increasing
$m$-symbols. Then
$\xi_{tr}=(\xi_1,\dots,\xi_n,\eta_n-2,\xi_{n+1}-2,\dots,\xi_{n-2+m}-2)$
and $\eta_{tr}=\eta-\eta_n$. Therefore \[ \xi^*_{tr}=(\xi^*_1,\dots,\xi^*_n,\eta^*_n-1,\xi^*_{n+1}-1,\dots,\xi^*_{n+m-2}-1)
{\rm \ and\ } \eta^*_{tr}=\eta^*-\eta^*_n.\]
Therefore, it follows that $\mu_{tr,i}=\mu_i$ for $i\in
\{1,2,\dots,2(n-1),2n+2,2n+3,\dots,m+2(n-2)\}$. It remains to
compare $\mu_i$ to $\mu_{tr,i}$ for $i=2n-1,2n,2n+1$. Both are
derived from $\xi_n,\xi_{n+1}$ and $\eta_n$.

We therefore have the following possibilities. Either we have
$e_m(\xi_n)=e_m(\eta_n)$, or we have
$e_m(\xi_n)<e_m(\eta_n)<e_m(\xi_{n+1})$, or we have
$e_m(\eta_n)=e_m(\xi_{n+1})$. Let us first assume that
$e_m(\xi_n)=e_m(\eta_n)$, so $\xi_n=\eta_n$,
$\mu^*_{2n-1}=2\eta^*_n,\mu^*_{2n}=2\xi^*_n+1$, and
$\mu_{2n-1}=2\eta^*_n-(2n-2)=\mu_{2n}=2\xi^*_n+1-(2n-1)$. In
$\l=\psi_m(\a,\beta)$ these entries remain the same. Now lets us
see what we find for $\psi_m(\a_{tr},\beta_{tr})$. First we
calculate
$\mu^*_{tr}=(\mu^*_1,\dots,\mu^*_{2n-2},2\xi^*_n+1,2\eta^*_n-1,2\xi^*_{n+1}-1,\dots,2\xi^*_{n+m-2}-1)$,
and so
$\mu_{tr}=(\mu_1,\dots,\mu_{2n-2},2\xi^*_n+1-(2n-2),2\eta^*_n-1-(2n-2),\mu_{2n+1},\mu_{2n+m-2})$.
Therefore, in $\l_{tr}$, the two entries
$(2\xi^*_n+1-(2n-2),2\eta^*_n-1-(2n-2))$ are replaced by the pair
$(2\xi^*_n-(2n-2),2\eta^*_n-(2n-2))$, which are the same entries
obtained in $\l$. This proves the claim in case
$e_m(\xi_n)=e_m(\eta_n)$. Next we consider the second possibility
$e_m(\xi_n)<e_m(\eta_n)<e_m(\xi_{n+1})$. Then in $\mu^*$ we get
$\mu^*_{2n-1}=2\xi^*_{n}+1$ and $\mu^*_{2n}=2\eta_n^*$, whereas
$\mu^*_{tr}=(\mu^*_1,\dots,\mu^*_{2n-1},2\eta_n^*-1,2\xi^*_{n+1}-1,\dots,2\xi^*_{n+m-2}-1)$,
and so
$\mu_{tr}=(\mu_1,\dots,\mu_{2n-1},2\eta^*_n-1-(2n-1),2\xi^*_{n+1}-1-2(n-1),\dots,2\xi^*_{n+m-2}-1-2(n-1))=(\mu_1,\dots,\mu_{2n+m-2})$.
Hence also $\l_{tr}=(\l_1,\dots,\l_{2n+m-2})$. The last
possibility is that $e_m(\eta_n)=e_m(\xi_{n+1})$, so
$\eta_n=\xi_{n+1}+2$ and $\eta^*_n=\xi^*_{n+1}+1$. In that case we
find in $\mu^*$ that $\mu^*_{2n}=2\xi^*_{n+1}+1$ and
$\mu^*_{2n+1}=2\eta^*_n$, and so in $\mu$ we get
$\mu_{2n}=\mu_{2n+1}=2\eta^*_n-2n$. On the other hand,
$\mu^*_{tr}=(\mu^*_1,\dots,\mu^*_{2n-1},2\eta^*_n-1,2\xi^*_{n+1}-1,2\xi^*_{n+2}-1,\dots,2\xi^*_{2n+m-2}-1)$
and so
$\mu_{tr}=(\mu_1,\dots,\mu_{2n-1},2\eta^*_n-2n+1,2\eta^*_n-2n-1,\mu_{2n+2},\dots,\mu_{2n+m-2})$.
This means that in $\l_{tr}$ we find
$\l_{tr,2n}=\l_{tr,2n+1}=2\eta^*_n-2n$, which are also the
corresponding entries of $\l$. So we find that in every case, the
statement of the lemma holds.

(ii-b) Let $(\tilde{\xi},\tilde{\eta}) \sim_m (\a,\beta)$ and
$(\tilde{\xi}_{tr},\tilde{\eta}_{tr}) \sim_m (\a_{tr},\beta_{tr})$
have increasing $m$-symbol. Then $l(\tilde{\xi})=n+m+\frac{1}{2}$
and $l(\eta)=n$. We readjust the lengths by putting
$\xi_i=\tilde{\xi}_{i+1}$ for $i=1,2,\dots,n+m-\half$, and taking
$\eta=\tilde{\eta}$. Similarly we let
$\xi_{tr,i}=\tilde{\xi}_{tr,i+1}$ for $i=1,2,\dots,n+m-3/2$, and
$\eta_{tr}=\tilde{\eta}_{tr}$. Then we have
\[ \xi_{tr}=(\xi_1,\dots,\xi_{n-1},\eta_n-1,\xi_{n}-2,\dots,\xi_{n+m-\half-2}-2),\eta_{tr}=(\eta_1,\dots,\eta_{n-1}).\]
For the partitions $\xi^*$ and $\eta^*$ this implies
\[\xi^*_{tr}=(\xi^*_1,\dots,\xi^*_{n-1},\eta^*_n-1,\xi^*_{n}-1,\dots,\xi_{n+m-\half-2}-1),\eta^*_{tr}=(\eta^*_1,\dots,\eta^*_{n-1}).\]
The length of $\l$ is $2n+m+\frac{1}{2}$. Therefore,
$\mu_i=\mu^*_i-(i-1)$ for $i=1,2,\dots,2n$ and
$\mu_i=\mu^*_i-(2n-1)$ for $i=2n+1,2n+2,\dots,2n+m+\frac{1}{2}$.
On the other hand, since $l(\l_{tr})=2(n-1)+m+\half$, we have
$\mu_{tr,i}=\mu^*_{tr,i}-(i-1)$ for $i=1,2,\dots,2(n-1)$ and
$\mu_{tr,i}=\mu^*_{tr,i}-(2(n-1)-1)$ for
$i=2n-1,2n,\dots,2(n-1)+m+\half$. To compare $\l$ to $\l_{tr}$ we
need to compare the triple of entries
$(\l_{2n-2},\l_{2n-1},\l_{2n})$ to
$(\l_{tr,2n-2},\l_{tr,2n-1},\l_{tr,2n})$.

Suppose first that $e_m(\tilde{\xi}_n)=e_m(\tilde{\eta}_n)$. Then
$\tilde{\xi}_n=\eta_n+1$, i.e. $\xi_{n-1}=\eta_n+1$, and
$\xi^*_{n-1}=\eta^*_n$. Thus
$\mu^*_{2n-2}=2\eta_n^*,\mu^*_{2n-1}=2\xi_{n-1}^*+1,\mu^*_{2n}=2\xi^*_n+1$.
Thus, $\mu_{2n-2}=2\eta_n^*-(2n-3),
\mu_{2n-1}=2\xi^*_{n-1}+1-(2n-2),\mu_{2n}=2\xi^*_n+1-(2n-1)$ and
$\l_{2n-2}=\mu_{2n-2},\l_{2n-1}=\mu_{2n-1},\l_{2n}=\mu_{2n}$. On
the other hand, for $\l_{tr}$ we find that
$\mu^*_{tr,2n-2}=2\xi^*_{n-1}+1,\mu^*_{tr,2n-1}=2\eta^*_n-1$ and
$\mu^*_{2n}=2\xi^*_n-1$, hence
$\mu_{tr,2n-2}=2\xi^*_{n-1}+1-(2n-3),
\mu_{tr,2n-1}=2\eta_n^*-1-(2n-3)$ and
$\mu_{tr,2n}=2\xi^*_n-1-(2n-3)$. This means that
$\l_{tr,2n-2}=\l_{tr,2n-1}=2\xi^*_{n-1}-(2n-3)=\l_{2n-2}=\l_{2n-1}$
and also $\l_{tr,2n}=2\xi^*_n-1-(2n-3)=\l_{2n}$.

The second possibility is that
$e_m(\tilde{\xi}_n)<e_m(\tilde{\eta}_n)<e_m(\tilde{\xi}_{n+1})$.
In this case, $\xi_{n-1}<\eta_n+1$ and $\eta_n<\xi_n+1$. This
means that $\xi^*_{n-1}<\eta^*_n$ and $\eta^*_n<\xi^*_n+1$. In
$\mu^*$ we thus get
$\mu^*_{2n-2}=2\xi^*_{n-1}+1,\mu^*_{2n-1}=2\eta^*_n,
\mu^*_{2n}=2\xi^*_n+1$. Hence,
$\mu_{2n-2}=2\xi^*_{n-1}+1-(2n-3),\mu_{2n-1}=2\eta^*_n-(2n-2),\mu_{2n}=2\xi^*_n+1-(2n-1)$.
But by the expressions for corresponding parts of $\mu_{tr}$ that
we have found in the preceding case, we see that $\mu_{tr}$ and
$\mu$ agree.

Finally consider the case where
$e_m(\tilde{\eta}_n)=e_m(\tilde{\xi}_{n+1})$. Then
$\eta_n=\xi_n+1$ and $\eta^*_n=\xi^*_n+1$. Therefore, we find
parts $\mu^*_{2n-2}=2\xi_{n-1}^*+1,\mu^*_{2n-1}=2\xi^*_{n}+1$ and
$\mu^*_{2n}=2\eta_n^*$. Thus,
$\mu_{2n-2}=2\xi_{n-1}^*+1-(2n-3),\mu_{2n-1}=2\xi^*_n+1-(2n-2)$
and $\mu_{2n}=2\eta^*_n-(2n-1)$. On the other hand, we now have
$\mu_{tr,2n-1}=\mu_{tr,2n}+2$, hence
$\l_{tr,2n-1}=\l_{tr,2n}=\l_{2n-1}=\l_{2n}=2\eta^*_n-(2n-1)$.
Hence for half-integer $m$ as well the statement of the Lemma is
now proven.\qed

\begin{prop}
\label{totaal}
Suppose $(\a,\beta)\in\P_{n,2}$. Then $\psi_m(\a,\beta)\in {\mathcal U}_m(n)$.
\end{prop}

\pf: We may suppose that $m \geq 0$.

(i) Let $m \in \Z$. Let the 2-composition $(\xi,\eta)$ of $n$ be
the one described in part (i) of the definition of $\psi_m$.
Following the procedure $\psi_m$, we find that, since we consider
$\xi$ as having $n+m$ parts and $\eta$ as having $n$ parts:
\[ |\psi_m(\a,\beta)|= \sum_{i=1}^{n+m} 2(\xi_i+(i-1))+1 + \sum_{i=1}^n 2(\eta_i+(i-1)) - \sum_{k=1}^{2n-1} k - 2nm = 2n+m^2, \]
so indeed $\l:=\psi_m(\a,\beta) \vdash 2n+m^2$. It remains to be
shown that even parts in $\l$ occur with even multiplicity and
that there are at least $m$ odd parts in $\l$ with odd
multiplicity. Notice that although $(\xi,\eta)$ is in general a
2-composition rather than a 2-partition, $\eta$ is a partition, as
is $\xi_{tr}=(\xi_1,\dots,\xi_{n+1})$. On this part
$(\xi_{tr},\eta)$, $\psi_m$ acts as $\psi_1$ and therefore the
first $2n+1$ parts of $\l$ form a partition where even parts occur
with even multiplicity. Also, since $\mu^*_{2n+i}=2\xi^*_{n+i}+1$
and $\mu_{2n+i}=\mu^*_{2n+i}-2n$, the last $m$ parts of $\mu$ are
all odd. In the final step $\mu \to \l$, we create only pairs of
equal parts, which shows that indeed in $\psi_m(\a,\beta)$ even
parts have even multiplicity.

Now let us show that $\psi_m(\a,\beta)$ contains at least $m$ odd
parts with odd multiplicity. We have seen that the parts
$\mu_{2n+1},\dots,\mu_{2n+m}$ are all odd, and at least equal to
$\mu_{2n}$. If $\l=\mu$, then we are done. If not, then there are
parts $\xi_{n+i}=\xi_{n+i+1}+2$ in $\xi$, for some $1\leq i \leq
m-1$, giving rise to entries in the $m$-symbol of $(\xi,\eta)$
which are equal. Hence we have to show that for any such pairs of
entries, there are two odd parts among $\mu_1,\dots,\mu_{2n}$ with
odd multiplicity. Suppose that in $(\a,\beta)$ one has
$e_m(\beta_n) =e_m(\a_{n+i})$. In view of Lemma \ref{psilemma}
(i), we may then assume that $i=m$, since otherwise the last $m-i$
entries of $\l$ will be different odd entries occurring with
multiplicity one. But then we may use \ref{psilemma}(ii) to remove
the last parts of both $\a$ and $\beta$. We can continue doing
this until we reach a 2-partition of the form $(\kappa,\nu)$ with
increasing $m$-symbol, for which $\psi_m(\kappa,\nu)$ will indeed
have at least $m$ odd parts with odd multiplicity, but then so
does $\psi_m(\a,\beta)$.

(ii) Suppose that $m \notin \Z$. Let $(\xi,\eta)$ be the
2-composition whose $m$-symbol is similar to the one of
$(\a,\beta)$ and increasing. Then, by definition of $\psi_m$ we
have
\begin{eqnarray*}
|\psi_m(\a,\beta)|&=&\sum_{i=1}^{n+m-\half}(2(\a_i+(i-1))+1)+2\sum_{i=1}^n2(\beta_i+(i-1))-\sum_{i=1}^{2n}(i-1)\\&&-(2n-1)(m-\half)\\&=&2n+m^2-\frac{1}{4}.\end{eqnarray*}
Again, it remains to be seen that odd parts have even multiplicity
and that there are at least $m-1/2$ distinct even parts with odd
multiplicity. Define $\xi_{tr}=(\xi_1,\dots,\xi_n)$. Then $\psi_m$
acts as $\psi_{\half}$ on $(\xi_{tr},\eta)$. Therefore the first
$2n$ parts of $\mu$ (in the definition of $\psi_m$) form a
partition in which odd parts have even multiplicity. The last
$m-\half$ parts of $\mu$ are of the form
$\mu_{2n+i}=\mu^*_{2n+i}-(2n-1)$. Since the last $m-\half$ parts
of $\mu^*$ are of the form $2\xi^*_{n+i}+1$, these parts of $\mu$
are all even. It remains to be seen that for every replacement
$(x,x)\to (x+1,x-1)$ in $\mu \to \l$, there are two distinct even
parts among $\mu_1,\dots,\mu_{2n}$ with even multiplicity. As in
the case where $m \in \Z$, this follows from Lemma \ref{psilemma}.
\qed

Now we can show that the maps we have defined lead to the desired result:
\begin{thm}
\label{commdiag}For $k_2=mk_1$ with $k_1 \neq 0$ and $m \in \half \Z$,  we have,

(i)  bijections
\[ \centra \longleftrightarrow \{\Sigma_m(W_0c_L) \mid W_0L \in \L_m(n)\} \longleftrightarrow \hat{W}_0/\sim_m;\]

(ii) the identity \[ (f_m^{\mbox{\tiny{BC}}} \circ \psi_m \circ \Sigma_m)(W_0c_L)=W_0c_L \mbox{ for all } W_0c_L \in \centra. \]
\end{thm}
\pf: (i) In view of Theorem \ref{tellen2}, we need to show that if $W_0c_L \neq W_0c_{L'}$, then $\Sigma_m(W_0c_L) \neq \Sigma_m(W_0c_{L'})$. However, since $\U_m(n) \leftrightarrow \centra$ via $f_m^{BC}$, and $\U_m(n) \leftrightarrow \hat{W}_0/\sim_m$ via $\phi_m$, the statement follows from (ii).

(ii) Consider $L \in \L_m(n)$ of type $A_\kappa \times (B_l,\nu)$. We prove the Proposition by induction on $l(\kappa)$.

(a) If $l(\kappa)=0$, then we consider a residual point
$c=c(\nu,k,mk)$, whose $m$-tableau $T_m(\nu)$ we may assume to be in
standard position, i.e., of the form $\J_m(\xi(c),\eta(c))$ as in \ref{fig:begin}. We have to show that
\[f^{\mbox{\tiny{BC}}}_m(\psi_m(\S_m(\nu)))=W_0c(\nu,k,mk).\] Suppose that $J(c)=\{j_1,j_2,\dots,j_{2r+m}\}$. Then
$\S_m(\nu)=(\xi,\eta)$ where
$\xi=(j_1,j_2,\dots,\allowbreak j_{2r+1},j_{2r+2}-1,j_{2r+3}-2,\dots,j_{2r+m}-(m-1))$
and $\eta=(j_2+1,j_4+1,\dots,j_{2r}+1)$. This leads to
$\xi^*=(j_1,j_3+1,j_5+2,\dots,j_{2r+1}+r,j_{2r+2}+r,\dots,j_{2r+m}+r)$
and $\eta^*=(j_2+1,j_4+2,\dots,j_{2r}+r)$, hence
$\mu^*=(2j_1+1,2j_2+2,2j_3+3,2j_4+4,\dots,2j_{2r+1}+2r+1,2j_{2r+2}+2r+1,\dots,2j_{2r+m}+2r+1)$,
and thus $\l=\psi_m(\xi,\eta)$ indeed consists of the distinct odd
parts $(2j_1+1,2j_2+1,\dots,2j_{2r+m}+1)$.

(b) Now suppose that the proposition is
true for every $L' \in \L_m(n)$ of type $(\a,\beta)$ with $l(\a)<l(\kappa)$, i.e. in particular for $M$ of type $(\kappa-\kappa_{l(\kappa)},\nu)$. Let the $m$-symbols of $\Sigma_m(W_0c_M)$ have size $(r+m,r)$. Put $t=\kappa_{l(\kappa)}$, and consider \[\Sigma_m(W_0c_L)= \trm_{S_t \times W_0(B_{n-t})}^{W_0(B_n)}((t)\otimes \Sigma_m(W_0c_M)).\]
Let $(\a,\beta) \in \Sigma_m(W_0c_M)$ and $(\a',\beta') \in \trm((t)\otimes (\a,\beta)) \subset \Sigma_m(W_0c_L)$. Recall from (the proof of) \ref{tellenkorrespondenten} that  $(\a',\beta')=(\a,\beta) \cup (t_1,t_2)$ for some decomposition $t=t_1+t_2$. Let $(\a'_i,\beta'_j)$ be as in part (iii) of the proof of \ref{tellenkorrespondenten}. By the induction hypothesis, we assume that $f^{BC}_m(\psi_m(\a,\beta))=W_0c_M$.

Now consider $\psi_m(\a,\beta)$. Let $(\xi',\eta') \sim_m (\a',\beta')$ and $(\xi,\eta) \sim_m (\a,\beta)$ be the 2-compositions with increasing symbol. There are three possibilities for the position of $e_m(\a'_i),e_m(\beta'_j)$ in the $m$-symbol of $(\xi',\eta')$: 
either they are found in position $(\xi'_i,\eta'_i)$ for some $0\leq
i \leq r$, either in position $(\xi'_{i+1},\eta'_i)$ for some $0\leq
i \leq r$, or in position $(\xi'_i,\xi'_{i+1})$ for some $i>r$. This
last possibility does not occur in the equal label cases.

Let $\psi_m(\xi,\eta)=\l$ and $\psi_m(\xi',\eta')=\l'$. If $\xi'_i$ gives rise to part $\l'_j$ we put $j:=\psi_m(\xi'_i)$, and similarly for $\eta'$. In each of the three cases one easily checks that $_psi_m(\xi',\eta')=\psi_m(\xi,\eta) \cup \psi_m(\xi'_i) \cup \psi_m(\eta'_j)$. Therefore, it is sufficient to show that $\psi_m(\xi'_i)=\psi_m(\eta'_j)=t$.

 Let us
consider these possibilities one by one. In the
first case, if $t$ is odd, we have
$(\xi'_i,\eta'_i)=(\lfloor\frac{t}{2}\rfloor,\lceil\frac{t}{2}\rceil)$.
This implies that $\xi^{'*}_i=\lfloor\frac{t}{2}\rfloor+(i-1)$ and
$\eta^{'*}_i=\lceil\frac{t}{2}\rceil+i-1$, hence in $\mu^*$ we find
the parts
$\mu^*_{2i-1}=2\lfloor\frac{t}{2}\rfloor+2(i-1)+1=t+2(i-1)$ and
$\mu^*_{2i}=2\lceil\frac{t}{2}\rceil+2(i-1)=t+2i-1$. This means
that indeed $\mu_{2i-1}=\mu^*_{2i-1}-(2i-2)=t$ and
$\mu_{2i}=\mu^*_{2i}-(2i-1)=t$. If $t$ is even, then 
$(\xi'_i,\eta'_i)=(\frac{t}{2},\frac{t}{2}),
(\xi^{'*}_i,\eta^{'*}_i)=(\frac{t}{2}+i-1,\frac{t}{2}+i-1)$, and so
$\mu^*_{2i-1}=t+2i-2$, $\mu^*_{2i}=t+2i-1$ which again yields
$\mu_{2i-1}=\mu_{2i}=t$. This proves the claim. 

The second case where $(\xi'_{i+1},\eta'_i)=(\lfloor\frac{t}{2}\rfloor,\lceil\frac{t}{2}\rceil)$
or $(\frac{t}{2}-1,\frac{t}{2}+1)$ proceeds in the same way.

Finally consider the third possibility. This happens if in the
proof of \ref{tellenkorrespondenten}, we have found that $(\a',\beta')=(\a,\beta)\cup(\a'_i,\beta'_r)$ (for some $i \geq m+r+1$), where $(\a'_i,\beta'_r)=(\lceil\frac{t}{2}\rceil-(i-r),\lfloor\frac{t}{2}\rfloor+(i-r))$,
with symbol entries $e_m(\a'_i)=\lceil\frac{t}{2}\rceil+i+r, e_m(\beta'_r)=\lfloor\frac{t}{2}\rfloor+i+r$. In $(\xi',\eta')$, these entries take position $(i-1,i)$ in the top row, hence $\xi'_{i-1}=\lfloor\frac{t}{2}\rfloor+i+r-2(i-1)=\lfloor\frac{t}{2}\rfloor-i+r+2$
and
$\xi'_i=\lceil\frac{t}{2}\rceil+i+r-2i=\lceil\frac{t}{2}\rceil-i+r$.
This implies $\xi^{'*}_{i-1}=\lfloor\frac{t}{2}\rfloor+r+1$ and
$\xi^{'*}_i=\lceil\frac{t}{2}\rceil+r$. In $\mu^*$, this yield two consecutive
parts $t+2r$, hence in $\mu$  two consecutive parts $t$. If
$t$ is even, then we find after rearranging the $m$-symbol in
increasing order that $e_m(\xi'_{i-1})=e_m(\xi'_i)=\frac{t}{2}+i+r$,
hence $\xi'_{i-1}=\frac{t}{2}-i+r+2,\xi'_i=\frac{t}{2}-i+r$, so
$\xi^*_{i-1}=\frac{t}{2}+r+1$ and $\xi^*_i=\frac{t}{2}+r$. This
means that in $\mu^*$ we find two consecutive parts $t+2r+3$ and
$t+2r+1$, and in $\mu$ we thus get $t+1$ and $t-1$. Therefore
again $\l$ contains 2 new parts $t$ and indeed
$f^{\mbox{\tiny{BC}}}_m$ finds a factor $A_{t-1}$. Since one
easily shows that the other parts of $\l$ have remained the same,
this proves the induction step.

(c) If $m$ is half-integer, the proof is analogous. \qed

\begin{cor} If $m\in\{\half,1\}$, and the central character $W_0c_L$ for $L \in \L_m(n)$ corresponds to the unipotent class $C_\l$, then $\Sigma_m(W_0c_L)=\Sigma_m(C_\l)$. Thus, we retrieve the classical Springer correspondence.
\end{cor}

\pf: For residual points, this has been remarked already. The Corollary follows since the maps $f_m^{BC}, \phi_m, \psi_m$ which we have defined for general $m$, reduce to the classical ones in the equal label cases.\qed

We can now refine \ref{conj} to include also the central character of the modules $M_\chi^m$.
\begin{conj}\label{con2}
In the notation of Conjecture \ref{conj}, the central character of the irreducible tempered $\Hgr$-module $M^m_\chi$ is equal to $W_0c_L$ if and only if $\chi \in \Sigma_m(W_0c_L)$.
\end{conj}
\begin{rem}
For $B_3$ and $B_4$, and special parameters, these conjectures have been shown to be true, thanks to the explicit calculations of both the appropriate Green functions by Gunter Malle, and the determination of $\Hgrrcc$.
\end{rem}

We can thus reformulate the Conjectures \ref{conj} and \ref{con2} as follows.

\begin{conj}\label{finalconj}
Let $\Hgr$ be the graded Hecke algebra associated to the root system of type $B_n$ with parameters $k_1,k_2$ such that $k_1 \neq 0$ and $k_2=mk_1$ for $m \in \half \Z$. Then the set $\Hgrrcc$ has the following description. Let $\U_m(n)$ be the set of partitions generalizing the unipotent conjugacy classes, and let $\Sigma_m(\l)=[\phi_m(\l)]_m$ be the set of Springer correspondents of $\l \in \U_m(n)$.

Then $\Hgrrcc$ has a set of representing modules
\[ \{ M^m_{\l;(\xi,\eta)} \mid \l \in \U_m(n),\ (\xi,\eta) \in \Sigma_m(\l) \}\] 
where this indexation is uniquely characterized by the requirements
\begin{itemize}
\item The central character of $M^m_{\l;(\xi,\eta)}$ is $W_0c_L$, where $L$ is the residual subspace such that $\Sigma_m(W_0c_L)=\Sigma_m(\l)$, given by $W_0c_L=f^{BC}_m(\l)$ ;
\item The $W_0$-module $M^m_{\l,(\xi,\eta)}$ contains the irreducible $W_0$-character $\chi_{(\xi,\eta)}$. 
\end{itemize}

Moreover, given this parametrization, the multiplicity of $(\xi,\eta)$ in the $W_0$-module $M^m_{\l;(\xi,\eta)}$ is equal to one.

In the notation of Conjecture \ref{conj}, we have $M^m_A=M^m_{\l;A}$. Thus, the $W_0$-decomposition of $M^m_{\l,(\xi,\eta)}$ is given by the Green functions $P^m_{BA}$ ($A,B \in \P_{n,2}$). It therefore inherits a grading, with respect to which the irreducible character $\chi_{(\xi,\eta)}$ occurs in the top degree of $M^m_{\l,(\xi,\eta)}$. This degree is equal to $a_m(\xi,\eta)$.

\end{conj}

In this formulation, the modules $H^m_{\l,(\xi,\eta)}$ are the analogue of the Springer modules $H(\B_u)_\rho$, where $\l$ replaces $u$ and $(\xi,\eta)$ replaces $\rho$. Indeed, in the classical case, the character $\rho$ can be read off the symbol of $(\xi,\eta)$.

\section{Examples}

We end by giving some examples of the conjectures and the constructions in the previous sections.
\subsubsection{Splitting and joining}
\begin{figure}[htb]
\begin{center}
{\includegraphics[angle=0,scale=0.27]{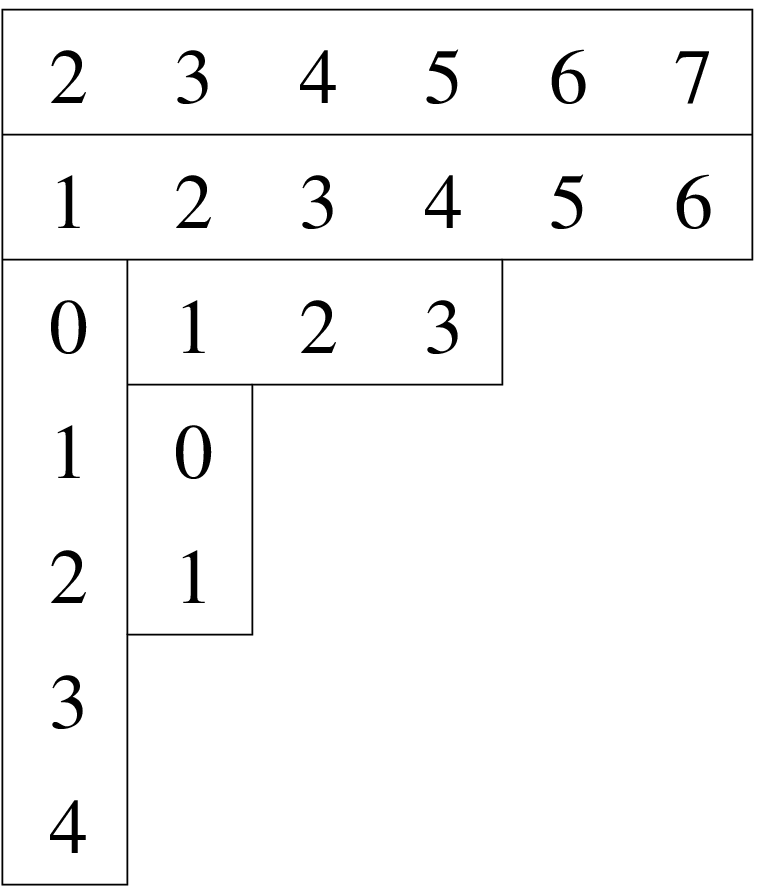}}
\end{center}
\caption{$\J_2((366),(25))$}
\label{(0366)(25)}
\end{figure}

Let $n=22,m=2$ and let $c$ be the residual point with coordinates $(7,6,6,5,5,4,4,4,3,3,3,3,2,2,2,2,1,1,1,1,0,0)k$, then $J(c)=(0,1,3,\allowbreak 4,6,7)$. Therefore the set $\Sigma_m(W_0c)$ is, according to \eqref{sigmarespunt}, equal to the $m$-equivalence class of
\begin{equation} \label{startsymbool}
 \left( \begin{array}{lllllll} 0& &3& &10& &12\\
                                  &2& &7&  & &   \end{array} \right)
\end{equation}
The 2-symbol \eqref{startsymbool} is the symbol of $(366,25)$. There are ${6 \choose 2}$ possible $2$-symbols with these entries, i.e., there are 15 generically residual points
$c(\l_i,k_1,k_2)$ ($i=0,1,\dots,14$) which, when $k_2=2k_1$,
become equal to $c$. The partitions $\l_i$ can be determined by
applying $\J_2$ to all double partitions whose $2$-symbols is similar to \eqref{startsymbool}. We construct first
$\l(c)$, obtained from \eqref{startsymbool}. It is obtained as $\l(c)=\J_2((366,25))$ and shown in Figure \ref{(0366)(25)}.
Now let us demonstrate how we can, given an arbitrary $\l \in
\P_2(c)$, reach $\l(c)$ through a series of flips. We take
$\l=(2^83^2)$. We need to rearrange the $2$-symbol of $\S_2(\l)$ into \eqref{startsymbool}, by interchanging the first entry which is bigger than
its successor with the unique entry in the other row with which
it can be interchanged (the entries in each row must remain increasing). In the Figure on page \pageref{voorbeeld}, we carry
this out. Observe that each permutation in the symbol
corresponds to a flip of one part of $\xi$ to a part of $\eta$,
until we reach the initial partition $\l(c)$. Notice also that we
have at every step 4 horizontal and 2 vertical extremities.

\subsubsection{Calculation of Springer correspondents} 
\begin{figure}[h]
\begin{center}
{\includegraphics[angle=0,scale=0.27]{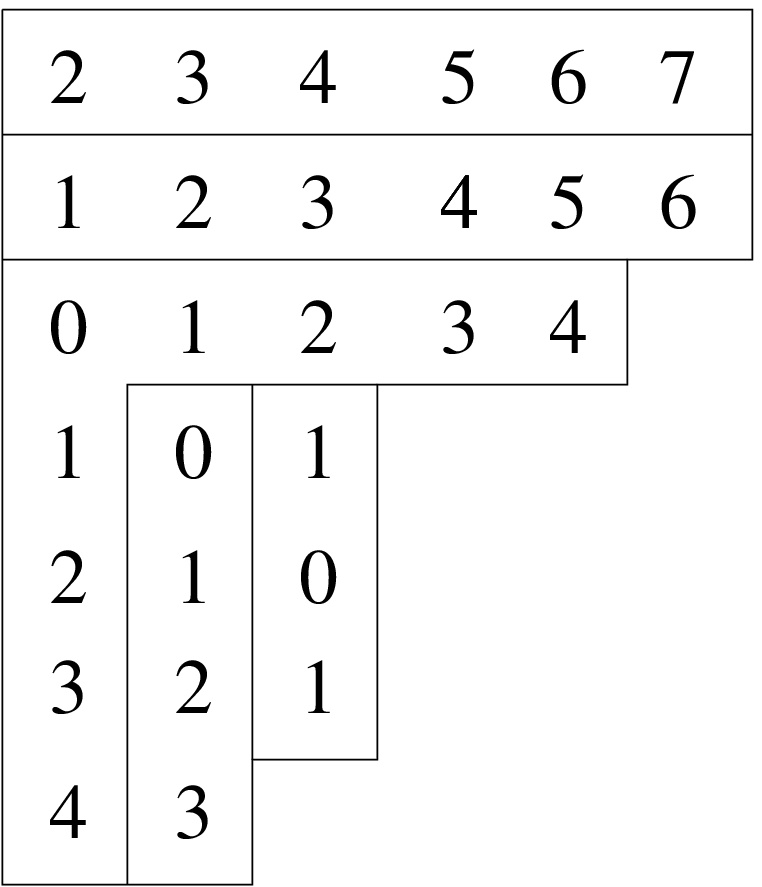}}
\caption{$T_2(\mu)$}
\end{center}
\label{fig:mu}
\end{figure}

Let $n=20,m=2$ and consider the residual subspace $L$ of type $A_2 \times A_6 \times A_8 \times B_{16},(\nu)$ where $\nu=(1^46^2)$. We will describe the sets $C_m(W_0c_L)$ and $\Sigma_m(W_0c_L)$, as well as the partition $\l \in \U_m(n)$ corresponding to $L$. We start with constructing $\Sigma_m(W_0c_L)$ by using the definition of truncated induction. We observe by drawing $T_2(\nu)$ that one can insert strips of length 3 and 9 into $T_2(\nu)$ to obtain a partition $\mu$. Then $T_2(\mu)$ is as is Figure \ref{fig:mu}.

We calculate $\Sigma_m(W_0c_L)$. Put $\kappa=(39), \a=(7), \rho=\a\cup\kappa=(379)$ as in \ref{algemeen}. We adopt the notation of Proposition \ref{tellenkorrespondenten}. We start with 
\[ \trm_{S_9 \times W_0(B_{16})}^{W_0(B_{25})}((9) \otimes (66,4)).\]
The 2-symbol of $(\xi,\eta)=(66,4)$, written with $l(\xi)=4,l(\eta)=2$ is
\[ \left( \begin{array}{lllllll}0&&2&&10&&12 \\ &0&&6&&& \end{array} \right) \]
Recall that $\xi=(\xi_0,\dots,\xi_3)$ and $\eta=(\eta_0,\eta_1)$. With $t=9$, we have
\[ e_2(\xi_3)>e_2(\xi_2)>e_2(\eta_1)>e_2(\xi_1), \]
and $\xi_2+\eta_1=10>t>\xi_1+\eta_1=4$. Thus, we are in situation \eqref{mog2} with $k=1,l=1,f=1$. Since $t>2\eta_1$, it follows that \eqref{verandering2} applies with $p=0$ and we obtain $(\xi'_1,\eta'_1)=(4,5)$ or $(5,4)$. The resulting 2-partitions $(0466,45)$ and $(0566,44)$ have $2$-symbols
\[ \left( \begin{array}{lllllll}0&&6&&10&&12 \\ &4&&7&&& \end{array} \right) \mbox{ and }  \left( \begin{array}{lllllll}0&&7&&10&&12 \\ &4&&6&&& \end{array} \right). \]
Notice that they are indeed similar. We proceed with the second induction, choosing $(0466,45)$ (the outcome does not depend on this choice), i.e. we calculate 
\[  \trm_{S_7 \times W_0(B_{25})}^{W_0(B_{32})}((7) \otimes (466,45)).\]
We put $t=7$ and consider the 2-symbol of $(\xi,\eta)$, where we have put $l(\xi)=5,l(\eta)=3$ and both indexations starting with zero:
\[ \left( \begin{array}{lllllllll}0&&2&&8&&12&&14 \\ &0&&6&&9& \end{array} \right). \]
Then we have 
\[ e_2(\xi_4)>e_2(\xi_3)>e_2(\eta_2)>e_2(\xi_2)>e_2(\eta_1)>e_1(\xi_1)>e_2(\eta_0)\geq e_2(\xi_0),\]
while $\xi_2+\eta_1=8>t>\xi_1+\eta_1=4$. Thus, we are in situation \eqref{mog2} with $f=1,k=l=1$. We have $t<2\eta_1$, hence we find $(\xi'_1,\eta'_0)=(3,4)$ or $(2,5)$. However, the latter is impossible since $\eta_1=4$. This is an illustration of how the impossibility of fitting the strip of length seven into $T_2(\nu)$, implies that this induction does not change the number of Springer correspondents. Thus we find the 2-partition $(03466,445)$. 
Finally , we consider the induction 
\[  \trm_{S_3 \times W_0(B_{32})}^{W_0(B_{35})}((3) \otimes (3466,445)).\]
We start from $(\xi,\eta)=(003466,0445)$ with 2-symbol
\[ \left( \begin{array}{lllllllllll}0&&2&&7&&10&&14&&16 \\ &0&&6&&8&&11& \end{array} \right). \]
Observe also the interval 6,7,8 which has been formed due to the induction of the factor $A_6$ which could not be fitted into $T_2(\nu)$.
We have $t=3$,
\[ e_2(\xi_5)>e_2(\xi_4)>e_2(\eta_3)>e_2(\xi_3)>e_2(\eta_2)>e_1(\xi_2)>e_2(\eta_1)> e_2(\xi_1)>e_2(\eta_0)\geq e_2(\xi_0),\]
and $\xi_1+\eta_1=4>t>\xi_1+\eta_0$. Thus we are in situation \eqref{mog1} with $f=0,k=l=1$. We have $2\xi_1+2<t$ and hence $(\xi'_1,\eta'_0)=(0,3)$ or $(1,2)$. Indeed, both are possible. Thus we obtain the 2-partitions $(3466,3445)$ and $(13466,2445)$ with 2-symbols
\[ \left( \begin{array}{lllllllllll}0&&2&&7&&10&&14&&16 \\ &3&&6&&8&&11& \end{array} \right)\]and \[\left( \begin{array}{lllllllllll}0&&3&&7&&10&&14&&16 \\ &2&&6&&8&&11& \end{array} \right). \]

We conclude that $\Sigma_m(W_0c_L)=[(13466,2445)]_m$. Let us check that this is also what we find by applying $\phi_m$ to $\l$, the partition corresponding to $L$. This partition consists of the parts $2j_i+1$ where the $j_i$ are the jumps of $c(\nu,k,mk)$, and the parts $(t,t)$ for every $A_{t-1}$, hence $\l=(1^2,3^2,7^3,9^2,13,15)$. Using the definition of $\phi_m$, one easily checks that $\phi_m(\l)=(13466,2445)$.

Finally we check the bijection between $\Sigma_m(W_0c_L)$ and $C_m(W_0c_L)$.
Since $l(\kappa)=2,m=2,r=1$, we have $|\Sigma_m(W_0c_L)|=2^2{2+2\cdot 1 \choose 1}=16$.

\begin{figure}[htb]
\begin{center}
{\includegraphics[angle=0,scale=0.27]{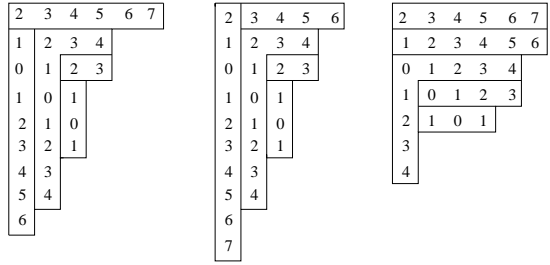}}
\caption{$\mu_1,\mu_2,\mu_3$}
\end{center}
\label{fig:others}
\end{figure}

The corresponding 2-symbols are found as follows. It is not hard to see that there are exactly three other partitions $\nu_1,\nu_2,\nu_3$ such that $\S_2(\nu_i)\sim_2 \S_2(\nu)$. We display them in Figure \ref{fig:others} where we insert also the blocks of lenght 3 and 9 to obtain $\mu_1,\mu_2,\mu_3$.

One easily checks that $\{\S_{2+\e}(\mu_i);i=1,2,3\}=\{(1236,268),(1234,26\ 10),(34566,4)\}$.
Applying
\[ \trm_{S_7\times W_0(B_{28})}^{W_0(B_{35})}((7)\otimes (\a,\beta))\]
where $(\a,\beta)=\S_{2\pm\e}(\mu_i)$ yields the other 2-partitions in $\Sigma_2(L)$. One easily checks that to $\mu_1$ the following 2-partitions are associated: we have
\[ \left( \begin{array}{lllllllllll}0&&3&&6&&8&&11&&16 \\ &2&&7&&10&&14& \end{array} \right),\]
and those resulting by interchanging 2 and 3, and/or 10 and 11; similarly from $\mu_2$ we find
\[ \left( \begin{array}{lllllllllll}0&&3&&6&&8&&11&&14 \\ &2&&7&&10&&16& \end{array} \right),\]
and those resulting by interchanging 2 and 3, and/or 10 and 11; and finally from $\mu_2$ we find
\[ \left( \begin{array}{lllllllllll}3&&6&&8&&11&&14&&16\\ &0&&2&&7&&10& \end{array} \right),\]
and those resulting by interchanging 2 and 3, and/or 10 and 11. One checks that these 16 2-partitions exhaust a similarity class.

\subsubsection{Calculation of Springer correspondents}
\begin{figure}[htb]
\begin{center}
{\includegraphics[angle=0,scale=0.27]{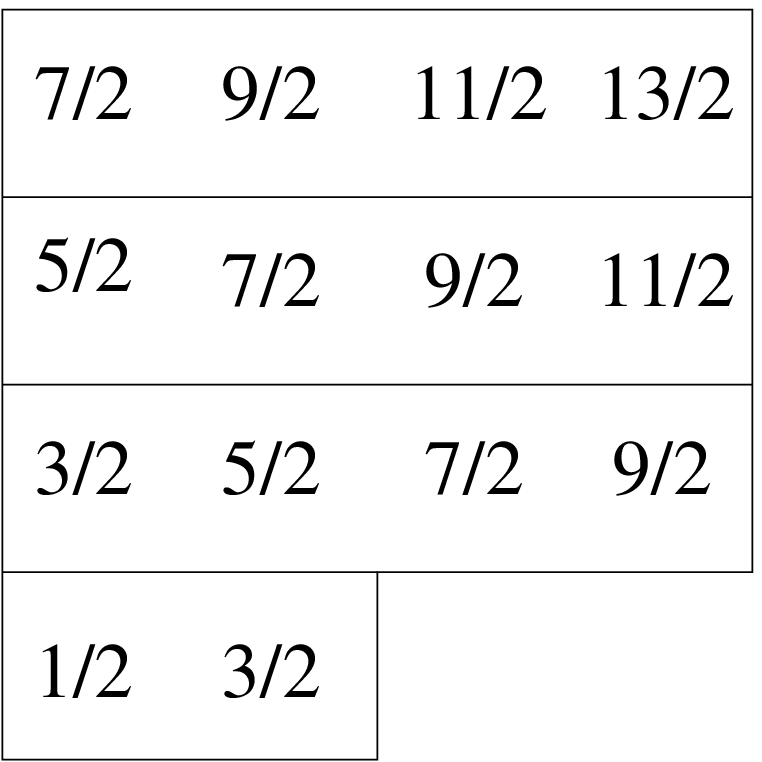}}
\caption{$T_{7/2}(\nu)$}
\end{center}
\label{fig:nu2}
\end{figure}

We give another example for non-integer $m$. Therefore, let $m=\frac{7}{2}$, let $n=35$ and consider $\U_m(n) \ni \l=(1\ 1\ 4\ 7\ 7\ 10\\\ 12\ 13\ 13\ 14)$. An application of $f^{BC}_m$ yields that $\l$ corresponds to $L \in \L_m(n)$ of type $A_6 \times A_{12} \times c$ where $c$ is a residual point for $B_{14}$ with jumps $J(c)=\{\frac{3}{2},\frac{9}{2},\frac{11}{2},\frac{13}{2}\}$. This is the point $c=c(\nu,k,mk)$ where $\nu$ is as in Figure \ref{fig:nu2}.
We apply $\phi_m$. According to the definition, first we need to add a zero to $\l$ and to replace the parts $(13\ 13)$ by $(14\ 12)$. Then we have $\mu=(0\ 1\ 1\ 4\ 7\ 7\ 10\ 12\ 14\ 12\ 14)$, $\mu^*=(0\ 2\ 3\ 7\ 11\ 12\ 16\ 19\ 21\ 19\ 21)$, $(\xi^*,\eta^*)=(1\ 3\ 5\ 9\ 10\ 9\ 10,0\ 1\ 6\ 8)$, and $(\xi,\eta)=(1236644,0045)$ whose lengths are readjusted to form $(\xi,\eta)=(123644,045)$ in the $\frac{7}{2}$-symbol, which reads:
\[ \left( \begin{array}{lllllllllllll}1&&4&&7&&12&&14&&14&&16\\ &1&&7&&10&&&&& \end{array} \right).\]
We replace it by
\[ \left( \begin{array}{lllllllllllll}1&&4&&7&&10&&12&&14&&16\\ &1&&7&&14&&&&& \end{array} \right),\]
which is the symbol of $(1234444,49)=\phi_m(\l)$.

We now check that we also find $\phi_m(\l)$ by computing $\Sigma_m(W_0c_L)$ using truncated induction. Thus, using again the notation of section \ref{algemeen} and Proposition \ref{tellenkorrespondenten}, we have $\nu=\mu$, $\kappa=\emptyset,\a=\rho=(1\ 7\ 13)$.
We start with computing
\[ \trm_{S_{13} \times W_0(B_{14})}^{W_0(B_{27})}((13) \otimes (2444,-)).\]
The symbol of $(02444,0)$ is
\[ \left( \begin{array}{lllllllll}0&&4&&8&&10&&12 \\ &1&&&&&&& \end{array} \right), \]
such that with $t=13$ we have
\[ e_m(\xi_4)>e_m(\xi_3)>e_m(\xi_2)>e_m(\xi_1)>e_m(\eta_0)>e_m(\xi_0),\]
and $t>\xi_1+\eta_0$.
That is, we are in situation \eqref{mog1} with $f=0,k=4,l=3$. In \eqref{i3}, we get $p=3$ since $2\xi_3+6-1=13$. According to \eqref{verandering3} we thus find $(\xi'_3,\eta'_0)=(4,9)$. We thus obtain the 2-partition $(24444,9)$. We proceed to calculate 
\[ \trm_{S_7 \times W_0(B_{27})}^{W_0(B_{34})}((7) \otimes (24444,9)).\]
The symbol of $(024444,09)$ is
\[ \left( \begin{array}{lllllllllll} 0&&4&&8&&10&&12&&14 \\ &1&&12&&&&&&& \end{array} \right), \]
such that 
\[ e_m(\xi_5)>e_m(\eta_1)\geq e_m(\xi_4)>e_m(\xi_3)>e_m(\xi_2)>e_m(\eta_0)>e_m(\xi_0),\]
and $\eta_1+\xi_4=13>7=t>\xi_1+\eta_0$. Analogous to the first induction we find this time that the $a$-value is maximal for $(234444,49)$. 
The last step, the induction
\[ \trm_{S_1 \times W_0(B_{34})}^{W_0(B_{35})}((1) \otimes (234444,49))\]
yields analogously that indeed \[ \Sigma_m(W_0c_L)=[(1234444,49)]_m=\phi_m(\l)\]
as expected.

\subsubsection{Calculation of tempered modules}
As an illustration to Conjectures \ref{conj} and \ref{con2} we consider a root system of type $B_3$. We consider the special parameters $k_2=mk_1$ with $m=2$.
 Then we order the double
partitions of $3$, refining the order given by the $a$-value, and
such that similarity classes form intervals. In this case, we have
\[ \begin{array}{cc} (\xi,\eta) & a_2(\xi,\eta) \\ \hline (-,111)
& 12 \\ (-,12) & 8\\ (1,11) & 6 \\ (-,3) & 5 \\ (1,2) &
4\end{array} \begin{array}{cc} (\xi,\eta) & a_2(\xi,\eta) \\
\hline (11,1) & 3 \\ (111,-) & 3 \\ (2,1)& 2 \\ (12,-) & 1 \\
(3,-) & 0
\end{array}
\]
and also $(11,1) \sim_2 (111,-)$; we therefore keep this ordering.
The matrix $P^2_{\vet{\a},\vet{\beta}}$ computed for us by G. Malle.
The results are as follows. For all $A=(\xi,\eta) \in \P_{n,2}$, we write for $l=1,2,\dots$ the entry $\sum_BP^{2,l}_{B,A}\chi_B\otimes\e$, i.e. the degree-$i$ component of $M_A^m$, viewed as $W_0$-module. The subspace $L$ is calculated as
$f_2^{\mbox{\tiny{BC}}}(\psi_2(\xi,\eta))$.
\[ \begin{array}{ccccl} (\xi,\eta) & \l=\psi_2(\xi,\eta)\in\U_2(3) & \mbox{type of }R_L,\mbox{res. pt.} & l & M_{(\xi,\eta)}^{2,l} \\
\hline  (-,111) & (1^73)& \emptyset & 3 & (3,-) \\ &&& 4& (2,1) \\&& & 5 & (1,2)+(12,-)\\
&&& 6 & (-,3)+(11,1)+(2,1) \\& && 7 & (1,11)+(1,2)+(12,-)\\&& & 8 &
(-,12)+(11,1)+(2,1) \\&& & 9 & (1,11)+(1,2)+(111,-) \\&& & 10&
(-,12)+(11,1)
\\ &&& 11 & (1,11) \\&& & 12 & (-,111) \vspace{0.1cm}\\ (-,12)& (1^32^23) & A_1 & 3 & (3,-)
\\ &&& 4 & (2,1)\\&& & 5 & (1,2)+(12,-)\\&&& 6 & (-,3)+(11,1)+(2,1) \\& && 7 &
(1,11)+(1,2) \\& && 8 & (-,12) \vspace{0.1cm}\\(1,11)& (1^55) & B_1& 2 &
(12,-)+(3,-)
\\ &&& 3 & (11,1)+(2,1) \\&& & 4& (1,2)+(111,-)+(12,-) \\&& & 5 & (11,1)+(2,1)\\&&& 6 &
(1,11) \vspace{0.1cm}\\ (-,3) & (13^3) & A_2&2 & (3,-) \\ &&& 3 & (2,1) \\&&
& 4 & (1,2)
\\&& & 5 & (-,3)\vspace{0.1cm}
 \\ (1,2) & (12^25)& A_1 \times B_1 & 2 & (12,-)+(3,-) \\ 
&&&3&(11,1)+(2,1)\\&&&4&(1,2)\vspace{0.1cm}\\(11,1)&(1^235)& B_2,(11) & 2 & (12,-) \\
&&& 3&(11,1)\vspace{0.1cm} \\ (111,-) & (1^235)& B_2,(11) & 3 & (111,-) \vspace{0.1cm}\\
(2,1) & (1^37)& B_2,(2) & 1 & (12,-)+(3,-) \\& & & 2&(2,1)\vspace{0.1cm} \\
(12,-) & (37) & B_3,(12) & 1 & (12,-)\vspace{0.1cm}
\\ (3,-)&(19) & B_3,(3) & 0 & (3,-)
\end{array}\]

Notice that the decomposition into graded parts of the regular
representation (which we obtain for center 0 of the subspace $L=\real$ 
which has $R_L=\emptyset$), is the one of the coinvariant algebra with a degree shift, as seen in section \ref{sign}.


We can now check that indeed Conjecture \ref{conj} holds. To do this one needs to calculate explicitly $\Hgrrcc$. This is not very difficult but somewhat tedious, for the results we refer the reader to \cite{proefschrift}.

\begin{figure}[t]
\begin{center}
{\includegraphics[angle=0,scale=0.27]{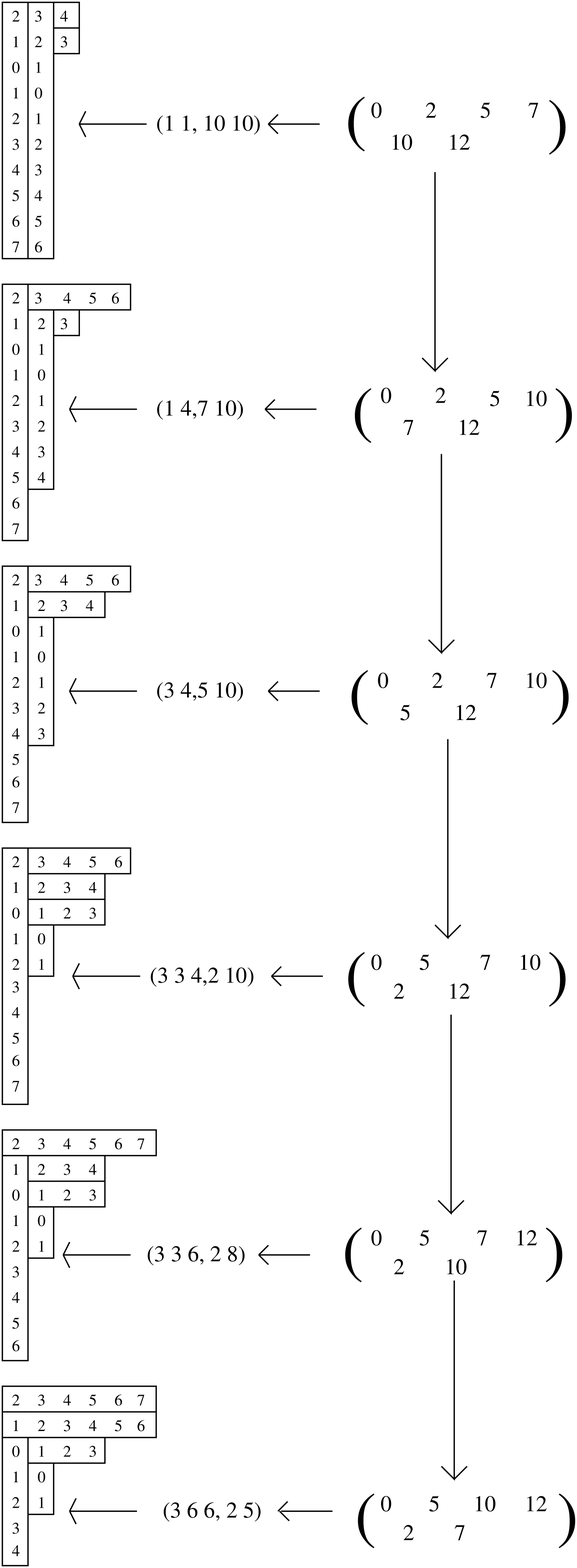}}
\end{center}
\label{voorbeeld}
\end{figure}

\begin{small}
\bibliographystyle{amsplain}
\bibliography{biblio}
\end{small}
\end{document}